\documentclass[oneside]{amsbook}
\def\VDB{}
\usepackage{amssymb}
\usepackage{psfrag}
\usepackage{verbatim}
\usepackage{graphicx}
\usepackage{amscd}
\usepackage[all]{xy}
\numberwithin{equation}{section}
 \numberwithin{figure}{section}
\numberwithin{section}{chapter}
\CompileMatrices
\newcommand{\h}{\ensuremath{\mathcal{H}}}

\newcommand{\X}{\ensuremath{\mathcal{X}}}

\newcommand{\Z}{\ensuremath{\mathbb{Z}}}

\newcommand{\la}{\ensuremath{\Lambda}}

\newcommand{\mc}{\ensuremath{\mathcal{C}}}
\newcommand{\ra}{\ensuremath{\rightarrow}}
\newcommand{\map}{\operatorname{map}}

\DeclareMathOperator{\m}{mod} \DeclareMathOperator{\Kdim}{Kdim}
\DeclareMathOperator{\mo}{Mod}
\DeclareMathOperator{\gd}{gl.dim}
 \DeclareMathOperator{\cok}{Cok}
\DeclareMathOperator{\HEIGHT}{ht}
\DeclareMathOperator{\height}{ht} \DeclareMathOperator{\op}{opp}
 
\DeclareMathOperator{\QMod}{QMod}
\DeclareMathOperator{\qmod}{qmod}
\def\dirlim{\mathop{\vtop{\baselineskip -100pt\lineskip 0.4ex
\setbox0\hbox{\upshape lim}\copy0\hbox to \wd0{\rightarrowfill}}}\limits}
\def\invlim{\mathop{\vtop{\baselineskip -100pt\lineskip 0.4ex
\setbox0\hbox{\upshape lim}\copy0\hbox to \wd0{\leftarrowfill}}}\limits}

\let\projlim\invlim

\let\dlim\dirlim

\let\cal\mathcal
\def\Ascr{{\cal A}}
\def\Bscr{{\cal B}}
\def\Cscr{{\cal C}}
\def\Dscr{{\cal D}}
\def\Escr{{\cal E}}
\def\Fscr{{\cal F}}

\def\Hscr{{\cal H}}
\def\Iscr{{\cal I}}

\def\Mscr{{\cal M}}
\def\Nscr{{\cal N}}
\def\Oscr{{\cal O}}
\def\Pscr{{\cal P}}
\def\Qscr{{\cal Q}}
\def\Rscr{{\cal R}}

\def\Tscr{{\cal T}}
\def\Uscr{{\cal U}}

\def\Xscr{{\cal X}}

\def\Zscr{{\cal Z}}
\let\blb\mathbb

\def\XX{{\blb X}}

\def \PP{{\blb P}}

\def \ZZ{{\blb Z}}

\def \NN{{\blb N}}

\def\id{\text{id}}
\def\Id{\operatorname{id}}

\def\Mod{\operatorname{Mod}}
\def\mod{\operatorname{mod}}
\def\Gr{\operatorname{Gr}}

\def\QGr{\operatorname{QGr}}
\def\qgr{\operatorname{qgr}}

\def\PC{\operatorname{PC}}
\def\pc{\operatorname{pc}}
\def\gr{\operatorname{gr}}

\def\length{\mathop{\text{length}}}

\def\coh{\mathop{\text{\upshape{coh}}}}

\def\rad{\operatorname {rad}}
\def\gr{\operatorname {gr}}
\def\Spec{\operatorname {Spec}}
\def\Rep{\operatorname {Rep}}

\def\Ext{\operatorname {Ext}}
\def\Hom{\operatorname {Hom}}
\def\End{\operatorname {End}}
\def\RHom{\operatorname {RHom}}
\def\uHom{\operatorname {\mathcal{H}\mathit{om}}}

\def\im{\operatorname {im}}

\def\ker{\operatorname {ker}}

\def\End{\operatorname {End}}

\def\Soc{\operatorname {soc}}

\def\rep{\operatorname {rep}}
\def\wrep{\widetilde{\operatorname {rep}}}
\def\Rep{\operatorname {Rep}}

\def\gldim{\operatorname {gl\,dim}}

\def\r{\rightarrow}
\def\l{\leftarrow}

\DeclareMathOperator{\Proj}{Proj}

\DeclareMathOperator{\tors}{tors}

\DeclareMathOperator{\dis}{dis} \DeclareMathOperator{\qpc}{qpc}

\newtheorem{lemma}{Lemma}[section]
\newtheorem{proposition}[lemma]{Proposition}
\newtheorem{theorem}[lemma]{Theorem}
\newtheorem{corollary}[lemma]{Corollary}

\newtheorem{lemmas}{Lemma}[subsection]
\newtheorem{propositions}[lemmas]{Proposition}

\newtheorem*{sublemma}{Sublemma}  %\def\thesublemma{}

\theoremstyle{definition}

\newtheorem{example}[lemma]{Example}

{

\newtheorem{step}{Step}
\newtheorem{case}{Case}

}

\theoremstyle{remark}
\newtheorem{remark}[lemma]{Remark}

\newdimen\uboxsep \uboxsep=1ex
\def\uboxn#1{\vtop to 0pt{\hrule height 0pt depth 0pt\vskip\uboxsep
\hbox to 0pt{\hss #1\hss}\vss}}

\def\uboxs#1{\vbox to 0pt{\vss\hbox to 0pt{\hss #1\hss}
\vskip\uboxsep\hrule height 0pt depth 0pt}}

\frontmatter

\title[Noetherian hereditary abelian categories satisfying Serre duality]
{Noetherian hereditary abelian categories satisfying Serre duality}
\author{I. Reiten}
\address{Department of Mathematical Sciences \\ Norwegian
University of Science and Technology \\ 7491 Trondheim \\ Norway}
\email{idunr@math.ntnu.no}
\author{M. Van den Bergh}
\thanks{The second author is a senior researcher at the Fund for
  Scientific Research. The second author also wished to thank the Clay
  Mathematics Institute
  for material support during the period that this paper was written.}
\address{Department WNI, Limburgs
Universitair Centrum \\ Universitaire Campus, Building D
\\3590 Diepenbeek \\ Belgium}
\email{vdbergh@luc.ac.be}


\begin{document}

\maketitle

\begin{abstract}
  In this paper  we classify $\Ext$-finite noetherian hereditary abelian categories over an
  algebraically closed field $k$
  satisfying Serre duality in the sense of Bondal and Kapranov. As
  a consequence we obtain a classification of saturated noetherian
  hereditary abelian categories.

As a side result we show that when our hereditary abelian categories have no nonzero
projectives or injectives, then the Serre duality property is equivalent to the existence of
almost split sequences.
\end{abstract}
%\maketitle
\tableofcontents

\chapter*{Notations and conventions}
  Most notations will be introduced locally. The few global ones are
  given below.  Unless otherwise specified $k$ will be an algebraically
  closed field and all rings and categories in this paper will be
  $k$-linear.

If $A$ is a ring then $\mod(A)$ will be the category of finitely generated right \VDB
$A$-modules. Similarly if $R$ is a $\ZZ$-graded ring then $\gr(R)$ will be the category of
finitely generated graded right modules with degree zero morphisms, and $\Gr (R)$ the category
of all graded right $R$-modules. If $R$ is noetherian
then following \cite{AZ} $\tors(R)$ will be the full subcategory
of $\gr(R)$ consisting of graded modules
with right bounded grading.  Also following \cite{AZ}
we put $\qgr(R)=\gr(R)/\tors(R)$.

For an abelian category \mc\ we denote by $D^b(\mc)$ the bounded derived category of \mc.

\chapter*{Introduction} One of the goals of non-commutative
algebraic geometry is to obtain an understanding of $k$-linear
abelian categories $\Cscr$, for a field $k$, which have properties
close to those of the category of coherent sheaves over a non-singular
proper scheme. Hence some obvious properties one may impose on
$\Cscr$ in this context are the following:
\begin{itemize}
\item $\Cscr$ is \emph{$\Ext$-finite}, i.e.\ $\dim_k \Ext^i(A,B)<\infty$
  for all $A,B\in\Cscr$ and for all $i$.
\item $\Cscr$ has \emph{homological dimension}
  $n<\infty$, i.e. $\Ext^i(A,B)=0$ for $A,B\in\Cscr$ and $i>n$, and $n$ is
  minimal with this property.
\end{itemize}
Throughout this paper $k$ will denote a field, and even though it is not always necessary we
will for simplicity assume that $k$ is algebraically closed. All categories will be
$k$-linear. When we say that $\Cscr$ is $\Ext$-finite, it will be understood that this is with
respect to the field $k$.

In most of this paper we will assume that $\Cscr$ is an $\Ext$-finite abelian category of
homological dimension at most 1, in which case we say that $\Cscr$ is
\emph{hereditary}.

A slightly more subtle property of non-singular proper schemes is Serre
duality. Let $X$ be a non-singular proper scheme over $k$ of dimension
$n$ and let $\coh(X)$ denote the category of coherent
$\Oscr_X$-modules. Then the classical Serre duality theorem
asserts that for $\Fscr\in\coh(X)$ there are natural isomorphisms
\[
H^i(X,\Fscr)\cong \Ext^{n-i}(\Fscr,\omega_X)^\ast
\]
 where $(-)^\ast = \Hom_k (-,k)$.

A very elegant reformulation of Serre duality was given by Bondal
and Kapranov in \cite{Bondal4}. It says that for any
$\Escr,\Fscr\in D^b(\coh(X))$ there exist natural isomorphisms
\[
\Hom_{D^b(\coh(X))}(\Escr,\Fscr)\cong
\Hom_{D^b(\coh(X))}(\Fscr,\Escr\otimes\omega_X[n])^\ast
\]
Stated in this way the concept of Serre duality can be generalized to certain abelian
categories.
\begin{itemize}
\item $\Cscr$ satisfies \emph{Serre duality} if it has a so-called
  \emph{Serre functor}. The latter is by definition an autoequivalence
  $F:D^b(\Cscr)\r D^b(\Cscr)$ such that there are isomorphisms
\[
\Hom(A,B)\cong \Hom(B,FA)^\ast
\]
which are natural in $A,B$.
\end{itemize}
On the other hand hereditary abelian $\Ext$-finite categories $\Cscr$ with the additional
property of having a tilting object  have been important for the representation theory of
finite dimensional algebras. Recall that $T$ is a tilting object in $\Cscr$ if
$\Ext^1(T,T)=0$, and if $\Hom (T,X)=0=\Ext^1(T,X)$ implies that $X$ is 0. These categories
$\Cscr$ are important in the study of quasitilted algebras, which by definition are the
algebras of the form $\End_{\Cscr} (T)$ for a tilting object $T$ \cite{HRS}, and which contain
the important classes of tilted and canonical algebras. A prominent property in the
representation theory of finite dimensional algebras is having almost split sequences, and
also the $\Ext$-finite hereditary abelian categories with tilting object have this property
\cite{HRS}.

In view of the above it is interesting, and useful, to investigate the relationship between
Serre duality and almost split sequences. In fact, this relationship is very close in the
hereditary case. The more general connections are on the level of triangulated categories,
replacing almost split sequences with Auslander-Reiten triangles. In fact one of our first
results in this paper is the following (see Chapter \ref{ref:I-0} for more complete results):
{
\def\thelemma{A}
\begin{theorem}
\begin{enumerate}
\item
$\Cscr$ has a Serre functor if and only if
  $D^b(\Cscr)$ has Auslander-Reiten triangles (as defined in \cite{Happel}).
\item If $\Cscr$ is hereditary, then $\Cscr$ has a Serre functor  if and only if
$\Cscr$ has almost split sequences and there is a one-one correspondence between the
indecomposable projective objects $P$ and the indecomposable injective objects $I$, such that
the simple top of $P$ is isomorphic to the socle of $I$.
\end{enumerate}
\end{theorem}
} Hence $\Ext$-finite hereditary abelian categories with Serre duality are of interest both
for non-commutative algebraic geometry and for the representation theory of finite dimensional
algebras. The main result of this paper is the classification of the  noetherian ones.

To be able to state our result we first give a list of hereditary abelian categories
satisfying Serre duality.
\begin{itemize}
\item[(a)] If $\Cscr$ consists of the finite dimensional nilpotent
  representations of the quiver $\tilde{A}_n$ or of the quiver $A_{\infty}^{\infty}$,
   with all arrows oriented in the same direction, then it
  is classical that $\Cscr$ has almost split sequences, and hence Serre
  duality.
\item[(b)] Let $X$ be a non-singular projective connected curve over $k$
  with function field $K$, and let $\Oscr$ be a sheaf of hereditary
  $\Oscr_X$-orders in $M_n(K)$ (see \cite{reiner}). Then one proves
  exactly as in the commutative case that $\coh(\Oscr)$ satisfies
  Serre duality.
\item[(c)] Let $Q$ be either $A_\infty^\infty$ or $D_\infty$ with
  zig-zag orientation. It is shown in \S\ref{ref:III.3-73} that
  there exists a noetherian hereditary abelian category $\Cscr$ which is
  derived equivalent to the category $\rep(Q)$ of finitely presented
  representations of $Q$, and which has no nonzero projectives or
  injectives. Depending on $Q$ we call this category the $\ZZ
  A^\infty_\infty$ or the $\ZZ D_\infty$ category.  Since Serre duality is
  defined in terms of the derived category, it follows that $\Cscr$
  satisfies Serre duality.

  If $Q=A_\infty^\infty$ then $\Cscr$ is nothing but the category
  $\gr_{\ZZ^2}(k[x,y])/(\text{finite length})$ considered in
  \cite{SmithZhang}. If $Q=D_\infty$ and $\text{char} k\neq 2$ then
  $\Cscr$ is a skew version of the $\ZZ A_\infty^\infty$ category (see
  \S\ref{ref:III.3-73}). The $\ZZ
  A_\infty^\infty$ category and the $\ZZ D_\infty$ category have
  also been considered by Lenzing.
\item[(d)] We now come to more subtle examples (see Chapter \ref{ref:II-23}).
  Let $Q$ be a connected quiver. Then for a vertex $x\in Q$ we have a
  corresponding projective representation $P_x$ and an injective
  representation $I_x$.  If $Q$ is locally finite and there is no infinite
  path ending at any vertex, the
  functor $P_x\mapsto I_x$ may be derived to yield a fully faithful
  endofunctor $F:D^b(\rep(Q))\r D^b(\rep(Q))$. Then $F$ behaves like a
  Serre functor, except that it is not in general an autoequivalence. We
  call such $F$ a right Serre functor (see \S\ref{ref:I.1-1}).
  Luckily given a right Serre functor there is a formal procedure to
  invert it so as to obtain a true Serre functor (Theorem
  \ref{ref:II.1.3-29}). This yields a hereditary abelian category
  $\wrep(Q)$ which satisfies Serre duality. Under the additional
  hypotheses that $Q$ consists of a subquiver $Q_o$ with no path of
  infinite length, with rays attached to vertices of $Q_o$, then  $\wrep(Q)$ turns
  out to be noetherian (see Theorem \ref{ref:II.4.3-56}). Here we mean by a ray an $A_{\infty}$
  quiver with no vertex which is a
  sink. An interesting feature of the noetherian categories  $\wrep(Q)$, exhibiting a new type
  of behavior, is that they are generated by the preprojective objects, but not necessarily by
  the projective objects.

\end{itemize}
Now we can state our main result. Recall that an abelian category
$\Cscr$ is connected if it cannot be non-trivially written as a
direct sum $\Cscr_1\oplus\Cscr_2$. {\def\thelemma{B}
\begin{theorem}
\label{theoremb} Let $\Cscr$ be a connected noetherian $\Ext$-finite
  hereditary abelian category satisfying Serre duality. Then  $\Cscr$ is one
  of the categories described in (a)-(d) above.
\end{theorem}
}The cases (a)(b)(c) are those where there are no nonzero projective objects. Those in (a) are
exactly the $\Cscr$ where all objects have finite length. For the $\Cscr$ having some objects
of infinite length, then either  all indecomposable objects of finite length have finite
$\tau$-period (case (b)) or all have infinite $\tau$-period (case (c)). Here the object $\tau
C$ is defined by the almost split sequence $0\ra \tau C \ra B \ra C\ra 0$ for $C$
indecomposable in $\Cscr$, and we have $\tau C=F(C)[-1]$.

Under the additional assumption that $\Cscr$ has a tilting object, such a classification was
given in \cite{Lenzing}. The only cases are the categories of finitely generated modules $\mod
A$ for a finite dimensional indecomposable hereditary $k$-algebra $A$ and the categories $\coh
\XX$ of coherent sheaves on a weighted projective line $\XX$ in the sense of \cite{GeigLen}.
{}{}From the point of view of the above list of examples the first case corresponds to the finite
quivers $Q$ in (d), in which case $\wrep Q =\rep Q$ is equivalent to $\mod (kQ)$, where $kQ$
is the path algebra of $Q$ over $k$. The categories $\coh \XX$ are a special case of (b),
corresponding to the case where the projective curve is $\PP^1$ (see \cite{IRMB}).

Our proof of Theorem \ref{theoremb} is rather involved and covers the first 4 chapters. The
main steps are as follows:
\begin{enumerate}
\item In the first two sections of Chapter \ref{ref:II-23} we construct
  the categories $\wrep(Q)$, and we show that they are characterized
  by the property of having noetherian injectives and being generated
  by preprojectives.
\item In \S\ref{ref:II.4-53} we give necessary and sufficient
  conditions for $\wrep(Q)$ to be noetherian, and furthermore we prove
  a decomposition theorem which states that an $\Ext$-finite  noetherian hereditary
  abelian category with Serre functor can
  be decomposed as a direct sum of a hereditary abelian category which is
  generated by preprojectives and a hereditary abelian category which doesn't
  have nonzero projectives or injectives.
\item We are now reduced to the case where there are no nonzero projectives or
  injectives. The case where all objects have finite length is
  treated in \S\ref{ref:III.1-64}.
\item The case where there are no nonzero projectives or injectives and at least one object of infinite length
  is covered in Chapter \ref{ref:IV-107}. It turns out that this
  case naturally falls  into two subcases:
\begin{itemize}
\item[($\alpha$)] The simple objects are $\tau$-periodic. In that case, using the results in \cite{AZ}, we show that $\Cscr$
  is of the form $\qgr(R)$ for $R$ a two-dimensional commutative graded ring,
  where $\qgr (R)$ is the quotient category $\gr  R$/ finite length.
  Using
  \cite{Staf5} it then follows that $\Cscr$ is of the form
  (b).
\item[($\beta$)] The simple objects are not $\tau$-periodic. We
  show that if such $\Cscr$ exists, then it is characterized by the fact
  that it has either one or two $\tau$-orbits of simple objects. Since the
  $\ZZ A_\infty^\infty$ and $\ZZ D_\infty$
  category have this property, we are done.
\end{itemize}
\end{enumerate}
Our methods for constructing the new hereditary abelian categories $\wrep(Q)$  are somewhat
indirect although we believe  they are interesting. After learning about our results Claus
Ringel has recently found a more direct construction for these categories \cite{Ri1}.

All the hypotheses for Theorem \ref{theoremb} are necessary. For example the non-commutative
curves considered in \cite{SmithZhang} are noetherian hereditary abelian categories of
Krull-dimension one which in general do not satisfy Serre duality (except for the special case
listed in (c)). If $\Cscr$ is the opposite category to one of the categories (b)(c)(d), then
it is not  noetherian, but it  satisfies the other hypotheses.
%Other examples of non-noetherian hereditary
%categories with Serre duality which have recently become of
%interest are the so-called non-commutative projective spaces
%constructed by Kontsevich and Rosenberg \cite{KoRo}.

Nevertheless it is tempting to ask whether a result similar to Theorem \ref{theoremb} remains
valid without the noetherian hypothesis if we work up to derived equivalence. In particular,
is any such category derived equivalent to a noetherian one? Under the additional assumption
that $\Cscr$ has a tilting object, this has been proved by Happel in \cite{HappelReiten}, and
it has recently been shown by Ringel that this is not true in general \cite{Ri2}.

\bigskip

In the final chapter we use Theorem \ref{theoremb} to draw some conclusions on the structure
of certain hereditary abelian categories.

To start with we discuss the ``saturation'' property. This is a subtle property of certain
abelian categories which was discovered by Bondal and Kapranov \cite{Bondal4}.  Recall that a
cohomological functor $H:D^b(\Cscr)\r \mod(k)$ is of finite type if for every $A\in
D^b(\Cscr)$ only a finite number of $H(A[n])$ are non-zero. We have already defined what it
means for $\Cscr$ to have homological dimension $n<\infty$. It will be convenient to say more
generally that $\Cscr$ has finite homological dimension if for any $A,B$ in $\Cscr$ there is
at most a finite number of $i$ with $\Ext^i(A,B)\neq0$. In particular, the
analogue of this definition makes sense for triangulated categories.
\begin{itemize}
\item Let $\Cscr$ be an $\Ext$-finite abelian category of finite
  homological dimension.  Then  $\Cscr$ is \emph{saturated} if  every
  cohomological functor $H:D^b(\Cscr)\r \mod(k)$ of finite type is of
  the form $\Hom(A,-)$ (i.e. $H$ is representable).
\end{itemize}
It is easy to show that a saturated category satisfies Serre duality. It was shown in
\cite{Bondal4} that $\coh(X)$ for $X$ a non-singular projective scheme is saturated and that
saturation also holds for categories of the form $\mod(\Lambda)$ with $\Lambda$ a finite
dimensional algebra. Inspired by these results we prove the following result in
\S\ref{ref:V.1-154}. {\def\thelemma{C}
\begin{theorem}
\label{theoremc}
 Assume that $\Cscr$ is a saturated connected
  noetherian $\Ext$-finite hereditary abelian category. Then
  $\Cscr$ has one of the following forms:
\begin{enumerate}
\item $\mod(\Lambda)$ where $\Lambda$ is a connected finite dimensional
  hereditary $k$-algebra.
\item $\coh(\Oscr)$ where $\Oscr$ is a sheaf of hereditary $\Oscr_X$-orders
 (see (b) above)   over a non-singular connected
  projective curve $X$.
\end{enumerate}
\end{theorem}
} It is easy to see that the hereditary abelian categories listed in the above theorem are of
the form $\qgr(R)$. We refer the reader to \cite{BondalVdb} (see also \S\ref{ref:V.2-157})
where it is shown in reasonable generality that abelian categories of the form $\qgr(R)$ are
saturated.
\begin{comment}
Next we discuss in \S\ref{ref:V.2-157} which of
our categories are of the form $\qgr(R)$. Our main result is the
following: {\def\thelemma{D}
\begin{theorem}
Assume that $\Cscr$ is a connected noetherian $\Ext$-finite
  hereditary abelian category with Serre functor. Assume furthermore that
  $\Cscr$ is of the form $\qgr(R)$ where $R$ is a noetherian
  $\NN$-graded ring finite dimensional in every degree which satisfies
  $\chi$ (see \cite{AZ} and \S\ref{ref:V.2-157}). Then $\Cscr$ has one of the following forms.
\begin{enumerate}
\item $\mod(\Lambda)$ where $\Lambda$ is a connected finite dimensional
  hereditary algebra.
\item $\coh(\Oscr)$ where $\Oscr$ is a sheaf of hereditary
$\Oscr_X$-orders (see \S\ref{ref:III.2-68})
over a non-singular connected
  projective curve $X$.
\end{enumerate}
\end{theorem}
} In view of Theorems C,D the following question appears to be
very natural. {\def\thelemma{F}
\bigskip

{\sc question.} \label{questionquestion} Is it possible to prove
in reasonable generality that
  abelian categories of the form $\qgr(R)$ are saturated?
\bigskip

} As remarked above saturatedness implies Serre duality. Theorem
\ref{gradedserredualitytheorem}, which is
proved using \cite{YZ}, shows that under suitable hypotheses
$\qgr(R)$ will have Serre duality. This gives some partial
evidence for a positive reply to the question.
\end{comment}

There are also applications to the relationship between existence of tilting objects and the
Grothendieck group being finitely generated. This was one of the original motivations for this
work, and is dealt with in another paper \cite{IRMB}.

We would like to thank Claus Ringel for helpful comments on the presentation of this paper.

\mainmatter

\chapter{Serre duality and almost split sequences}
\label{ref:I-0} It has been known for some time that there is a connection between classical
Serre duality and existence of almost split sequences.  There is a strong analogy between the
Serre duality formula for curves and the formula $D (\Ext ^1_\la (C,A)) \simeq \overline{\Hom}
(A,\operatorname{DTr}C)$ for artin algebras (where $D=\Hom_k (-,k)$), on which the existence of
almost split sequences is based (see \cite{ARS}). Actually, existence of almost split
sequences in some sheaf categories for curves can be proved either by using an analogous
formula for graded maximal Cohen--Macaulay modules or by using Serre duality
\cite{MAIR}\cite{ASch}. The notion of almost split sequences was extended to the notion of
Auslander--Reiten triangles in triangulated categories \cite{Happel}, and existence of such
was proved for $D^b(\mod \la)$ when $\la$ is a $k$-algebra of finite global dimension
\cite{Happel}. In this case the corresponding translate is given by an equivalence of
categories.  On the other hand an elegant formulation of Serre duality in the bounded derived
category, together with a corresponding Serre functor, was given in \cite{Bondal4}. These
developments provide the basis for further connections, which turn out to be most complete in
the setting of triangulated categories. For abelian categories we obtain strong connections in
the hereditary case. In fact, we show that when $\Cscr$ is hereditary then $\Cscr$ has Serre
duality if and only if it has almost split sequences and there is a one-one correspondence
between indecomposable projective objects $P$ and indecomposable injective objects I, such
that  $P$ modulo its unique maximal subobject is isomorphic to the socle of I.

\section{Preliminaries on Serre duality}
\label{ref:I.1-1}
 Let $\Ascr$ be a $k$-linear
$\Hom$-finite additive category. A \emph{{{right Serre}} functor}
is an additive functor $F:\Ascr\r \Ascr$ together with  isomorphisms
\begin{equation}
\label{ref:I.1.1-2} \eta_{A,B}:\Hom(A,B)\r \Hom(B,FA)^\ast
\end{equation}
for any $A,B\in\Ascr$ which are natural in $A$ and $B$. A left Serre functor is a functor
$G:\Ascr\r\Ascr$ together with  isomorphisms
\begin{equation}
\zeta_{A,B}:\Hom(A,B)\r \Hom(GB,A)^\ast
\end{equation}
for any $A,B\in\Ascr$ which are natural in $A$ and $B$. Below we state and prove a number of
properties of right Serre functors. We  leave the proofs of the corresponding properties for
left Serre functors to the reader.

Let $\eta_A:\Hom(A,FA)\r k$ be given by $\eta_{A,A}(\Id_A)$ and
let $f\in\Hom(A, B)$. Looking at the following commutative diagram
(which follows from the naturality of $\eta_{A,B}$ in $B$)
\[
\begin{CD}
\Hom(A,A) @>\eta_{A,A} >> \Hom(A,FA)^\ast\\ @V\Hom(A,f)VV @V
\Hom(f,FA)^\ast VV\\ \Hom(A,B) @>\eta_{A,B}>> \Hom(B,FA)^\ast
\end{CD}
\]
we find for $g\in\Hom(B,FA)$
\begin{equation}
\label{ref:I.1.3-3} \eta_{A,B}(f)(g)=\eta_A  (g f)
\end{equation}
 Similarly by the naturality of $\eta_{A,B}$ in $A$ we obtain a commutative diagram
\[
\begin{CD}
\Hom(B,B) @>\eta_{B,B}>> \Hom(B,FB)^\ast\\ @V \Hom(f,B) VV
@V\Hom(B,Ff)^\ast VV\\ \Hom(A,B) @> \eta_{A,B}>> \Hom(B,FA)^\ast
\end{CD}
\]
This yields for $g\in \Hom(B,FA)$ the formula
\begin{equation}
\label{ref:I.1.4-4} \eta_{A,B}(f)(g)=\eta_B(F(f) g)
\end{equation}
and we get the following description of the functor $F$.
\begin{lemma}
The following composition coincides with $F$.
\[
\Hom(A,B)\xrightarrow{\eta_{A,B}} \Hom(B,FA)^\ast
\xrightarrow{(\eta^\ast_{B,FA})^{-1}} \Hom(FA,FB)
\]
\end{lemma}
\begin{proof}
To prove this we need to show that for $f\in \Hom(A,B)$ and $g\in \Hom(B,FA)$ one has
$\eta_{A,B}(f)(g)= \eta^{\ast}_{B,FA}(Ff)(g)$. Thanks to the formulas
\eqref{ref:I.1.3-3}\eqref{ref:I.1.4-4} we obtain $\eta_{A,B}(f)(g)=\eta_B(Ff \circ g)$ and
also $\eta^{\ast}_{B,FA}(Ff)(g) =\eta_{B,FA}(g)(Ff)=\eta_B(F(f) g)$. Thus we obtain indeed the
correct result.
\end{proof}
We have the following immediate consequence.
\begin{corollary}
If $F$ is a {{right Serre}} functor, then $F$ is fully faithful.
\end{corollary}
Also note the following basic properties.
\begin{lemma}
\begin{enumerate}
\item If $F$ and $F'$ are {{right Serre}} functors, then they are naturally isomorphic.
\item $\Ascr$ has a {{right Serre}} functor if and only if
  $\Hom(A,-)^\ast$ is representable for all $A\in\Ascr$.
\end{enumerate}
\end{lemma}
{}{}From the above discussion it follows that there is a lot of redundancy in the data
$(F,(\eta_{A,B})_{A,B})$. In fact we have the following.
\begin{proposition}
\label{ref:I.1.4-5} In order to give $(F,(\eta_{A,B})_{A,B})$
it is necessary and sufficient
  to give the action of $F$ on objects, as well as $k$-linear maps
  $\eta_{A}:\Hom(A,FA)\r k$ such that the composition
\begin{equation}
\label{ref:I.1.5-6} \Hom(A,B)\times \Hom(B,FA)\r
\Hom(A,FA)\xrightarrow{\eta_A} k
\end{equation}
yields a non-degenerate pairing for all $A,B\in\Ascr$. If we are
given $\eta_A$, then $\eta_{A,B}$ is obtained from the pairing
\eqref{ref:I.1.5-6}. Furthermore the action $F$ on maps
\[
F:\Hom(A,B)\r \Hom(FA,FB)
\]
is defined by the property that for $f\in\Hom(A,B)$ we have
$\eta_A(gf)=\eta_B(F(f)g)$ for all $g\in \Hom(B,FA)$.
\end{proposition}
\begin{proof} It is clear from the previous discussion that the data
  $(F,(\eta_{A,B})_{A,B})$  gives rise to $(\eta_A)_A$ with the required
  properties. So conversely assume that we are given $(\eta_A)_A$
  and the action of $F$ on objects. We define $(\eta_{A,B})_{A,B}$ and
  the action of $F$ on maps as in the statement of the proposition.

\noindent

  We first show that $F$ is a functor. Indeed let $A,B,C\in\Ascr$
   and assume that there are maps $g:A\r B$ and $h:B\r C$.
  Then for all $f\in \Hom(C,FA)$ we have $ \eta_A(f hg)=\eta_C(F(hg)
  f)$, but also $\eta_A(fhg)= \eta_B(F(g) fh)=\eta_C(F(h) F(g) f)$.
  Thus by non-degeneracy we have $F(hg)=F(h)F(g)$.

\noindent

It is easy to see that the pairing \eqref{ref:I.1.5-6} defines an isomorphism
\[
\eta_{A,B}:\Hom(A,B)\r\Hom(B,FA)^\ast: f\mapsto \eta_A(-f)
\]
which is natural in $A$ and $B$.

The proof is now complete.
\end{proof}
A Serre functor  is by definition a {{right Serre}} functor which
is essentially surjective. The following is easy to see.
\begin{lemma} $\Ascr$ has a Serre functor if and only it has both a
  right and a left Serre functor.
\end{lemma}
{}{}From this we deduce the following  \cite{Bondal4}.
\begin{lemma}
$\Ascr$ has a Serre functor if and only if the functors
  $\Hom(A,-)^\ast$ and $\Hom(-,A)^\ast$ are representable for all
  $A\in\Ascr$.
\end{lemma}
\begin{remark}
\label{ref:I.1.7-7} In the sequel $\Ascr$ will always be a Krull--Schmidt category (in the
sense that indecomposable objects have local endomorphism rings). In that case it is clearly
sufficient to specify $\eta_{A,B},\eta_A,F$, etc\dots on the full subcategory of $\Ascr$
consisting of indecomposable objects.
\end{remark}
If one is given a {{right Serre}} functor, then it is possible to
invert it formally in such a way that the resulting additive
category has a Serre functor. The next result is stated in
somewhat greater generality.
\begin{proposition}
\label{ref:I.1.8-8} Let $\Ascr$ be an additive category as
above and let $U:\Ascr\r \Ascr$ be a fully faithful additive
endofunctor.

Then there exists an additive category $\Bscr$
 with the following properties.
\begin{enumerate}
\item There is a fully faithful functor $i:\Ascr\r \Bscr$.
\item There is an  autoequivalence ${{\bar{U}}}:\Bscr\r \Bscr$
  together with a natural isomorphism $\nu:\bar{U}i\r iU$.
\item For every object $B\in\Bscr$ there is some $b\in\NN$ such that
  $\bar{U}^bB$ is isomorphic to $i(A)$ with $A\in \Ascr$.
\end{enumerate}
Furthermore a quadruple $(\Bscr,i,{{\bar{U}}},\nu)$ with these
properties is unique (in the appropriate sense).
\end{proposition}
\begin{proof} Let us sketch the construction of $\Bscr$. The
  uniqueness will be clear.

  The objects in $\Bscr$ are formally written as $U^{-a}A$ with
  $A\in\Ascr$ and $a\in \ZZ$. A morphism $U^{-a}A\r U^{-b}B$ is
  formally written as $U^{-c}f$ with $f\in
  \Hom_{\Ascr}(U^{c-a}A,U^{c-b}B)$ where $c\in\ZZ$ is such that $c\ge a$,
  $c\ge b$. We identify
  $U^{-c}f$ with $U^{-c-1}(Uf)$.

The functor $i$ is defined by $i(A)=U^0A$ and the functor
${{\bar{U}}}$ is defined by  ${{\bar{U}}}(U^{-a})(A)=U^{-a+1}(A)$.
It is clear that these have the required properties.
\end{proof}
The following lemma provides a complement to this proposition in
the case that $\Ascr$ is triangulated.
\begin{lemma}
Assume that in addition to the usual hypotheses one has that
$\Ascr$ is triangulated. Let $(\Bscr,i,{{\bar{U}}},\nu)$ be as in
the previous proposition.
 Then there is a unique way to
make $\Bscr$ into a triangulated category such that $i$ and ${{\bar{U}}}$ are exact.
\end{lemma}
\begin{proof}
If we require exactness of $i$ and ${{\bar{U}}}$, then there is only one way to make $\Bscr$
into a triangulated category. First we must define the shift functor by
$(U^{-a}A)[1]=U^{-a}(A[1])$ and then the  triangles in $\Bscr$ must be those diagrams that are
isomorphic to
\[
U^{-c}A\xrightarrow{U^{-c}f} U^{-c}B\xrightarrow{U^{-c}g} U^{-c}C
\xrightarrow{U^{-c}h}U^{-c}A[1]
\]
where
\[
A\xrightarrow{f}B\xrightarrow{g} C\xrightarrow{h} A[1]
\]
is a triangle in  $\Ascr$ (note that the exactness of $\bar{U}$ is
equivalent to that of $\bar{U}^{-1}$).

To show that this yields indeed a triangulated category one must
check the axioms in \cite{Verdier}. These all involve the
existence of certain objects/maps/triangles. By applying a
sufficiently high power of ${{\bar{U}}}$ we can translate such
problems into ones involving only objects in $\Ascr$. Then we  use
the triangulated structure of $\Ascr$ and afterwards we go back to
the original problem by applying a negative power of
${{\bar{U}}}$.
\end{proof}
In the sequel we will denote by $U^{-\infty} \Ascr$ the category
$\Bscr$ which was constructed in Proposition
\ref{ref:I.1.8-8}. Furthermore we will consider $\Ascr$ as
subcategory of $U^{-\infty}\Ascr$ through the functor $i$. Finally
we will usually write $U$ for the extended functor ${{\bar{U}}}$.

Below we will only be interested in the special case where $U=F$ is a {{right Serre}} functor
on $\Ascr$. In that case we have the following.
\begin{proposition}
\label{ref:I.1.10-9} The canonical extension of  $F$ to
$F^{-\infty}\Ascr$ is
  a Serre functor.
\end{proposition}
\begin{proof} By construction $F$ is an automorphism on $F^{-\infty}\Ascr$. To prove that
  $F$ is a Serre functor
    we have to construct suitable maps $\eta_{F^{-a}A,F^{-b}B}$. Pick $c\ge a$,
    $c\ge b$. Then we have
\begin{align*}
  \Hom_{F^{-\infty}\Ascr}(F^{-a}A,F^{-b}B)&=
  \Hom_{\Ascr}(F^{c-a}A,F^{c-b}B)\\
  &\overset{\eta_{F^{c-a}A,F^{c-b}B}}{\cong}\Hom_{\Ascr}(F^{c-b}B,F^{c-a+1}A)^\ast\\
  &=\Hom_{F^{-\infty}\Ascr}(F^{-b}B,F^{-a+1}A)^\ast
\end{align*}
We define $\eta_{F^{-a}A,F^{-b}B}$ as the composition of these maps. It follows easily from
Lemma I.1.1 that the constructed map is independent of $c$, and it is clear that
$\eta_{F^{-a}A,F^{-b}B}$ has the required properties.
\end{proof}
We shall also need the following easily verified fact.
\begin{lemma}
\label{ref:I.1.11-10} If $\Ascr=\Ascr_1\oplus\Ascr_2$ is a
direct sum of additive categories, then  a (right) Serre functor
on $\Ascr$ restricts to (right) Serre functors on $\Ascr_1$ and
$\Ascr_2$.
\end{lemma}

\section{Connection between Serre duality and Auslander--Reiten
  triangles}
\label{ref:I.2-11} In this section we prove that existence of a right Serre functor is
equivalent to the existence of right Auslander--Reiten triangles, in triangulated
$\Hom$-finite Krull--Schmidt $k$-categories. Hence the existence of a Serre functor is
equivalent to the existence of Auslander--Reiten triangles.

In the sequel $\Ascr$ is a $\Hom$-finite $k$-linear Krull-Schmidt triangulated category.
Following \cite{Happel} a triangle $A\xrightarrow{f} B\xrightarrow{g} C\xrightarrow{h} A[1]$
in $\Ascr$ is called an \emph{Auslander--Reiten triangle} if the following conditions are
satisfied.
\begin{itemize}
\item[(AR1)] $A$ and $C$ are indecomposable.
\item[(AR2)] $h\neq 0$.
\item[(AR3)] If $D$ is indecomposable then for every
  non-isomorphism $t:D\r C$ we have $ht=0$.
\end{itemize}
It is shown in \cite{Happel} that, assuming (AR1)(AR2), then (AR3)
is equivalent to
\begin{itemize}
\item[(AR3)'] If $D'$ is indecomposable then for every non-isomorphism
  $s:A\r D'$ we have $sh[-1]=0$.
\end{itemize}
We say that  right Auslander--Reiten triangles exist in $\Ascr$ if for all indecomposables
$C\in\Ascr$ there is a triangle satisfying the conditions above. Existence of left
Auslander--Reiten triangles is defined in a similar way, and we say that $\Ascr$ has
Auslander--Reiten triangles if it has both right and left Auslander--Reiten triangles. (Note
that in \cite{Happel} one says that $\Ascr$ has Auslander--Reiten triangles if it has right
Auslander--Reiten triangles in our terminology.)

It is shown in \cite[\S4.3]{Happel} that given $C$ the corresponding Auslander--Reiten
triangle $A\ra B \ra C\ra A [1]$ is unique up to isomorphism of triangles. (By duality a
similar result holds if  $A$ is given.) For a given indecomposable $C$ we let $\tilde{\tau} C$
be an arbitrary object in $\Ascr$, isomorphic to $A$ in the Auslander--Reiten triangle
corresponding to $C$.

The following characterization of Auslander--Reiten triangles is analogous to the
corresponding result on almost split sequences (see \cite{ARS}).
\begin{proposition}
\label{ref:I.2.1-12} Assume that $\Ascr$ has right Auslander--Reiten
  triangles, and assume that we have a triangle in $\Ascr$
\begin{equation}
\label{ref:I.2.1-13} A\xrightarrow{f} B\xrightarrow{g}
C\xrightarrow{h} A[1]
\end{equation}
with $A$ and $C$ indecomposable and $h\neq 0$. Then the following are equivalent:
\begin{enumerate}
\item The triangle \eqref{ref:I.2.1-13} is an Auslander--Reiten triangle.
\item The map $h$ is in the socle of $\Hom(C,A[1])$ as right $\End(C)$-module
  and $A\cong \tilde{\tau} C$.
\item The map $h$ is in the socle of $\Hom(C,A[1])$ as left
  $\End(A)$-module and $A\cong \tilde{\tau} C$.
\end{enumerate}
\end{proposition}
\begin{proof} We will show that 1. and 2. are equivalent. The
  equivalence of 1. and 3. is similar.
\begin{itemize}
\item[$1.\Rightarrow 2.$]
By definition we have $A\cong \tilde{\tau} C$. Assume that $t\in\End(C)$ is a
  non-automorphism. Then by (AR3) we have $ht=0$.
\item[$2.\Rightarrow 1.$]
Let
\begin{equation}
\label{ref:I.2.2-14} A\xrightarrow{f'}B'\xrightarrow{g'} C
\xrightarrow{h'} A[1]
\end{equation}
be the Auslander--Reiten triangle associated to $C$. {}From the properties of Auslander--Reiten
triangles it follows that there is a morphism of triangles
\[
\begin{CD}
A @>f>> B @>g>> C @>h>> A[1]\\ @| @A sAA @A t AA @|\\ A @>f' >> B'
@>g'>> C @>h'>>A[1]
\end{CD}
\]
The fact that $h'\neq 0$ together with the fact that $h$ is in the socle of $\Hom(C,A[1])$
implies that $t$ must be an isomorphism. But then by the properties of triangles $s$ is also
an isomorphism. So in fact the triangles \eqref{ref:I.2.2-14} and \eqref{ref:I.2.1-13} are
isomorphic, and hence in particular \eqref{ref:I.2.1-13} is an Auslander--Reiten triangle.
\qed
\end{itemize}
\def\qed{}\end{proof}
\begin{corollary} Assume that $\Ascr$ has right Auslander--Reiten
  triangles. Then the socle of $\Hom(C,\tilde{\tau} C[1])$ is one dimensional
  both as right $\End(C)$-module and as left $\End(\tilde{\tau} C)$-module.
\end{corollary}
\begin{proof} It is easy to see that $u$ linearly independent
  elements of the (left
  or right) socle define different triangles. However Auslander--Reiten
  triangles are unique. This is a contradiction.
\end{proof}
The following is the basis for the  main result of this section.
\begin{proposition}
\label{ref:I.2.3-15} The following are equivalent.
\begin{enumerate}
\item $\Ascr$ has a {{right Serre}} functor.
\item $\Ascr$ has right Auslander-Reiten triangles.
\end{enumerate}
If either of these properties holds, then the action of the Serre functor on objects coincides
with $\tilde{\tau}[1]$.
\end{proposition}
\begin{proof} \strut

\begin{itemize}
\item[$1.\Rightarrow 2.$] Let $C\in \Ascr$ be an
  indecomposable object. By Serre  duality
  there is a natural isomorphism
\[
\Hom(C,FC)\cong \Hom(C,C)^\ast
\]
as $\End(C)$-bimodules. In particular $\Hom(C,FC)$ has a one dimensional socle which
corresponds to the map $\theta_C:\Hom(C,C)\r \Hom(C,C)/\rad(C,C)=k$.  Define $\tilde{\tau}
C=FC[-1]$ and let $h$ be a non-zero element of the socle of $\Hom(C,\tilde{\tau} C[1])$. We
claim that the associated triangle
\[
\tilde{\tau} C \r X \r C \xrightarrow{h} \tilde{\tau} C[1]
\]
is an  Auslander--Reiten triangle. Let $D$ be indecomposable and let $t:D\r C$ be a
non-isomorphism. We have to show that the composition
\[
D\xrightarrow{t} C\xrightarrow{h} FC
\]
is zero. Using Serre duality this amounts to showing that the composition
\[
\Hom(C,D)\xrightarrow{\Hom(C,t)} \Hom(C,C)\xrightarrow{\theta_C} k
\]
is zero. Since $t$ is a non-isomorphism, this is clear.
\item[$2.\Rightarrow 1.$] This is the interesting direction.
As pointed out in remark \ref{ref:I.1.7-7} it is sufficient to construct the Serre functor on
the full subcategory of $\Ascr$ consisting of the indecomposable objects. For $A$ an
indecomposable object in $\Ascr$ we let  $FA$ be the object $\tilde{\tau} A[1]$. Let $h_A$ be
a non-zero element of $\Hom(A,\tilde{\tau} A[1])$
 representing the Auslander--Reiten triangle
\[
\tilde{\tau} A\r X \r  A\xrightarrow{h_A} \tilde{\tau} A[1]
\]
\begin{sublemma} Let $A$ and $B$ be indecomposable objects in
  $\Cscr$. Then the following hold.
\begin{enumerate}
\item For any non-zero $f\in \Hom(B,\tilde{\tau} A[1])$ there exists $g\in
  \Hom(A,B)$ such that $fg=h_A$.
\item For any non-zero $g\in \Hom(A,B)$ there exists $f\in
  \Hom(B,\tilde{\tau} A[1])$ such that $fg=h_A$.
\end{enumerate}
\end{sublemma}
\begin{proof}
\begin{enumerate}
\item
Using the properties of Auslander-Reiten triangles there is a
morphism between the triangle determined by $h_A$ (the
AR-triangle) and the triangle determined by $f$.
\[
\begin{CD}
\tilde{\tau} A @>>> Y @>>>  A @>h_A>> \tilde{\tau} A[1]\\ @| @VVV @V gVV @|\\ \tilde{\tau} A
@>>> X @>>>  B @>f>> \tilde{\tau} A[1]\\
\end{CD}
\]
The morphism labeled $g$ in the above diagram has the required
properties.
\item Without loss of generality we may assume that $g$ is not an
  isomorphism. We complete $g$ to a triangle
\[
Z\xrightarrow{s} A\xrightarrow{g} B \xrightarrow{t} Z[1]
\]
Then $Z\neq 0$ and since $g$ is non-zero, $s$ will not be split.
Now we look at the following diagram
\[
\begin{CD}
@. @. B @.\\ @. @. @A gAA @.\\ \tilde{\tau} A @>>> Y @>>> A @>h_A>> \tilde{\tau} A[1]\\ @. @.
@AsAA @.\\ @. @. Z
\end{CD}
\]
Since $s$ is not split we have by (AR3) that $h_As=0$. Hence by the properties of triangles we
have $h_A=fg$ for a map $f:B\r \tau A[1]$. This proves what we want. \qed
\end{enumerate}
\def\qed{}\end{proof}
Having proved the sublemma we return to the main proof. For any
indecomposable object $A\in\Cscr$ choose a linear map
\[
\eta_A:\Hom(A,\tilde{\tau} A[1])\r k
\]
such that $\eta_A(h_A)\neq 0$. It follows from the sublemma that
the pairing
\begin{equation}
  \Hom(A,B)\times \Hom(B,\tilde{\tau} A[1])\r k\; \text{given by } (g,f)\r
  \eta_A(fg)\label{ref:I.2.3-16}
\end{equation}
is non-degenerate. We can now finish our proof by invoking
Proposition \ref{ref:I.1.4-5}.\qed
\end{itemize}
\def\qed{}\end{proof}
The following is now a direct consequence.
\begin{theorem} The following are equivalent:
\begin{enumerate}
\item
$\Ascr$ has  a Serre functor.
\item
$\Ascr$ has  Auslander--Reiten triangles.
\end{enumerate}
\end{theorem}
\begin{proof}
This follows from applying Theorem \ref{ref:I.2.3-15} together
with its dual version for left Serre functors.
\end{proof}
{}{}From now on we shall denote by $\tilde{\tau}$ also the equivalence $F[-1]$, where $F$ is the
Serre functor.

\section{Serre functors on hereditary abelian categories}
\label{ref:I.3-17} In this section we investigate the
relationship between existence of Serre functors and of almost
split sequences for hereditary abelian categories.

If $\Cscr$ is an $\Ext$-finite $k$-linear abelian category  we  say that it has a right Serre
functor  if this is the case for $D^b(\Cscr)$.

If $\Cscr$ has a right Serre functor $F$ and $A$ and $B$ are in $ \Cscr$, then from the fact
that $\Ext^i(A,B)=\Hom(A,B[i])=\Hom(B[i],FA)^\ast$ we deduce that only a finite number of
$\Ext^i(A,B)$ can be non-zero.

Before we go on we recall some basic definitions (see [4]). For an indecomposable object $C$ in $\Cscr$ a $\map g:B\to
C$ is right almost split if for any nonisomorphism $h:X\to C$ with $X$ indecomposable in $\Cscr$
there is some $\map t: B\to C$ with $gt=h$. The $\map g:B\to C$ is minimal right almost split
if in addition $g:B\to C$ is right minimal, that is, any $\map s:B\to B$ with $gs=g$ is an
isomorphism. The concepts of left almost split and minimal left almost split are defined
similarly. A nonsplit exact sequence $0\to A\to B\to C\to 0$ is almost split if $A$ and $C$
are indecomposable and $g:B\to C$ is (minimal)right almost split (or equivalently, $f:A\to B$
is (minimal) left almost split). We say that $\Cscr$ has right almost split sequences if for
every non-projective indecomposable object $C\in\Cscr$ there exists an almost split sequence
ending in $C$, and for each indecomposable projective object $P$ there is a minimal right
almost split map $E\ra P$. Possession of left almost split sequences is defined similarly. We
say that $\Cscr$ has almost split sequences if it has both left and right almost split
sequences.

Now let  $\Cscr$ be an $\Ext$-finite $k$-linear
  hereditary abelian  category.

The following characterization of when we have a minimal right almost split map to a
projective object or a minimal left almost split map from an injective object is easy to see.
\begin{lemma}
\label{ref:I.3.1.1ny}
\begin{enumerate}
\item There is some minimal right almost split map to an indecomposable projective object $P$
if and only if $P$ has a unique maximal subobject rad $P$. (If the conditions are satisfied,
then the inclusion map $i\colon \rad P\ra P$ is minimal right almost split).
\item There is a minimal left almost split map from an indecomposable injective object I if
and only if I has a unique simple subobject $S$. (If these conditions are satisfied, then the
epimorphism $j\colon I\ra I/ S$ is minimal left almost split.)
\end{enumerate}
\end{lemma}
Note that if an indecomposable projective object has a maximal subobject, then it must be
unique, and if $P$ is noetherian, then $P$ has a maximal subobject. Similarly if an
indecomposable injective object $I$ has a simple subobject, then it must be unique.

\begin{lemma}
\label{ref:I.3.1-18} Let $\Ascr=D^b(\Cscr)$. The following are
equivalent.
\begin{enumerate}
\item $\Ascr$ has right Auslander--Reiten triangles.
\item $\Cscr$  has right almost split sequences and for
every indecomposable
  projective $P\in\Cscr$ the simple object $S=P/\rad(P)$
  possesses an  injective hull $I$ in $\Cscr$.
\item $\Cscr$ has a right Serre functor.
\end{enumerate}
If any of these conditions holds, then the right Auslander--Reiten triangles in $\Ascr$ are
given by the shifts of the right almost split sequences in $\Cscr$ together with the shifts of
the triangles of the form
\begin{equation}
\label{ref:I.3.1-19} I[-1]\r  X  \r        P\xrightarrow{h_P} I
\end{equation}
where $I$ and $P$ are as in 2. and $h_P$ is the composition $P\r S\r I$. The middle term $X$
of this triangle is isomorphic to $I/S[-1]\oplus \rad P$.
\end{lemma}
\begin{proof} This follows using Proposition \ref{ref:I.2.3-15} and standard properties of almost
  split sequences and Auslander--Reiten triangles. We leave the proof
  to the reader.
\end{proof}
We now get the following main result on the connection between Serre duality and almost split
sequences.
\begin{theorem}
\label{ref:I.3.2-20} Let $\Cscr$ be an $\Ext$-finite hereditary abelian category.
\begin{enumerate}
\item $\Cscr$ has  Serre duality if and only if $\Cscr$ has almost split sequences, and there is a
one-one correspondence between indecomposable projective objects $P$ and indecomposable
injective objects $I $, via $P/ \rad P \cong \Soc I$.
\item If $\Cscr$  has no nonzero projective or injective objects, then $\Cscr$ has Serre
duality if and only if it has almost split sequences.
\end{enumerate}
\end{theorem}
Note that if $\Cscr$ is the category of finite dimensional representations over $k$ of the
quiver $\cdot \ra \cdot \ra \cdot \ra \cdots$, then $\Cscr$ is an $\Ext$-finite hereditary
abelian $k$-category with almost split sequences. Since $\Cscr$ has nonzero injective objects,
but no nonzero projective objects, it does not have Serre duality.

Let $\Pscr$ and $\Iscr$ be the full subcategories of $\Cscr$ whose objects are respectively
the projectives and the injectives in $\Cscr$. If $H$ is a set  of objects in $\Cscr$, then we
denote by $\Cscr_H$ the full subcategory of $\Cscr$ whose objects have no summands in $H$.

When $\Cscr =\mod A$ for a finite dimensional hereditary $k$-algebra $A$, we have an
equivalence $\tau:\Cscr_{\Pscr} \ra \Cscr_{\Iscr}$, where for an indecomposable object $C$ in
$\Cscr$ the object $\tau C$ is the left hand term of the almost split sequence with right hand
term $C$. Also we have the Nakayama functor $N:\Pscr \ra \Iscr$, which is an equivalence of
categories, where $N(P)=\Hom_k (\Hom_A (P,A),k)$ for $P\in \Pscr$. For the equivalence $F:
D^b(\Cscr) \ra D^b(\Cscr)$, where $F$ is the Serre functor, we have $F| \Cscr_{\Pscr} =\tau
[1]$ and $F|\Pscr =N$. Hence $F$ is in some sense put together from the two equivalences
$\tau$ and $N$ (see \cite{Happel}). Using Lemma \ref{ref:I.3.1-18} we see that the situation
is similar in the general case.
\begin{corollary}
\label{ref:I.3.3-21}
 Assume that $\Cscr$ has a {{right Serre}} functor $F$. Then
  the following hold:
\begin{enumerate}
\item $F$ defines a fully faithful functor $\Pscr\r \Iscr$. We denote this
  functor by $N$ (the ``Nakayama functor'').
\item $F [-1]=\tilde{\tau}$ induces a fully faithful functor $\Cscr_\Pscr\r
  \Cscr_\Iscr$, which we denote by $\tau$.
\item If $P\in\Pscr$ is indecomposable, then $\Hom(P,NP)$ has one
  dimensional socle, both as left $\Hom(NP,NP)$-module and as right
  $\Hom(P,P)$-module.  Let $h_P$ be a non-zero element in this
  socle. Then $S=\im h_P$ is simple, and furthermore $P$ is a
  projective cover of $S$ and $NP$ is an injective hull of $S$.
\end{enumerate}
If $F$ is a Serre functor then the functors $N$ and $\tau$ defined above are equivalences.
\end{corollary}
\begin{proof}
That $F$ takes the indicated values on objects follows from the nature of the
Auslander--Reiten triangles in $D^b(\Cscr)$ (given by lemma \ref{ref:I.3.1-18}). The assertion
about fully faithfulness of $\tau$ and $N$ follows from the corresponding property of $F$.
Part 3. follows by inspecting the triangle \eqref{ref:I.3.1-19}. Finally that $N$ and $\tau$
are equivalences in the case that $F$ is a Serre functor follows by considering $F^{-1}$,
which is a left Serre functor.
\end{proof}

In the next chapter we will  use the following result.
\begin{lemma}
\label{ref:I.3.4-22} Assume that $\Ascr$ is a triangulated category
with a
  $t$-structure $(\Ascr_{\le 0},\Ascr_{\ge 0})$ in such a way that
  every object in $\Ascr$ lies in some
  $\Ascr^{[a,b]}$.
Let $\Cscr$ be
  the heart of the $t$-structure, and assume that $\Hom(A,B[n])=0$ for
  $A,B\in\Cscr$ and $n\neq 0,1$. Then $\Cscr$ is a hereditary abelian
  category,
  and furthermore $\Ascr$ has a (right) Serre functor if and only if
  $\Cscr$ has a (right) Serre functor.
\end{lemma}
\begin{proof}
To show that $\Cscr$ is hereditary we have to show that
$\Ext^n_\Cscr(A,B)=\Hom_{\Ascr}(A,B[n])$. Now $\Ext^m_\Cscr(A,-)$
is characterized by the property that it is an effaceable
$\delta$-functor which coincides in degree zero with
$\Hom_\Cscr(A,-)$. Hence we have to show that
$\Hom_\Ascr(A,-[n])$ is effaceable. This is clear in degree $\ge 2$ since
there the functor is zero, and for $n\le 1$ it is also clear since
then $\Hom_{\Ascr}(A,-[n])=\Ext^n_{\Cscr}(A,-)$ \cite[p75]{BBD}.

Now standard arguments show that as additive categories $\Ascr$ and $D^b(\Cscr)$ are
equivalent to $\bigoplus_n \Cscr[n]$. We don't know if this equivalence yields an exact
equivalence between $\Ascr$ and $D^b(\Cscr)$, but recall that the definition of a (right)
Serre functor does not involve the triangulated structure. Hence if $\Ascr$ has a (right)
Serre functor then so does $D^b(\Cscr)$ and vice versa.
\end{proof}

\chapter[Hereditary noetherian abelian categories]{Hereditary noetherian
abelian categories with nonzero projective objects.} \label{ref:II-23} In this chapter we classify the
connected noetherian hereditary $\Ext$-finite abelian categories with Serre functor, in the
case where there are  nonzero projective objects.

Let $\Cscr$ be a connected hereditary abelian noetherian $\Ext$-finite category with Serre
functor and nonzero projective objects. The structure of the projective objects gives rise to
an associated quiver $Q$, which satisfies special assumptions, including  being locally finite
and having no infinite path ending at any vertex. The category $\Cscr$ contains the category
$\rep Q$ of finitely presented representations of $Q$ as a full subcategory, which may be
different from $\Cscr$. Actually $\rep Q$ is the full subcategory generated by projective
objects. Hence $\Cscr$ is not in general generated by the projective objects, but it is
generated by the preprojective objects. This provides a new interesting phenomenon.

Conversely, starting with a locally finite quiver $Q$ having no infinite path ending at any
vertex, the category $\rep Q$ is a hereditary abelian $\Ext$-finite category with a right
Serre functor (and right almost split sequences), but not necessarily having a Serre functor.
We construct a hereditary abelian $\Ext$-finite category with Serre functor $ \wrep(Q)$
containing $\rep Q$ as a full subcategory. We describe for which quivers $Q$ the category $
\wrep(Q)$  is noetherian, and show that in the noetherian case all connected noetherian
hereditary abelian $\Ext$-finite categories with Serre duality are of this form.

In Section \ref{ref:II.1-24} we explain the  construction of $\wrep(Q)$ via inverting a right
Serre functor. An alternative approach using derived categories is discussed in Section
\ref{ref:II.3-44}. After learning about our results Claus Ringel found yet another approach
towards the construction of $\wrep(Q)$ \cite{Ri1}.

In Section \ref{ref:II.2-32} we show that if $\Cscr$ is generated by the
preprojective objects, then $\Cscr$ is equivalent to some $\wrep(Q)$. That $\Cscr$ is
generated by preprojective objects when $\Cscr$ is connected and noetherian is proved in
Section \ref{ref:II.4-53}, along with showing for which quivers $Q$ the category $ \wrep(Q)$
is noetherian.

\section{Hereditary abelian categories constructed from
quivers}\label{ref:II.1-24} For a quiver $Q$ denote by $\rep(Q)$ the category of finitely
presented representations of $Q$. Under some additional assumptions on $Q$ we construct a
hereditary abelian $\Ext$-finite category $\wrep(Q)$ with Serre functor. The category
$\wrep(Q)$ contains $\rep(Q)$ as a full subcategory, and is obtained from $\rep(Q)$ by
formally inverting a right Serre functor.

Let $Q$ be a quiver with the following properties:
\begin{itemize}
\item[(P1)] $Q$ is locally finite, that is, every vertex in $Q$ is adjacent to only a finite number of other
  vertices.
\item[(P2)] There is no infinite path in $Q$ of the form $x_0\l x_1\l
  x_2\l\cdots$ (in particular there are no oriented cycles).
\end{itemize}
Note that (P1)(P2) imply the following.
\begin{lemma}
\label{ref:II.1.1-25}
\begin{enumerate}
\item
If $x$ is a vertex in $Q$, then there are only a finite number of
 paths ending in $x$.
\item
 If $x$ and $y$ are vertices in $Q$ then there are only a finite number
  of paths from $x$ to $y$.
\end{enumerate}
\end{lemma}
\begin{proof} Part 2. is an obvious consequence of 1., so we prove 1. Assume
  there is an infinite number of paths ending in $x$. By (P1) there
  must be an arrow $x_1\r x$ such that there is an infinite number of paths
  ending in $x_1$.  Repeating this we obtain an infinite path
  $\cdots\r x_2\r x_1\r x$ such that there is an infinite number of
  paths ending in every $x_n$.  The existence of such an infinite path
   contradicts (P2).
\end{proof}
For  a vertex $x$ in $Q$ we denote by $P_x$, $I_x$, $S_x$ respectively the corresponding
projective, injective  and simple object in $\Rep(Q)$, the category of all
$Q$-representations. By lemma \ref{ref:II.1.1-25} the objects $I_x$ have finite length. We
also have a canonical isomorphism
\begin{equation}
\label{ref:II.1.1-26} \Hom(P_x,P_y)\cong \Hom(I_x,I_y)
\end{equation}
since both of these vector spaces have as basis the paths from $y$ to $x$.

The category $\Rep(Q)$ is too big for what we want. For example it is not $\Ext$-finite. As
will be clear from the considerations in the next section the natural subcategory to consider
is $\rep(Q)$. It is easy to see that this is a hereditary abelian subcategory of $\Rep(Q)$,
and furthermore by lemma \ref{ref:II.1.1-25} it has finite dimensional $\Ext$'s. Thanks to
(P1) the simple and hence the finite dimensional representations are contained in $\rep(Q)$.
Hence in particular $I_x\in \rep(Q)$ for every vertex $x$.

Let $\Pscr$ and $ \Iscr$ be the full subcategories of $\rep(Q)$ consisting respectively of the
finite direct sums of the $P_x$ and the finite direct sums of the $I_x$. Let $N:\Pscr\r \Iscr$
be the equivalence obtained from additively extending \eqref{ref:II.1.1-26}. We denote also by
$N$ the corresponding equivalence $N:K^b(\Pscr)\r K^b(\Iscr)$. Finally denoting by $F$ the
composition
\begin{equation}
\label{ref:II.1.2-27} D^b(\rep(Q))\cong
K^b(\Pscr)\xrightarrow{N} K^b(\Iscr)\r D^b(\rep(Q))
\end{equation}
we have the following.
\begin{lemma} $F$ is a {{right Serre}} functor.
\end{lemma}
\begin{proof} Let $A,B\in K^b(\Pscr)$. Then we need to construct
  natural isomorphisms
\begin{equation}
\label{ref:II.1.3-28} \Hom_{K^b(\Pscr)}(A,B)\r
\Hom_{D^b(\rep(Q))}(B,FA)^\ast
\end{equation}
Since $A$ and $B$ are finite complexes of projectives we can reduce to the case $A=P_x$ and
$B=P_y$. So we need natural isomorphisms
\[
\Hom(P_x,P_y)\r \Hom(P_y,I_x)^\ast
\]
Again both of these vector spaces have a natural basis given by
the paths from $y$ to $x$. This leads to the required
isomorphisms.
\end{proof}
Under the assumptions (P1)(P2) on the quiver $Q$ we now know that $\rep Q$ is a hereditary
abelian $\Ext$-finite category with a right Serre functor $F$ such that the image of the
projective objects in $\rep Q$ are injective objects in $\rep Q$ of finite length. Hence the
following result applies to this setting.
\begin{theorem}
\label{ref:II.1.3-29} Let $\Bscr$ be an abelian $\Ext$-finite hereditary category with a right
Serre functor $F$ and enough projectives. Denote the full
subcategory of  projective objects in $\Bscr$ by $\Pscr$ and assume
 that $F(\Pscr)=\Iscr$ consists of (injective) objects of finite length. Then there
exists an $\Ext$-finite abelian hereditary category $\Cscr$ with the following properties.
\begin{enumerate}
\item There exists a full faithful exact embedding $i:\Bscr\r \Cscr$.
\item The injectives and projectives in $\Cscr$ are given by $i(\Iscr)$
  and $i(\Pscr)$.
\item $\Cscr$ possesses a Serre functor which extends
  $F:D^b(\Bscr)\r D^b(\Bscr)$ in such a
  way that there is a natural equivalence $\nu:Fi\r iF$ where we have
  denoted the derived functor of $i$ also by $i$ and the extended Serre functor also by $F$.

\item For every indecomposable object $X\in \Cscr$ there exists $t\ge 0$
  such that $\tau^t X$ is defined and lies in $\Bscr$, where $\tau$ denotes the  functor $\Cscr_\Pscr\r
  \Cscr_\Iscr$ (see Corollary \ref{ref:I.3.3-21}.2) induced by $F[-1]=\tilde{\tau}$.
\end{enumerate}
Furthermore a quadruple $(\Cscr,i,F,\nu)$ with these properties is
unique in the appropriate sense.
\end{theorem}
\begin{proof}
First we  show that if $\Cscr$ exists satisfying properties $1.,2.,3.,4.$ then it is unique.
This proof will in particular tell us how to construct $\Cscr$. We then show that this
construction always yields a category $\Cscr$ with the required properties.

 Since $D^b(\Bscr)=K^b(\Pscr)$ and
since the objects in $\Pscr$ remain projective in $\Cscr$ by property
2., it follows that the derived functor of $i$
is fully faithful. It is also clear that $D^b(\Bscr)$ is closed inside $D^b(\Cscr)$ under the
formation of cones.

Let $F:D^b(\Cscr)\r D^b(\Cscr)$ be the extended Serre functor. If $K\in D^b(\Cscr)$ then $K$
can be obtained by starting with objects in $\Cscr$ and repeatedly taking cones. It then
follows from 4. that for $t\gg 0$ we have $F^tK\in D^b(\Bscr)$. Hence the triple
$(D^b(\Cscr),i,F)$ satisfies the hypotheses of Proposition \ref{ref:I.1.8-8}. Thus we obtain
$D^b(\Cscr)=F^{-\infty}D^b(\Bscr)$.

\begin{comment}
To simplify the notation we will no longer explicitly write $i$.
The indecomposable objects in $F^{-\infty}D^b(\rep(Q))$ are of the
form $F^{-a}A[l]$ with $A$ an indecomposable object in $\rep(Q)$.
We now have to determine which of these indecomposable objects
actually lie in $\Cscr$. This will complete the description of
$\Cscr$.

Let us temporarily write $G=F[-1]$. It is sufficient to determine
which of the objects $G^{-a}A[l]$ lie in $\Cscr$. Following the
sequence $A,G^{-1}A,G^{-2}A,\ldots$ and using the fact that the
injectives are given by $\Iscr$, we find that the indecomposable
objects in $\Cscr$ are of the form
\end{comment}
According to 4. every indecomposable object in $\Cscr$ will be of the form $\tau^{-t} Y$ with
$Y\in \Bscr$. Implicit in the notation $\tau^{-t}Y$ is the assumption that $\tau^{-l}Y$ is
defined for $l=1,\ldots,t$, that is $\tau^{-i}Y\not\in \Iscr$ for $i=0,\ldots,t-1$. In
addition we may assume that $t$ is minimal. Thus either $t=0$ or else $\tau^{-l}Y\not\in
\Bscr$ for $l=1,\ldots,t$. The last case is equivalent to
\begin{equation}
\label{ref:II.1.4-30} \tau^{-i}Y\not \in \tau \Bscr_\Pscr\cup \Iscr
\end{equation}
 for
$i=0,\ldots,t-1$.

If $Y\not\in \tau \Bscr_\Pscr\cup \Iscr$ then the same is true for $\tau^{-1}Y$. Hence we have
to impose \eqref{ref:II.1.4-30} only for $i=0$.

We conclude that the indecomposable objects in $\Cscr$ are of the
following form:
\begin{itemize}
\item[(C1)] The indecomposable objects in $\Bscr$.
\item[(C2)] Objects of the form $\tilde{\tau}^{-t}Y$ where $t>0$ and $Y$ is an
  indecomposable object in $\Bscr\setminus (\tau \Bscr_{\Pscr}\cup
  \Iscr)$, where $\tilde{\tau}=F[-1]$.
\end{itemize}
This completes the determination of $\Cscr$ as an additive subcategory of $F^{-\infty}D^b(B)$,
and finishes the proof of the uniqueness.

\medskip

Let us now assume that $\Cscr$ is the additive subcategory of $F^{-\infty}D^b(\Bscr)$ whose
indecomposable objects are given by (C1)(C2). We have to show that $\Cscr$ is a hereditary
abelian category satisfying 1.-4.

Since $F^{-\infty}D^b(\Bscr)$ has a Serre functor, it has Auslander--Reiten triangles.  Below
we will need the triangle associated to $P_x[1]$. Using the criterion given in Proposition
\ref{ref:I.2.1-12} we can compute this triangle in $D^b(\Bscr)$. So according to
Lemma \ref{ref:I.3.1-18} the requested triangle  is of the form
\begin{equation}
\label{ref:II.1.5-31} I_x\r I_x/\Soc(I_x)\oplus  \rad(P_x)[1]\r
P_x[1]\r
\end{equation}
%(note that the maps in this triangle are not the obvious ones.)
 We now define a $t$-structure
on $F^{-\infty}D^b(\Bscr)$. Using the fact that $\Bscr$ is hereditary we easily obtain that as
additive categories $F^{-\infty}D^b(\Bscr)=\bigoplus_n \Cscr[n]$. We now define
\begin{align*}
F^{-\infty}D^b(\Bscr)_{\le 0}&=\bigoplus_{n\ge 0} \Cscr[n]\\ F^{-\infty}D^b(\Bscr)_{\ge
0}&=\bigoplus_{n\le  0} \Cscr[n]\\
\end{align*}
We claim that this is a $t$-structure. The only non-trivial axiom
we have to verify is that
\[
\Hom(F^{-\infty}D^b(\Bscr)_{\le 0},F^{-\infty}D^b(\Bscr)_{\ge 1})=0
\]
So this amounts to showing that $\Ext^{-i}(A,B)=0$ for
$A,B\in\Cscr$ and for $i>0$. We separate this into four cases.
\begin{case} $A$ and $B$ fall under (C1). This case is trivial.
\end{case}
\begin{case} $A$ falls under (C1) and $B$ falls under (C2). Thus
  we have $B=\tilde{\tau}^{-t} X$ for some $X$ in $\Bscr$ and $t>0$. We want to show that
  $\Ext^{-i} (A,B)=0$ for $i>0$ for any $A\in \Bscr$, by induction on $t$.

  Let first $t=1$. We can assume that $A$ is projective, since otherwise we could reduce to
  Case 1 by applying $\tilde{\tau}$. Then $\tilde{\tau} (A)=I[-1]$ for some object $I$ in
  $\Iscr$, so it is sufficient to show $\Ext^{-i} (I[-1],X)\cong \Ext^{-(i-1)} (I,X)=0$ for
  $i>0$. For $i>1$ this follows from Case 1. For $i=1$ we clearly have $\Hom (I,X)=0$ since
  $I$ is injective in $\Bscr$ and $X$ is not injective.

  Assume now that $t>1$ and that the claim has been proved for $t-1$ for all objects $A$ in
  $\Bscr$. Then we can again assume that $A$ is projective, and consider $\Ext^{-(i-1)}
  (I,\tilde{\tau}^{-t+1} X)$, where $\tilde{\tau} (A)=I[-1]$. We want to prove that
  $\Ext^{-(i-1)} (I,\tilde{\tau}^{-t+1} X)=0$ for $i>0$. By the induction assumption we only
  need to consider $i=1$. So we want to prove that $\Hom (J, \tilde{\tau}^{-t+1} X)=0$ for $J$
  indecomposable in $\Iscr$, and we do this by induction on the length of $J$ which by
  assumption is finite. We have the Auslander--Reiten triangle $J\stackrel{g}{\r} J/\Soc J \oplus (\rad
  P)[1] \r P[1] \r J[1]$, with $P$ indecomposable projective in $\Pscr$. Then any nonzero map
  $h:J\r \tilde{\tau}^{-t+1} X$ would factor through $g$ since $h$ is not an isomorphism. Note
  that by the induction assumption we have $\Hom ((\rad P)[1],\tilde{\tau}^{-t+1} X) =\Ext^{-1}
  (\rad P, \tilde{\tau}^{-t+1} X)=0$, so the claim follows.
\end{case}
\begin{case} $A$ falls under (C2) and $B$ falls under (C1). Thus
  $A=\tilde{\tau}^{-t}X$ for some $X$ in $\Bscr$ and $t>0$. We have $\Ext^{-i}(A,B)=\Ext^{-i}(X,\tilde{\tau}^tB)$. Since
  $\tilde{\tau}^tB\in\Bscr[l]$ for $l\le 0$, this case is trivial.
\end{case}
\begin{case} $A$ and $B$ both fall under (C2). This case can be
  reduced to one of the previous cases by applying powers of $\tilde{\tau}$.
\end{case}
So $\Cscr$ is indeed an abelian category. Since it is easily seen that $F(\Cscr)\subset
\Cscr\oplus \Cscr[1]$, we obtain by Serre duality that $\Hom(\Cscr,\Cscr[n])=0$ for $n\neq
0,1$. Hence  it follows by lemma \ref{ref:I.3.4-22} that $\Cscr$ is hereditary.

Now we verify properties 1.-4..
\begin{itemize}
\item[1.] This is clear from the construction.
\item[3.] This follows from lemma \ref{ref:I.3.4-22}.
\item [4.] This is clear from the construction.
\item [2.] By 3. we already know that $\Cscr$ has a Serre functor, so
  we can use its properties. An indecomposable object $P$ in $\Cscr$ is
  projective precisely when $FP\in\Cscr$. Using the construction of
  $\Cscr$ it is easy to see that this happens if and only if
  $P\in\Pscr$. Using the properties of Serre functors we find that the
  injectives in $\Cscr$ must be given by $F(\Pscr)=\Iscr$. This proves
  what we want.\qed
\end{itemize}
\def\qed{}\end{proof}
If $Q$ is a quiver satisfying (P1)(P2) and $\Bscr=\rep(Q)$, then in the sequel we will denote by $\wrep(Q)$ the
hereditary abelian category $\Cscr$ whose existence is asserted by Theorem \ref{ref:II.1.3-29}.

\begin{example}

Consider the example
\[
0\r 1\r 2\r\cdots
\]

Since the category $\rep (Q)$ has right almost split sequences, it makes sense to talk about
an associated right AR-quiver, which looks like

\begin{equation}
  \xymatrix@!0@=0.6cm{
    & & & & \\
    & & & \ar[ur]_>{P_0} \\
    & & \ar[ur]_>{P_1} \\
    & \ar[ur]_>{P_2} \\
    \ar@{.}[ur]_>{P_3}
    }
  \quad
  \xymatrix@!0@=0.6cm{
    & \ar[dr] & & \ar[dr] & & \ar[dr] & & \\
    \ar@{.}[r] & & \ar[dr] \ar[ur] & & \ar[dr] \ar[ur] & & \ar[ur]_>{I_0} \\
    & \ar[dr] \ar[ur] & & \ar[dr] \ar[ur] & & \ar[ur]_>{I_1} \\
    & & \ar[dr] \ar[ur] & & \ar[ur]_>{I_2} \\
    & & & \ar[ur] \\
    & & \ar@{.}[ur]
    }
\end{equation}

The AR-quiver for the new category $\wrep(Q)$ consists of the two components

\begin{equation}
  \xymatrix@!0@=0.6cm{
    & & & & \ar[dr]^<{P_0} & & \ar[dr] & & \\
    &&& \ar[ur] \ar[dr]^<{P_1} && \ar[ur] \ar[dr] && \ar[ur] & \ar@{.}[r] & \\
    && \ar[ur]_<{P_2} & & \ar[ur] \ar[dr] & & \ar[ur] & & \\
    & \ar@{.}[dl] \ar[ur] & & & & \ar[ur] & \\
    &
    }
  \quad
  \xymatrix@!0@=0.6cm{
    & \ar[dr] & & \ar[dr] & & \\
    \ar@{.}[r] & & \ar[ur] \ar[dr] & & \ar[ur]_>{I_0} \\
    & \ar[ur] \ar[dr] & & \ar[ur]_>{I_1} \\
    & & \ar@{.}[dl] \ar[ur]_>{I_2} \\
    &
    }
\end{equation}

The new objects are the $\tau^{-i} P_j$ for $j\geq 0$ and $i>0$.

\end{example}

\begin{remark}
It is not necessarily true that all injectives in $\rep(Q)$ are in
$\Iscr$.
E.g. in the above example the projective representation associated to $0$ is also
injective. However it is not in $\Iscr$.
\end{remark}

\section{Hereditary abelian categories generated by preprojectives}
\label{ref:II.2-32} Let $\Cscr$ be an $\Ext$-finite hereditary abelian category. In this
section we associate a quiver $Q$ with $\Cscr$. We show that if $\Cscr$ has Serre duality, the
injective objects in $\Cscr$ are noetherian and $\Cscr$ is generated by the preprojective
objects, then $\Cscr$ is equivalent to $\wrep(Q)$.

For an $\Ext$-finite hereditary abelian category $\Cscr$, denote as before by $\Pscr$ and $
\Iscr$  respectively the full subcategories of $\Cscr$ consisting of projective and injective
objects. By $\overline{\Pscr}$ we denote the full subcategory of $\Cscr$ whose objects are
quotients of objects in $\Pscr$. For later reference we state the following easily proved
fact.
\begin{lemma}
$\overline{\Pscr}$ is closed under subobjects and extensions. In particular $\overline{\Pscr}$
is an abelian category.
\end{lemma}
To $\Cscr$ we associate a quiver $Q$ whose vertices are in one-one correspondence with the
indecomposable projectives in $\Cscr$. If $x$ is a vertex then we denote the corresponding
projective object in $\Cscr$ by $\frak{p}_x$. For $x$ and $y$ vertices in $Q$ we let the
arrows from $x$ to $y$ index elements $(f_{y,x}^{i})_{i}\in\rad\Hom(\frak{p}_y,\frak{p}_x)$,
which are representatives for a basis of $\rad\Hom(\frak{p}_y,\frak{p}_x)
/\rad^2\Hom(\frak{p}_y,\frak{p}_x)$. Somewhat inaccurately we will say below that the
projectives in $\Cscr$ are given by $Q$. We now have a functor
\[
i:\rep(Q)\r \bar{\Pscr} : P_x\mapsto \frak{p}_x
\]
where $P_x $ is the projective representation of $Q$ associated with the vertex $x$. It is
easy to see that this functor is faithful, but it is not necessarily full, as the following
example shows.
\begin{example} Let $\Gr^\dagger k[x]$ be the category of graded
  $k[x]$-modules which are contained in a finite direct sum of
  indecomposable graded injective $k[x]$-modules, and let
  $\Cscr=(\Gr^\dagger k[x])^{\circ pp}$. {}From the fact that $k[x]$ is
  hereditary one deduces that the same holds for $\Cscr$. The indecomposable
  graded injective $k[x]$-modules are given by $E_\infty=k[x,x^{-1}]$
  and $E_a=x^ak[x^{-1}]$ for $a\in\ZZ$. These correspond to projective
  objects $P_a$ in $\Cscr$ for $a\in\ZZ\cup\{\infty\}$.
  We have $\Hom_\Cscr(P_a,P_\infty)=\Hom_{\Gr
  k[x]}(E_\infty,E_a)=k$. Thus $\rad\Hom_\Cscr(P_\infty,P_\infty)
  =0$. Furthermore if $a\neq \infty$ then every map $P_a\r P_\infty$
  factors through $P_{a+1}$. It follows that $\rad
  \Hom_\Cscr(P_a,P_\infty)=\rad^2\Hom_\Cscr(P_a,P_\infty)$ and
  thus $P_\infty$ is an isolated point in the quiver associated to
  $\Cscr$. This contradicts fullness.
\end{example}
We now have the following result.
\begin{proposition}
\label{ref:II.2.3-33}
 Assume that $\Cscr$ has a Serre functor. Then the
  following are equivalent.
\begin{enumerate}
\item $\Cscr$ has finite length injectives.
\item $\Cscr$ has noetherian injectives.
\item The quiver $Q$  satisfies (P1)(P2).
\item The quiver $Q$  satisfies (P1)(P2) and $i: \rep (Q)\r \overline{\Pscr}$ defines an
  equivalence between $\rep(Q)$ and $\overline{\Pscr}$. In addition
  the indecomposable injectives in $\Cscr$
  are of the form $i(I_x)$ for $x\in Q$, so that in particular $\Iscr\subset
  \overline{\Pscr}$.
\end{enumerate}
\end{proposition}
\begin{proof}\strut
\begin{itemize}
\item[$1.\Rightarrow 2.$] This is clear.
\item[$2.\Rightarrow 1.$] Assume that $I$ is an injective object in
  $\Cscr$. By decomposing $I$ into a direct sum of indecomposables and invoking
Theorem \ref{ref:I.3.2-20} we see that the socle $S$ of $I$ is non-zero and of finite length.

  Since $\Cscr$ is hereditary, $I/S$ is again injective. Repeating
  this we find a strictly ascending chain
\[
0=I_0\subset I_1\subset \cdots \subset I_n\subset\cdots\subset I
\]
such that $I_{n+1}/I_n$ has finite length. Since $I$ is
noetherian, this chain must stop, whence $I$ has finite length.
\item[$1.\Rightarrow 4.$] We first  show that $Q$ satisfies (P2).
In fact we show that for a vertex
$z$ in $Q$ there is  a bound on the length of a chain of
non-isomorphisms
\begin{equation}
\label{ref:II.2.1-34} \frak{p}_z\r \frak{p}_{z_1}\r \cdots \r
\frak{p}_{z_n}.
\end{equation}
Put $\frak{i}_x=F(\frak{p}_x)$. Applying $F$ we get non-isomorphisms
\[
\frak{i}_z\r \frak{i}_{z_1}\r \cdots \r \frak{i}_{z_n}
\]
Since these are indecomposable injectives and $\Cscr$ is
hereditary, all these maps must be surjective. Clearly the
required bound is now given by the length of $\frak{i}_z$.

Now we prove that $Q$ satisfies (P1). Since  every $\rad(\frak{p}_x)$ has only a finite number
of summands in a direct sum decomposition, it is sufficient to show that there exists only a
finite number of vertices $y$ such that there is a non-zero map $\frak{p}_x\r \frak{p}_y$.
Applying $F$ we find that $\Soc(\frak{i}_y)$ must be a subquotient of $\frak{i}_x$. Since
$\frak{i}_x$ has finite length, there is only a finite number of possibilities for $y$.

 The fact that there is a bound (depending
only on $z$) on chains of the form \eqref{ref:II.2.1-34} implies
that every map $\frak{p}_z\r \frak{p}_t$ is a linear combination
of products of the $f_{y,x}^i$. Thus $i$ is full, and from this we
easily obtain that $i$ yields an equivalence between $\rep(Q)$ and
$\overline{\Pscr}$.

Now let $I$ be indecomposable injective in $\Cscr$. Since $\Cscr$ has a Serre functor, it
follows from Theorem \ref{ref:I.3.2-20} that $I$ is the injective hull of some simple object
$S$ lying in $\bar{\Pscr}$. By considering the injective $I/S$ and using induction on the
length of $I$ we find that $I\in \bar{\Pscr}$. Now $I$ is clearly  the injective hull of $S$
also in $\bar{\Pscr}$. Using the fact that $i$ is an equivalence there must exist a vertex
$x\in Q$ such that $S=i(S_x)$. Furthermore $I$ must correspond under $i$ to the injective hull
of $S_x$. Since the latter is $I_x$, we are done.
\item[$4.\Rightarrow 3.$] This is clear.
\item[$3.\Rightarrow 1.$] This proof is similar to that of
  $2.\Rightarrow 1.$ Assume that $I$ is an indecomposable injective
  of infinite length, and let $S$ be its socle. There exists an indecomposable summand
  $I_1$ of $I/S$ which is of infinite length. Continuing this
  procedure we find an infinite sequence of irreducible maps
\[
I\rightarrow I_1\rightarrow I_2\rightarrow\cdots
\]
and applying $F^{-1}$ we find a corresponding infinite sequence of
irreducible maps between projectives
\[
P\rightarrow P_1\rightarrow P_2\r\cdots
\]
However such an infinite sequence cannot exist by (P2). \qed
\end{itemize}
\def\qed{}\end{proof}
We have the following result on reconstructing $\Cscr$ from the associated quiver $Q$.
\begin{corollary}
\label{ref:II.2.4-35} Assume that $\Cscr$ is an
  $\Ext$-finite hereditary abelian category possessing a Serre functor $F$. Let
  the projectives in $\Cscr$ be given by the quiver $Q$. Assume that
  $\Cscr$ satisfies any of the conditions of Proposition
  \ref{ref:II.2.3-33} and furthermore that every
  indecomposable in $\Cscr$ is of the form $\tau^{-t} X$ with $X\in
  \bar{\Pscr}$. Then $\Cscr$ is equivalent to $\wrep(Q)$.
\end{corollary}
\begin{proof} {}From the fact that $F(\Pscr)=\Iscr\subset
  \bar{\Pscr}$ and the fact that $D^b(\bar{\Pscr})$ is closed under
  cones in $D^b(\Cscr)$, it follows that $D^b(\bar{\Pscr})$ is closed under
  $F$. Thus $F$ defines a right Serre functor on $\bar{\Pscr}$, and
  hence via $i:\rep (Q)\r \Cscr$ a right Serre functor on $\rep(Q)$. Since {{right
      Serre}} functors are unique, $F$ coincides (up to a
  natural isomorphism) with the standard
  {{right Serre}} functor obtained from deriving the Nakayama functor
  which was introduced in \eqref{ref:II.1.2-27} .

It now follows that $\Cscr$ satisfies  properties 1.,2.,3.,4. of Theorem \ref{ref:II.1.3-29}.
Hence we get $\Cscr=\wrep(Q)$.
\end{proof}

We remind the reader of some elementary facts concerning preprojective objects. In general if
$\Cscr$ is a hereditary abelian category with almost split sequences, then a preprojective
object is by definition an object of the form $\oplus_{i=1}^n \tau^{-a_i} P_i$ where the
$P_i$'s are indecomposable projectives and $a_i\ge 0$. Similarly a preinjective object is an
object of the form $\oplus_{i=1}^n \tau^{a_i} I_i$ where the $I_i$ are indecomposable
injective and $a_i\ge 0$.

For simplicity we  restrict ourselves here to the case where $\Cscr$ is  $\Ext$-finite with a
Serre functor. In this case we have the functor $\tau:\Cscr_{\Pscr} \r \Cscr_{\Iscr}$ inducing
the correspondence between the end terms of an almost split sequence.

Note the following.
\begin{lemma}
\label{ref:II.2.5-36}
 Let $U\in\Cscr$ be indecomposable, and assume  there is a
  non-zero map $\phi:U\r \tau^{-a}P$ where $P$ is an indecomposable
  projective and $a\ge 0$. Then $U$ is of the form $\tau^{-b}Q$ where $Q$
  is an indecomposable projective and $b\le a$.
\end{lemma}
\begin{proof}
  If $U$ is projective, then we are done. If $a=0$ then the fact that
  $U$ is indecomposable plus the fact that $\Cscr$ is hereditary
  implies that $U$ is projective. Hence this case is covered also.
  Assume now that $U$ is not projective (and hence $a\neq 0$). Then by
  faithfulness of $\tau$ we obtain a non-zero map $\tau U\r \tau^{-a+1}P$.
  Induction on $a$ now yields that $\tau U=\tau^{-c}Q$ with $c\le a-1$ for
  some indecomposable projective $Q$. Thus $U=\tau^{-c-1}Q$ and we are
  done.
\end{proof}
\begin{corollary}
\label{ref:II.2.6-37}
\begin{enumerate}
\item
 Assume that we have a map $\phi:U\r Q$ where $Q$ is
  preprojective. Then $U=U'\oplus U''$ where $U'\subset\ker \phi$
  and $U''$ is preprojective.
\item Every subobject of a preprojective object is preprojective.
\end{enumerate}
\end{corollary}
\begin{proof} 2. follows trivially from 1., so we prove 1. Let $U=U_1\oplus\cdots \oplus U_m$ be a decomposition
  of $U$ into a direct sum of indecomposable objects. If $U_i\not\subset \ker\phi$
  then it maps non-trivially to $Q$, whence it is preprojective by the
  previous lemma. This proves what we want.
\end{proof}
We say that $\Cscr$ is generated by preprojectives if every object
in $\Cscr$ is a quotient of a preprojective object. One has the
following result.
\begin{lemma}
\label{ref:II.2.7-38} The following are equivalent for an $\Ext$-finite hereditary abelian
category with Serre functor.
\begin{enumerate}
\item $\Cscr$ is generated by preprojectives.
\item Every indecomposable object in $\Cscr$ is of the form $\tau^{-t}X$ with
  $X\in\bar{\Pscr}$ and $t\geq 0$.
\end{enumerate}
\end{lemma}
\begin{proof}
Let us first assume 2. and let $Y=\tau^{-t}X\in \Cscr$, with $X\in \overline{\Pscr}$ and
$t\geq 0$. We have to show that $Y$ is the quotient of a preprojective object. Take a minimal
projective presentation
\begin{equation}
\label{ref:II.2.2-39} 0\r Q\r P\r X\r 0
\end{equation}
 with $P,Q\in \Pscr$. For
$l<t$ the object $\tau^{-l}Q$ cannot contain an injective summand since then the resolution
\eqref{ref:II.2.2-39} wasn't minimal.
It follows that the object $\tau^{-l}P$ cannot contain an injective summand since then $\tau^{-l}X$
would be injective, contradicting the fact that it is in the essential image of $\tau$.
 Hence $\tau^{-t}Q$ and $\tau^{-t} P$ are defined, and by
the exactness of $F$ (see \cite{Bondal4}) we obtain an exact sequence
\begin{equation}
\label{ref:II.2.3-40} 0\r \tau^{-t}Q\r \tau^{-t}P\r \tau^{-t}X\r 0
\end{equation}
In particular $\tau^{-t}X$ is covered by a preprojective object.

Now assume 1. Let $Y$ be an indecomposable object in $\Cscr$, and assume that there is a
surjective map $\phi:R\r Y$ with $R$ preprojective. Then by Corollary \ref{ref:II.2.6-37} we
have that  $S=\ker \phi$ is also preprojective. Thus for large $t$ we have a triangle
\[
F^t S\r F^t R \r F^t Y \r
\]
with $F^tR$ and $F^t S$ in the image of $D^b(\bar{\Pscr})$. Since $D^b(\bar{\Pscr})$ is closed
under cones in $D^b(\Cscr)$ and since $F^tY$ is indecomposable, it follows that $F^tY\in
\bar{\Pscr}[l]$ for some $l$.

Thus if $\tau^pY$ is defined for $1\le p\le t$ then $\tau^tY\in \bar{\Pscr}$. If $\tau^pY$ is
not defined, then $\tau^{p-1}Y$ is projective and so lies in $\Pscr$. Hence in this case we
are done also.
\end{proof}
Combining the last lemma with Corollary \ref{ref:II.2.4-35} we
obtain the following:
\begin{corollary}
\label{ref:II.2.8-41} Assume that $\Cscr$ is an $\Ext$-finite hereditary abelian category
possessing a Serre functor $F$. Let the projectives in $\Cscr$ be given by the quiver $Q$.
Assume that $\Cscr$ is generated by preprojectives and satisfies any of the conditions of
Proposition \ref{ref:II.2.3-33}. Then $\Cscr$ is equivalent to $\wrep(Q)$.
\end{corollary}
We will also need the following result.
\begin{theorem}
\label{ref:II.2.9-42} Let $\Qscr$ be the full subcategory of $\Cscr$ whose objects are
quotients of preprojective objects. Then $\Qscr$ is closed under subquotients and
  extensions. If in addition $\Cscr$ satisfies any of the conditions
  of Proposition \ref{ref:II.2.3-33}, then the Serre functor on
  $D^b(\Cscr)$ restricts to a Serre functor on $D^b(\Qscr)$.
\end{theorem}
\begin{proof}
Closedness under quotients is clear. Let us now prove closedness
under subobjects. Clearly it is sufficient to show that a
subobject of a preprojective object is preprojective. But this is
precisely Corollary \ref{ref:II.2.6-37}.

Let us now prove closedness under extensions. So assume that we
have an exact sequence
\begin{equation}
\label{ref:II.2.4-43} 0\r A\r B\r C\r0
\end{equation}
where $A,C\in\Qscr$. Since we already know that $\Qscr$ is closed under quotients, we may
prove our result for pullbacks of \eqref{ref:II.2.4-43}. In particular we may assume that $C$
is preprojective. But then according to Corollary \ref{ref:II.2.6-37} we have $B=A_1\oplus Q$
where $A_1\subset A$ and $Q$ is preprojective. Hence by what we know already we have
$A_1\in\Qscr$. This proves what we want.

Assume now that $\Cscr$ satisfies any of the conditions
  of Proposition \ref{ref:II.2.3-33}. We have to show that
  $D^b(\Qscr)$ is closed under $F$ and $F^{-1}$ in $D^b(\Cscr)$. Since
  $D^b(\Qscr)$ is closed under the formation of cones, it suffices to
  show that $F(P)$ and $F^{-1}(P)$ are in $D^b(\Qscr)$ for a
  preprojective object $P$. The only case which isn't entirely obvious  is that
  $F(P)$ is in $D^b(\Qscr)$ when $P$ is projective, but this follows from Proposition
  \ref{ref:II.2.3-33}.
\end{proof}

For completeness let us include the following lemma.
\begin{lemma} Let $Q$ be a quiver satisfying (P1)(P2). Then $Q$ is
  connected if and only if $\wrep(Q)$ is connected.
\end{lemma}
\begin{proof}
Since clearly $\wrep(Q_1\coprod Q_2)=\wrep(Q_1)\oplus \wrep(Q_2)$ it suffices to prove the
implication that $Q$ connected implies $\wrep(Q)$ connected. So assume that $\Cscr=\wrep(Q)$
is a direct sum $\Cscr=\Cscr_1\oplus \Cscr_2$ with $\Cscr_1$ and $\Cscr_2$ non-trivial. Let
$\Pscr_i$ be the category of projectives in $\Cscr_i$ for $i=1,2$. Then  clearly
$\Pscr=\Pscr_1\oplus\Pscr_2$. This yields a corresponding decomposition $Q=Q_1\coprod Q_2$.
Since we had assumed that $Q$ is connected, it follows that for example $Q_2=0$.

By lemma \ref{ref:I.1.11-10} the Serre functor on $D^b(\Cscr)$
restricts to one on $D^b(\Cscr_1)$ and one on $D^b(\Cscr_2)$. It
is now easy to see that the conditions for Corollary
\ref{ref:II.2.4-35} descend to $\Cscr_1$ and $\Cscr_2$. Thus
$\Cscr_i=\wrep(Q_i)$ and hence $\Cscr_2=0$. This is a
contradiction.
\end{proof}

\section{Derived equivalences}\label{ref:II.3-44}
In this section we construct a derived  equivalence between $\wrep(Q)$ and $\wrep(Q')$ for
certain quivers $Q$ and $Q'$ satisfying (P1)(P2).  This result will be used in the proof of
Theorem \ref{theoremc}. At the same time we get an alternative method for the construction of
$\wrep Q$ from $\rep Q$.

If $Q$ is a quiver, then $\ZZ Q$ denotes  the quiver whose vertices are of the form $(n,a)$
with $n\in \ZZ$ and $a\in Q$ and whose arrows are of the form $(n,a)\r (n,b)$ and $(n,b)\r
(n+1,a)$ for any arrow $a\r b$.  Then $\ZZ Q$ is a translation quiver with translation given
by $\tau(n,a)=(n-1,a)$. The following is easy to see.
\begin{lemma}
\label{ref:II.3.1-45} Assume that $Q$ satisfies (P1)(P2). Let
$x,y\in
  \ZZ Q$. Then  the number of  paths in $\ZZ Q$ from $x$ to
  $y$ is finite.
\end{lemma}
For a quiver $Q$ satisfying (P1)(P2), denote by $\NN Q$  the part of $\ZZ Q$ given by the
vertices $(n,a)$ where $n\ge 0$. Let  $Q'\subset \ZZ Q$ be a section of $\ZZ Q$. By this we
mean the following: $Q'$ contains exactly one vertex from each $\tau$-orbit of $\ZZ Q$ and if
$x\in Q'$ and there is an arrow $x\r y$ in $\ZZ Q$, then either the arrow $x\r y$ or the arrow
$\tau y\r x$ is in $Q'$, and if there is an arrow $z\r x$, then either the arrow $z\r x$ or
the arrow $x\r \tau^{-1} z$ is in $Q'$. It is then clear that if $Q'$ is a section in $\ZZ Q$
then $\ZZ Q'=\ZZ Q$ and $Q$ is a section in $\ZZ Q'$. Our main result    is that the
categories $\wrep(Q)$ and $\wrep(Q')$ are derived equivalent when $Q'$ is a section in $\ZZ
Q$. This is well-known when $Q$ is  finite, in which case $\wrep(Q)=\rep(Q)$ and
$\wrep(Q')=\rep(Q')$, and one can go from $Q$ to $Q'$ by a finite sequence of reflections.
\begin{lemma}
\label{ref:II.3.2-46} Let $\Cscr$ be an $\Ext$-finite
hereditary
  abelian category possessing a Serre functor $F$. Assume that the
  projectives in $\Cscr$ are given by a quiver satisfying (P1)(P2).

Then the   preinjective
  objects in $\Cscr$ have finite length, and furthermore they are
  quotients of projectives.
\end{lemma}
\begin{proof}  Let $\Pscr$, $\Iscr\subset \Cscr$ have their usual
  meaning. Since $F(\Pscr)=\Iscr$, we easily obtain
  $F(D^b(\Pscr))=D^b(\Iscr)$, and by Proposition
  \ref{ref:II.2.3-33} we have that the objects in $\Iscr$ have
  finite length and furthermore $D^b(\Iscr)\subset D^b(\Pscr)$. Thus
  in particular $F(D^b(\Iscr))\subset D^b(\Iscr)$. Iterating we find
  $F^a\Iscr\subset D^b(\Iscr)$ for any $a>0$. This proves  that the preinjective
  objects have finite length. We also obtain $F^a\Iscr\subset
  D^b(\Pscr)$, which yields that the preinjective objects are quotients
  of projectives.
\end{proof}
Also recall that a quiver $Q$ is Dynkin if the underlying graph is a Dynkin diagram, that is,
of the form $A_n$, $D_n$ or $E_6$, $E_7$, $E_8$, where $n$ denotes the number of vertices. Note
that $\rep Q$ is of finite representation type for a connected quiver $Q$ if and only if $Q$
is a Dynkin quiver.

We shall also need the infinite graphs $A_\infty$, $A_\infty^\infty$ and $D_\infty$

\begin{align}
\xymatrix@=0.4cm{
A_\infty: & \cdot\ar@{-}[r]    & \cdot\ar@{-}[r] & \cdot\ar@{-}[r] & \cdot\ar@{-}[r] &\cdot\ar@{.}[r] & }\notag \\
\xymatrix@=0.4cm{
A^\infty_\infty: & \, \ar@{.}[r] & \cdot\ar@{-}[r] & \cdot\ar@{-}[r] & \cdot\ar@{-}[r] & \cdot\ar@{.}[r] & }\notag\\
\xymatrix@=0.4cm@R=0.3cm{
          &  \cdot \ar@{-}[rd] &                 &                 &                  &                 &\\
D_\infty: &    & \cdot\ar@{-}[r] & \cdot\ar@{-}[r] & \cdot\ar@{-}[r]  & \cdot\ar@{.}[r] & \\
          &  \cdot \ar@{-}[ru] &                 &                 &                  &         &}\notag
\end{align}

We also need the following. Here the preprojective component of the AR-quiver denotes the component
containing the projective objects and the preinjective component denotes the component
containing the injective objects.
\begin{lemma} Assume that $Q$ is a connected quiver which is not Dynkin but satisfies
  (P1)(P2), and let $\Cscr$ be an  $\Ext$-finite
  hereditary abelian category with Serre functor, whose projectives are
  given by $Q$.
\begin{enumerate}
\item
The preprojective component $\Pi$ of $\Cscr$
  contains no injective objects.
\item There is an isomorphism $\Pi\cong \NN Q^{\op}$ of quivers, compatible with the
  translation.
\end{enumerate}
\end{lemma}
\begin{proof}
For $x\in Q$ let $P_x$ and $I_x$ be the corresponding projective and injective object in
$\Cscr$. It is well known and easy to see that sending $\tau^{-n} P_x$ to $(n,x)$ defines an
injective morphism of quivers $\Pi\r \NN Q^{\op}$, where $\Pi$ is the preprojective component of $\Cscr$.
It is clear that this is surjective if and only if $\Pi$ contains no injectives. So it is
sufficient to show that $\Pi$ contains no injective objects.

Assume to the contrary that $\Pi$ has some injective object. The projectives are
given by the connected quiver $Q$, and the arrows of $Q$ correspond to irreducible maps between indecomposable
projectives, in the opposite direction. Because the
categories of injectives and projectives are equivalent, the same
holds for the injectives. Thus $\Pi$ contains all indecomposable
injectives.  In particular $\Pi$ is also the preinjective
component of $\Cscr$, and hence by lemma
  \ref{ref:II.3.2-46} all objects in $\Pi$ have finite length.

Assume  that there is some finite connected subquiver $Q_1$ of $Q$ which is not Dynkin. Choose
$x\in Q_1$ and let $Q_2\supset Q_1$ be a finite connected subquiver of $Q$ such that the
simple composition factors of all objects in any path from $P_x$ to $I_x$ correspond to
vertices in $Q_2$. In the subcategory $\rep(Q_2)$ of $\Cscr$ we clearly have that $P_x$ and $I_x$ are
still projective and injective objects respectively, and it is not hard to see that there is still a path of irreducible
maps from $P_x$ to $I_x$. But since $Q_2$ is not Dynkin, it is well known that the
preprojective component of $\rep Q_2$ has no injective objects, and we have a contradiction.

If all finite subquivers of $Q$ are Dynkin, then it is easy to see directly that $Q$ must be of type
$A_\infty$, $A^{\infty}_{\infty}$ or $D_\infty$. These cases can be taken care of by direct
computations (see III. 3).
\end{proof}
In order to prove the main result in this section we shall use the
principle of tilting with respect to torsion pairs \cite{HRS}.

We
recall the main features in the special case we need. If $\Hscr$
is a hereditary abelian category, then a \emph{split torsion pair}
$(\Tscr,\Fscr)$ is a pair of additive subcategories of
 $\Hscr$ satisfying
$\Hom(\Tscr,\Fscr)=\Ext^1(\Fscr,\Tscr)=0$, and having in addition the property that for every
$A\in\Hscr$ there exist $T\in\Tscr$ and $F\in\Fscr$ such that $A\cong T\oplus F$.

The theory of $t$-structures \cite{BBD} is used to obtain  inside $D^b(\Hscr)$  a hereditary
abelian subcategory $\Hscr'$ with split torsion pair $(\Fscr[1],\Tscr)$. This construction of
$\Hscr'$ from $\Hscr$ is called tilting with respect to the torsion pair $(\Tscr,\Fscr)$ (see
\cite{HRS}). Furthermore in this case we have $D^b(\Hscr)\cong D^b(\Hscr')$. This can be seen
by combining \cite[Appendix]{Beilinson1} (which defines a functor $D^b(\Hscr)\r D^b(\Hscr')$)
with \cite[p. 13]{HRS} (which shows that this functor is fully faithful).

It is easily seen that if $\Hscr$ is a Krull-Schmidt category, then in order to specify
$(\Tscr,\Fscr)$ it is necessary and sufficient to give a partition of the indecomposable
objects in $\Hscr$: $(T_i)_{i\in I}, (F_j)_{j\in
  J}$ with properties $\Hom(T_i,F_j)=\Ext^1(F_j,T_i)=0$. Then $\Tscr$
and $\Fscr$ are respectively  the additive categories with
indecomposable objects $(T_i)_i$ and $(F_j)_j$.
Note that in the presence of Serre duality the condition
$\Ext^1(F_j,T_i)=0$ is often automatic.

\begin{figure}
\psfrag{t}[][]{$\tau$} \psfrag{ac}[][]{Other components}
\psfrag{tt}[][]{$\Tscr[-1]$} \psfrag{ff}[][]{$\Fscr$}
\psfrag{ff1}[][]{$\Fscr$} \psfrag{qq}[][]{$Q$}
\psfrag{ss}[][]{Joined components ($\Sigma$)}
\begin{center}
\includegraphics[height=5cm]{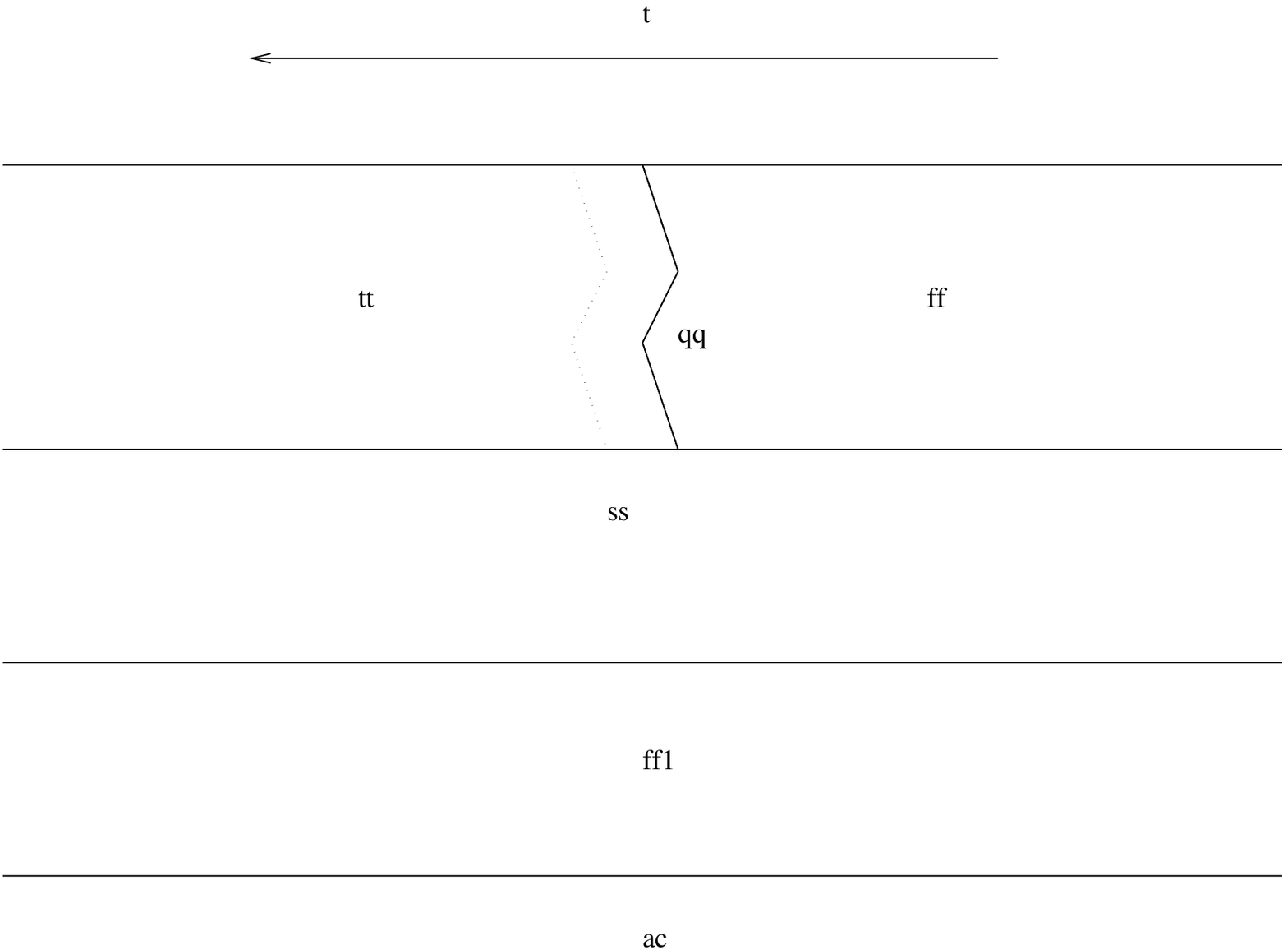}
\end{center}
\caption{The AR-quiver of $\Cscr_1$.} \label{ref:II.3.2-48}
\end{figure}

\begin{figure}
\psfrag{t}[][]{$\tau$} \psfrag{ac}[][]{Other components}
\psfrag{tt}[][]{$\Fscr_1$} \psfrag{ff}[][]{$\Tscr_1$}
\psfrag{ff1}[][]{$\Tscr_1$} \psfrag{qq}[][]{$Q'$}
\psfrag{ss}[][]{Joined components ($\Sigma$)}
\begin{center}
\includegraphics[height=5cm]{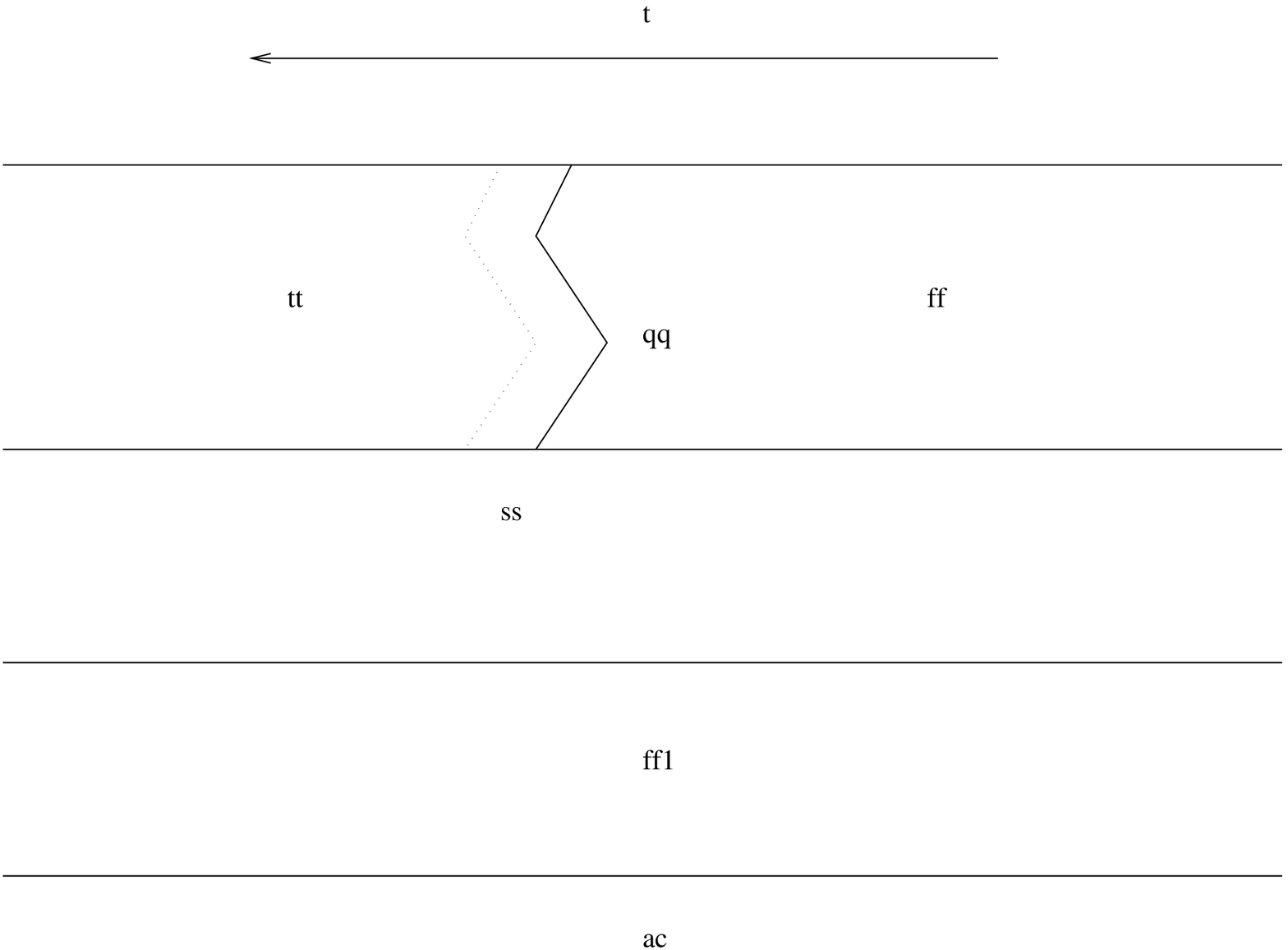}
\end{center}
\caption{A new torsion pair in $\Cscr_1$.} \label{ref:II.3.3-49}
\end{figure}

We now identify appropriate split torsion pairs in our hereditary abelian categories.
\begin{lemma}\label{lemII.3.4ny}
Let $\Cscr$ be a hereditary abelian $\Ext$-finite category with a Serre functor whose
projectives are given by a connected quiver $Q$ which is not Dynkin but satisfies (P1)(P2).
Let $Q'\subset \ZZ Q$ be a section.
\begin{enumerate}
\item
Let $\Tscr$ be the additive category generated by  the objects in
the preinjective component, and let $\Fscr$ be the additive
category whose indecomposable objects are not preinjective. Then
$(\Tscr,\Fscr)$ is a split  torsion pair in $\Cscr$.
\item
Let $\Cscr_1$ denote the hereditary abelian subcategory of
$D^b(\Cscr)$ obtained by tilting with respect to $(\Tscr,\Fscr)$
and shifting one place to the right (thus somewhat informally:
$\Cscr_1=\Tscr[-1]\oplus \Fscr$). The category $\Cscr_1$ has a
component  $\Sigma$ of type $\ZZ Q^{\op}$
  put together from the preprojective and preinjective component in
  $\Cscr$ (see Figure \ref{ref:II.3.2-48}).  In
particular $\Cscr_1$ is derived equivalent to $\Cscr$ and has no
non-zero projectives or  injectives.

Consider now the section $Q'$ in $\ZZ Q$ (see Figure \ref{ref:II.3.3-49}). Let $\Fscr_1$ be
the additive category generated by the indecomposable objects of the form $\tau^i C$ for $C\in
Q'$ and $i>0$ and let $\Tscr_1$ be the additive category generated by the other indecomposable
objects in $\Cscr_1$. Then $(\Tscr_1,\Fscr_1)$ is a split  torsion pair in $\Cscr_1$.
\item Let $\Dscr$ be the hereditary abelian category obtained by
  tilting $\Cscr_1$ with respect to $(\Tscr_1,\Fscr_1)$. Then $\Dscr$
  is derived equivalent to $\Cscr$ and the
  projectives in $\Dscr$ are given by $Q'$.
\end{enumerate}
\end{lemma}
\begin{figure}
\psfrag{t}[][]{$\tau$} \psfrag{pi}[][]{Preinjectives}
\psfrag{pp}[][]{Preprojectives} \psfrag{ac}[][]{Other components}
\psfrag{tt}[][]{$\Tscr$} \psfrag{ff}[][]{$\Fscr$}
\psfrag{ff1}[][]{$\Fscr$} \psfrag{qq}[][]{$Q$}
\begin{center}
\includegraphics[height=5cm]{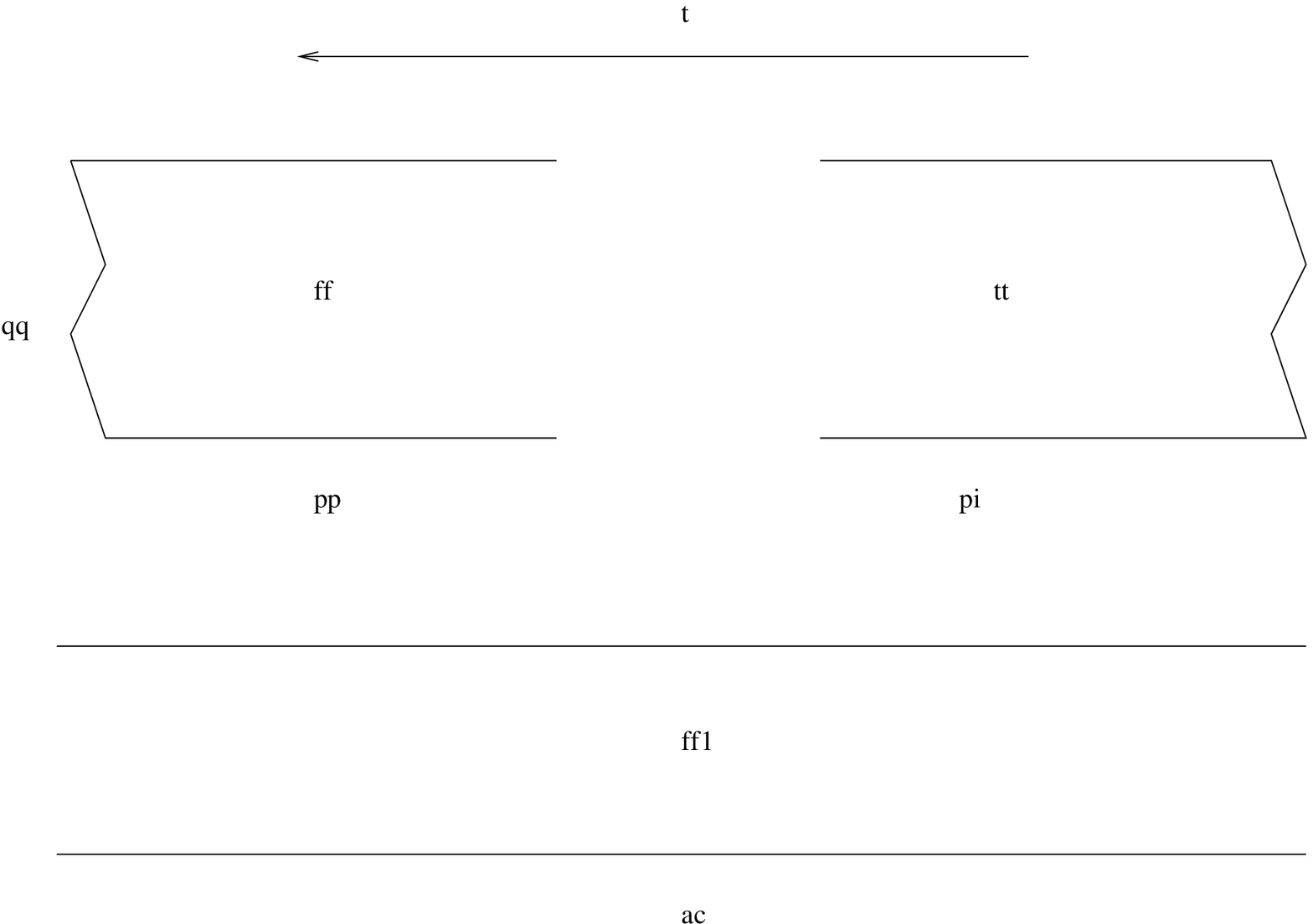}
\end{center}
\caption{The AR-quiver of $\Cscr$.} \label{ref:II.3.1-47}
\end{figure}
\begin{figure}
\psfrag{t}[][]{$\tau$} \psfrag{pi}[][]{Preinjectives}
\psfrag{pp}[][]{Preprojectives} \psfrag{ac}[][]{Other components}
\psfrag{tt}[][]{$\Fscr_1[1]$} \psfrag{ff}[][]{$\Tscr_1$}
\psfrag{ff1}[][]{$\Tscr_1$} \psfrag{qq}[][]{$Q'$}
\begin{center}
\includegraphics[height=5cm]{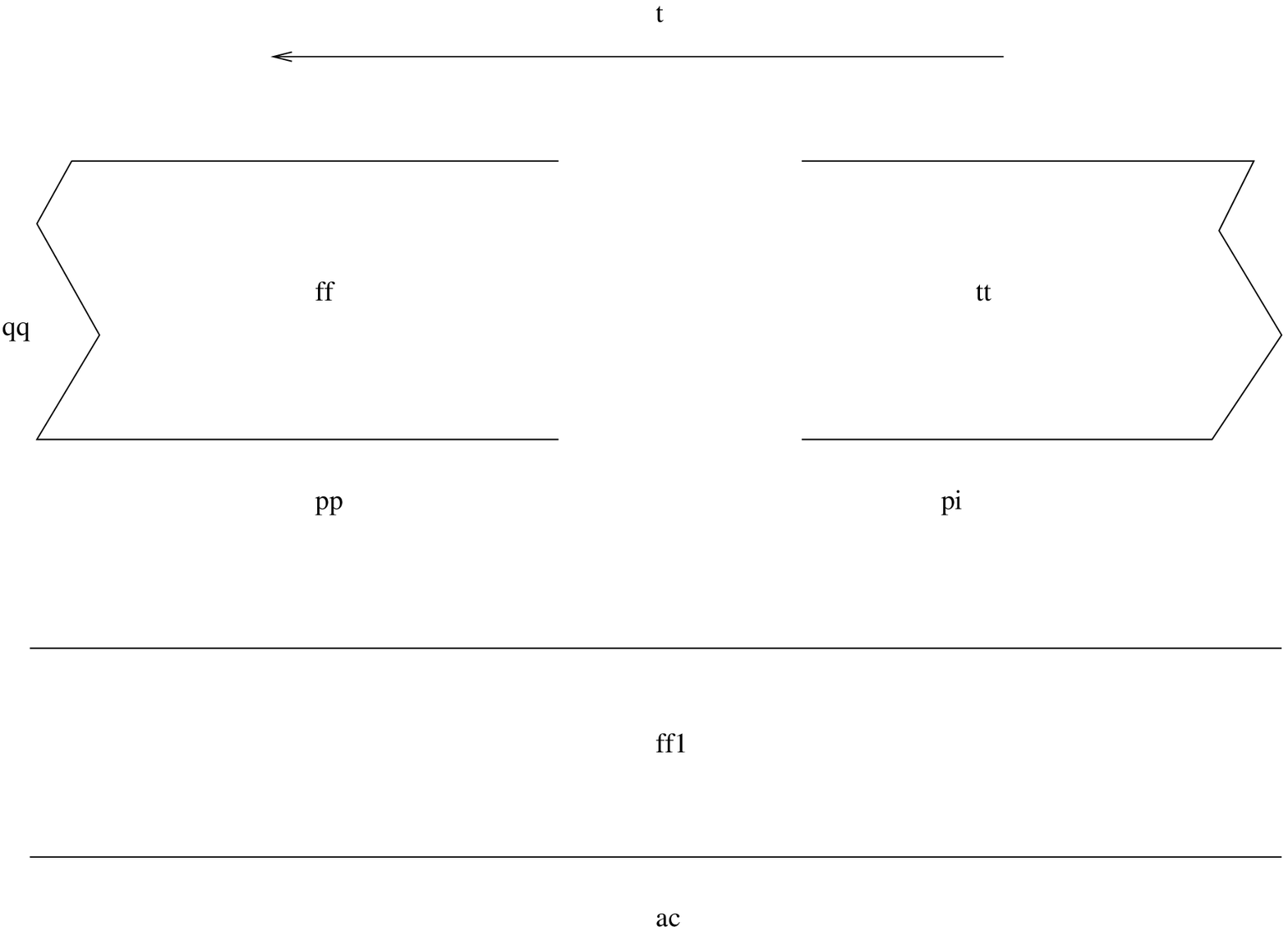}
\end{center}
\caption{The AR-quiver of $\Dscr$.} \label{ref:II.3.4-50}
\end{figure}
\begin{proof}
The proof amounts to verifying the conditions for a
  split torsion pair. In order to visualize this we have included some
  pictures of the various AR-quivers involved. Let the functor $\tau$ be
  defined as usual.
\begin{enumerate}
\item By the dual version of lemma \ref{ref:II.2.5-36} we find
  $\Hom(\Tscr,\Fscr)=0$. If  $T$ is a non-injective
  indecomposable object in $\Tscr$, then have that $\tau^{-1}T\in\Tscr$. Hence it
  follows by Serre duality that $\Ext^1(\Fscr,\Tscr)=0$.  Thus
  $(\Tscr,\Fscr)$ is a split torsion pair in $\Cscr$.
\item Now we need to show
  $\Hom(\Tscr_1,\Fscr_1)=0$.

  First let $B$ be  indecomposable in $\Sigma$ and let $A$ be an
  indecomposable which is not in $\Sigma$ (and hence in particular
  $A\in\Tscr_1$). Then looking at Figure
  \ref{ref:II.3.2-48} we see that either $B\in \Tscr[-1]$ or $B\in
  \Fscr$. In the first case it is clear that $\Hom(A,B)=0$ and in the
  second case this follows from the fact that $(\Tscr,\Fscr)$ is a
  split torsion pair.

Let now $A$ and $B$ be in $\Sigma$. By applying $\tau^{-n}$ for $n$ large enough so that
$\tau^{-n} A$ and $\tau^{-n} B$ are in the preprojective component of $\Cscr$, it is not hard
to see that if $\Hom(A,B)\neq 0$ then there is a path from $A$ to $B$ in the AR-quiver of
$\Cscr_1$. If $\Hom(A,B)\neq 0$ for some $A\in \Tscr_1$ and $B\in\Fscr_1$, we can assume that
there is an arrow from $A$ to $B$. Write $A=\tau^i X$ and $B=\tau^j Y$ where $X$ and $Y$ are
in $Q'$ and $i\le 0$ and $j>0$. Then we have an arrow $\tau^{-i} A=X\r \tau^{j-i} Y$. Since
$Q'$ is a section in $\ZZ Q$ we have $\tau^{j-i} Y\in Q'$ or $\tau^{j-i+1} Y\in Q'$. Since
$Y\in Q'$ we conclude that $j=i$ or $j=i-1$, so $j\le 0$ which is a contradiction.

Since $\Tscr_1$ is $\tau^{-1}$ stable, the fact that
$\Ext^1(\Fscr_1,\Tscr_1)=0$ follows from Serre duality.
\item  That
  the projectives in $\Dscr$ are given by $Q'$ follows from the
  shape of the AR-quiver of $\Dscr$. \qed
\end{enumerate}
\def\qed{}\end{proof}
We shall also need the following.
\begin{proposition}
\label{ref:II.3.5-51}
 Let $Q$ and $Q'$ be non-Dynkin connected quivers satisfying
  (P1)(P2) and assume that $Q'$ is a section in $\ZZ Q$. Let $\Cscr$
  be a hereditary $\Ext$-finite abelian category with a
  Serre functor whose projectives are given by the quiver $Q$, and
  assume that $\Cscr$ is generated
 by preprojectives.  Then the
  category $\Dscr$ constructed in Lemma \ref{lemII.3.4ny} is also generated by
  preprojectives.
\end{proposition}
\begin{proof}
Let $C$ be an indecomposable object in $\Dscr$. We have
  to show that it is a quotient of a preprojective object. If $C$ is
  itself preprojective, this is clear, and if $C$ is preinjective then
  we can invoke Lemma \ref{ref:II.3.2-46}.

So assume that $C$ is neither in the preinjective nor
preprojective component. Looking at Figures
\ref{ref:II.3.1-47}-\ref{ref:II.3.4-50} we see that $C$ is untouched by
the tilting, so we may consider $C$ as an object in $\Cscr$. Then
by assumption there is an exact sequence in $\Cscr$ of the form
\[
0\r K\r P\r C\r 0
\]
where $P$, and consequently $K$, is preprojective. Since all terms in the exact sequence lie
in the tilted category $\Cscr_1$, this is also an exact sequence in $\Cscr_1$. Let the functor
$\tau$ be defined as usual. Choose $n>0$ large enough such that $\tau^{-n}K$ and $\tau^{-n}P$
are preprojective in $\Dscr$, and consider the exact sequence
\[
0\r \tau^{-n}K\r \tau^{-n} P\r \tau^{-n} C\r 0.
\]
Thus $\tau^{-n}C$ is in the subcategory $\Qscr$ of $\Dscr$ generated by the preprojectives.
Since the Serre functor on $\Dscr$ restricts to one on $\Qscr$ by Theorem \ref{ref:II.2.9-42},
it follows that $C$ is also in $D^b(\Qscr)$ and hence in $\Qscr$.
\end{proof}

We now obtain the main result of this section.
\begin{theorem}
\label{ref:II.3.6-52}
 Let $Q$ and $Q'$ be connected quivers satisfying (P1)(P2),
  and assume that $Q'$ is a section in $\ZZ Q$. Then the categories
  $\wrep(Q)$ and $\wrep(Q')$ are derived equivalent.
\end{theorem}
\begin{proof}
If $Q$ and $Q'$ are Dynkin, this is well-known, so we assume that
$Q$, $Q'$ are non-Dynkin.

It follows from lemma \ref{ref:II.2.7-38} and Theorem \ref{ref:II.1.3-29} that $\wrep(Q)$ is
generated by preprojectives. By Proposition \ref{ref:II.3.5-51} the category $\Dscr$ whose
projectives are given by $Q'$, and which is obtained by a sequence of the two tilts described
above, also has the property that it is generated by preprojectives. Then it follows from
Theorem \ref{ref:II.1.3-29} that $\Dscr$ is equivalent to $\wrep(Q')$.  Since we already know
that $\wrep(Q)$ and $\Dscr$ are derived equivalent, this finishes the proof.
\end{proof}
The results in this section can be used to give an alternative approach to the construction of
the category $\wrep(Q)$ associated with the quiver $Q$ satisfying (P1)(P2), under the
following assumption:
\begin{itemize}
\item[(*)] There is a section $Q'\subset \ZZ Q$  with orientation such
  that there are no infinite paths in $Q'$.
\end{itemize}
Then we have $\wrep(Q')=\rep(Q')$. The category $\Dscr$ obtained by two tilts is also
generated by preprojectives (we also have that $Q\subset \ZZ Q'$ is a section)  by Proposition
\ref{ref:II.3.5-51}, and hence is uniquely determined by $Q$.

We have no example where (*) is \emph{not} satisfied by a quiver $Q$ satisfying (P1)(P2),
hence it is possible that the construction of $\wrep(Q)$ given in this section has the same
generality as the one given in Theorem \ref{ref:II.1.3-29}. In any case, as follows from the
results in the next section, it is general enough for the classification in the noetherian
case.

\section{The classification}
\label{ref:II.4-53} In this section we show that a connected
noetherian hereditary abelian $\Ext$-finite category with Serre
functor and nonzero projective objects is generated by the
preprojective objects, and hence is of the form $\wrep(Q)$ for
some quiver $Q$. Then we describe the quivers $Q$ such that
$\wrep(Q)$ is noetherian, completing our classification when there
are nonzero projective objects.

Our first step is  a decomposition theorem. We start with $\Cscr$ not necessarily connected.
As before let $\Qscr$ be the subcategory of $\Cscr$ generated by the  preprojectives. Let
$\Rscr$ be the full subcategory of $\Cscr$ consisting of objects $R$ such that $\Hom(Q,R)=0$
for any preprojective $Q$ (or equivalently for any $Q\in \Qscr$).
\begin{lemma}
\label{ref:II.4.1-54} $\Rscr$ is closed under subquotients and extensions.
\end{lemma}
\begin{proof}
 The only thing that is not entirely trivial is the fact that
  $\Rscr$ is closed under quotients. So let $\alpha:R\r S$ be a
  surjective map in $\Cscr$ with $R\in\Rscr$. Let $\theta:Q\r S$ be a
  non-trivial map with $Q$ preprojective. Finally let
\[
\begin{CD}
R @>\alpha>> S\\ @A \psi AA @A \theta AA\\ Z @>\beta >> Q
\end{CD}
\]
be the corresponding pullback diagram.  Since $\theta\beta$ is not
the zero map, there must exist  an indecomposable summand $Z_1$ of
$Z$ such that $\theta\beta(Z_1)\neq 0$. Since $Z_1$ maps
non-trivially to $Q$, it follows from Corollary
\ref{ref:II.2.6-37} that it must be preprojective. Now clearly
$\psi(Z_1)\neq 0$. This yields a contradiction.
\end{proof}
Now we come to our main result in this section. Let $\Cscr$ be an $\Ext$-finite
noetherian hereditary abelian category, and let the notation be as above.
\begin{theorem}
\label{ref:II.4.2-55}
\begin{enumerate}
\item
The inclusions $\Rscr\r \Cscr$ and $\Qscr\r \Cscr$ define
  an equivalence $\Rscr\oplus\Qscr\cong \Cscr$.
\item The hereditary abelian categories $\Qscr$ and $\Rscr$ satisfy Serre
  duality.  If the functor $\tau$ is defined as usual, then $\tau$ is
  everywhere defined on $\Rscr$ and is invertible.
\end{enumerate}
\end{theorem}
\begin{proof}

\begin{enumerate}
\item
We have by definition $\Hom(\Qscr,\Rscr)=0$. Furthermore it follows from lemma
\ref{ref:II.4.1-54} and Theorem \ref{ref:II.2.9-42} that $\Hom(\Rscr,\Qscr)=0$. Hence it is
sufficient to show that every $C\in \Cscr$ is  a direct sum $R\oplus Q$ with $R\in\Rscr$ and
$Q\in\Qscr$.

{}{}From the fact that $\Cscr$ is noetherian and $\Qscr$ is closed under extensions  (Theorem
\ref{ref:II.2.9-42}) it follows that there exists a maximal subobject $Q$ in $C$ which lies in
$\Qscr$, and so $R=C/Q$ lies in $\Rscr$. Hence it now suffices to prove that $\Ext^1(R,Q)=0$.
Since $\Cscr$ is hereditary, $\Ext^1(R,-)$ preserves epis, thus it is sufficient to show that
$\Ext^1(R,\tau^{-a}P)=0$ with $P$ an indecomposable projective and $a\ge 0$. In addition we
may and we will assume that $R$ is indecomposable.

Assume $\Ext^1(R,\tau^{-a}P)\neq 0$. Then clearly $R$ is not projective and $\tau^{-a}P$ is
not injective. Hence by Serre duality $\Hom(\tau^{-a-1}P,R)\neq 0$. This contradicts the fact
that $R\in\Rscr$. Hence $\Ext^1(R,Q)=0$, which finishes the proof.
\item That $\Qscr$ and $\Rscr$ satisfy Serre duality follows for
  example from lemma \ref{ref:I.1.11-10}. Any projective or
  injective object in $\Rscr$ clearly has the same property in
  $\Cscr=\Rscr\oplus \Qscr$ and therefore lies in $\Qscr$. Since a non-zero
  object cannot be both in $\Rscr$ and $\Qscr$, we conclude that
  $\Rscr$ does not contain nonzero projectives or injectives. Hence $\tau$ is defined everywhere on
  $\Rscr$ and is invertible.
\end{enumerate}
\end{proof}

For an example showing that this theorem fails for non-noetherian
hereditary abelian categories (even when they have noetherian injectives),
see \S\ref{ref:III.3.6-105}.

Theorem \ref{ref:II.4.2-55} together with Corollary
\ref{ref:II.2.8-41} tells us that $\Cscr$ is of the form
$\wrep(Q)\oplus \Rscr$ where $\Rscr$ has a Serre functor. When
$\Cscr$ is connected, it follows that $\Cscr = \wrep(Q)$. It
remains to describe the $\wrep(Q)$ which are noetherian.

Let us first define a particular kind of quivers. We call a \emph{ray} an infinite quiver of
the following form
\[
x:\quad x_0\r x_1\r x_2\r\cdots
\]
We call $x_0$ the \emph{starting vertex}.

We call a quiver $Q$ \emph{strongly locally finite} if it is connected and if for every vertex in $Q$ the
associated projective and injective representation have finite length.
Note that a strongly locally finite
quiver automatically satisfies (P1)(P2).

If $Q$ is a quiver, then by attaching a ray we mean identifying a vertex of $Q$ with the
starting vertex of the ray.

Finally we call a quiver a \emph{star} if it consists of a strongly locally finite quiver $Q_0$ to which  a
set of rays is simultaneously attached in such a way that  only a finite number of rays are
attached to every vertex. Note that a star is connected and satisfies (P1)(P2).

We  now prove the following result.
\begin{theorem}
\label{ref:II.4.3-56}
 Let $Q$ be a connected quiver satisfying the properties
  (P1)(P2). Then $\wrep(Q)$ is noetherian if and only if
  $Q$ is a star.
\end{theorem}
The proof will consist of a number of lemmas.

We will be analyzing  \emph{infinite paths} in $Q$. These are by
definition  subquivers of $Q$ of the form
\[
x:\quad x_0\r x_1\r x_2\r\cdots
\]
(so actually infinite paths are the same as rays but we use different terminology to avoid
confusion).

Two infinite  paths are said to be equivalent if they are equal from some vertex on.
\begin{lemma} Assume that $\rep(Q)$ is noetherian.
 If $x$ is an infinite  path in $Q$, then there
  exists $n_0$ such that for $n\ge n_0$ the only arrow in $Q$ starting
 in  $x_n$ is the arrow in $x$ going from $x_{n}$ to $x_{n+1}$.
\end{lemma}
\begin{proof}
Assume that $n_0$ does not exist. We will show that then $\rep(Q)$ is not noetherian. The
non-existence of $n_0$ implies that there exists a proper strictly ascending array of integers
$(n_i)_{i\ge 1}$ as well a corresponding array of vertices $(y_i)_i$ such that there is an
arrow $e_i:x_{n_i}\r y_{i}$ which is not in $x$.  The edges $e_i$ induce maps
$\phi_i:P_{y_i}\r P_{x_0}$ which are the composition of the maps $P_{y_i}\xrightarrow{e_i}
P_{x_{n_i}}\r P_{x_0}$.

We now claim that the image of $(\phi_i)_{i\le
  m}:\oplus_{i\le m} P_{y_i}\r P_{x_0}$ defines a proper ascending chain of
  subobjects of $P_{x_0}$. This is clear since the path in $Q$
  corresponding to the map $\phi_{m+1}$ in $Q$ does not pass through
  $e_i$ for $i\le m$. Hence $\rep(Q)$ is not noetherian, and we have
  obtained a contradiction.
\end{proof}
In order to give a more precise analysis when $\wrep(Q)$ is noetherian we need the following.

\begin{lemma}
\label{ref:II.4.5-57} Let $\Cscr$ be a $k$-linear $\Ext$-finite hereditary abelian category
with a Serre functor $F$. Let $\tau$ be defined as usual.

Let $Q,P_1$ and $P_2$ be non-injective indecomposable projectives in $\Cscr$. Assume we have
an irreducible map $\alpha :Q\r \tau^{-1} P_1$ and a map $\beta :P_1\r P_2$ which does not
factor through $Q$. Let $\varphi$ be the composition of maps
\[
Q\xrightarrow{\alpha} \tau^{-1}P_1\xrightarrow{\tau^{-1}\beta} \tau^{-1}P_2
\]

Then  for any factorization of $\phi$ through a projective object $X$
\[
Q\xrightarrow{\gamma} X\xrightarrow{\psi} \tau^{-1}P_2
\]
we have that  $\gamma$ splits.
\end{lemma}
\begin{proof}
Assume that $\gamma$ is not split. Consider the almost split
sequence
\[
0\r
Q\xrightarrow{(\delta_i)_i}\bigoplus_{i=1,\ldots,m}Y_i\xrightarrow{(\epsilon_i)_i}
 \tau^{-1}Q\r 0
\]
{}{}From the theory of almost split sequences it follows that we may take $\delta_1=\alpha$ and
$Y_1=\tau^{-1}P_1$. Also by the theory of almost split sequences it follows that
$\gamma=\sum_i\mu_i\delta_i$ for some maps $\mu_i:Y_i\r X$. Since $X$ is projective we have
$\Hom(\tau^{-1}P_1,X)=0$, and hence $\mu_1=0$. Now consider the map $\theta:\bigoplus_i Y_i \r
\tau^{-1}P_2$ given by
\[
(-\tau^{-1}\beta,\psi\mu_2,\ldots,\psi\mu_m)
\]
Clearly $\sum \theta_i\delta_i=0$, whence there exists $\omega:\tau^{-1}Q\r \tau^{-1}P_2$ such
that
\begin{align*}
 -\tau^{-1}\beta&=\omega\epsilon_1\\
\psi\mu_i&=\omega\epsilon_i\qquad (i\ge 2)
\end{align*}
Applying $\tau$ to the first of these equations yields a factorization of $\beta$,
contradicting the hypotheses.
\end{proof}

\begin{lemma} Assume that $\wrep(Q)$ is noetherian. If $x$ is an
  infinite path in $Q$, then there exists $n_0$ such that for $n\ge
  n_0$ the only arrows in $Q$ adjacent to $x_n$ are those in $x$.
\end{lemma}
\begin{proof} If $\wrep(Q)$ is noetherian, then so is $\rep(Q)$ since
  the latter
  is a full subcategory closed under subobjects. Therefore it follows
  from the previous lemma that there is some $n'_0$ such that for $n\ge
  n'_0$ the only arrow starting in $x_n$ is the one which lies in
  $x$. By dropping some initial vertices in $x$ we may assume $n'_0=0$.

Assume now that $n_0$ (as in the statement of the lemma) does not exist. Since injectives have
finite length, it is clear that $P_{x_0}$ is not injective. We will show that
$\tau^{-1}P_{x_0}$ is not noetherian in $\wrep(Q)$.

The non-existence of $n_0$ implies that  there exists a strictly
ascending array of integers $(n_i)_{i\ge 1}$ as well a
corresponding array of vertices $(y_i)_i$ such that there is an
arrow $e_i:y_i\r x_{n_i}$ which is not in $x$.

The edges $e_i$ correspond to an irreducible map $P_{x_{n_i}}\r P_{y_i}$. Using the theory of
almost split sequences there is a corresponding map $e'_i:P_{y_i}\r \tau^{-1}P_{x_{n_i}}$.
Also from the path $x$ we obtain a map $P_{x_{n_i}}\r P_{x_0}$ which gives rise to a
corresponding map $\tau^{-1}P_{x_{n_i}}\r \tau^{-1}P_{x_0}$. We denote by $\phi_i$ the
composition
\begin{equation}
\label{ref:II.4.1-58} P_{y_i}\xrightarrow{e'_i} \tau^{-1} P_{x_i}\r \tau^{-1} P_{x_0}.
\end{equation}

We now claim that the image of $(\phi_i)_{i\le m}:\oplus_{i\le m} P_{y_i}\r \tau^{-1}P_{x_0}$
defines an ascending chain of subobjects of $\tau^{-1}P_{x_0}$. If this were not the case,
then for some $m$ all $\phi_{m'}$ and $m'> m$ must factor through $\bigoplus_{i\le m}P_{y_i}$.
It follows that the resulting map $P_{y_{m'}}\r \bigoplus_{i\le
  m}P_{y_i}$ must be split since otherwise by lemma \ref{ref:II.4.5-57}
and \eqref{ref:II.4.1-58} the map $P_{x_{m'}}\r P_{x_0}$
factors through $P_{y_{m'}}$, which is impossible given the
construction of this map.

Thus it follows that $y_{m'}\in\{y_i\mid i\le m\}$. We conclude
that some vertex $y$  occurs infinitely often among the $y_i$, but
this contradicts condition (P1).

Thus we have shown that $\wrep(Q)$ is not noetherian, and in this
way we have obtained a contradiction to the hypotheses.
\end{proof}
\begin{corollary}
Assume that $\wrep(Q)$ is noetherian. Then $Q$ is a star.
\end{corollary}
\begin{proof}
Let $\Omega$ be the set of equivalence classes of infinite paths
in $Q$. For every $\omega\in\Omega$ we choose a representative
$x_\omega$. By the previous lemma we can without loss of
generality assume that the only edges in $Q$ adjacent to
$x_{\omega,n}$ ($n\ge 1$) are those in $x_\omega$.

Now let $Q_0$ be obtained from $Q$ by removing for all $\omega\in\Omega$ the vertices
$x_{\omega,n}$ for $n\geq 1$ as well as the edges in $x_\omega$. It is clear that $Q$ is
obtained from $Q_0$ by adjoining the strings $(x_\omega)_{\omega\in\Omega}$. Furthermore $Q_0$
itself cannot contain any infinite paths since such an infinite path would have to be
equivalent to one of the $x_\omega$. In particular it would have to contain vertices outside
$Q_0$, which is of course a contradiction. It follows that $Q_0$ is
strongly locally finite, and so $Q$ is
indeed a star.
\end{proof}
The following lemma completes the proof of Theorem
\ref{ref:II.4.3-56}.
\begin{lemma} If $Q$ is a star, then $\wrep(Q)$ is noetherian.
\end{lemma}
\begin{proof}
  Assume that $Q$ is a star. Since by construction $\wrep(Q)$ is generated by preprojectives, it
  suffices to show that the indecomposable preprojectives are
  noetherian. So we need to consider objects of the form $\tau^{-b}P_x$
  with $x$ a vertex in $Q$ and $b\geq 0$.

  Taking the sum of all irreducible maps (with multiplicities) going
  into $\tau^{-b}P_x$, we obtain a map $\oplus_i P_i\r \tau^{-b}P_x$, where the sum is finite, which is
  either surjective or has simple cokernel (the latter happens if
  $b=0$). In any case it is sufficient to show that the $P_i$ are
  noetherian.

  By lemma \ref{ref:II.2.5-36} we know what the $P_i$ can be.
They are
  of the form $\tau^{-c}P_y$ such   that either $c<b$, or else $c=b$ and the map $\tau^{-b}P_y\r \tau^{-b}P_x$
  is obtained from an irreducible map $P_y\r P_x$, which in turn corresponds
  to an arrow $x\r y$ in the quiver $Q$.

  Induction on $b$ now yields that it is sufficient to show that
  $\tau^{-b}P_z$ is noetherian, where $z$ lies on one of the rays
  contained in $Q$, and furthermore for a given $b$ we may assume that
  $z$ lies arbitrarily far from the starting vertex of the ray.

Let $w$ be the ray on which $z$ lies. We will assume that $z=w_{b}$ (if $z$ happens to be
farther away, then we just drop the initial vertices in $w$.)

To check that $\tau^{-b}P_z$ is noetherian we need to understand the additive category $\Uscr$
whose objects are direct sums of indecomposables which have a non-zero map to $\tau^{-b}P_z$.
The theory of almost split sequence easily yields that $\Uscr$ is the path category of the
part of the AR-quiver of $\wrep(Q)$ spanned by the vertices that have a non-trivial path to
$\tau^{-b}P_z$. Furthermore, also by the theory of almost split sequences, it is easy to work out
what this quiver is. The result is as in Figure \ref{ref:II.4.1-61}. By convention we will
assume that the relations on the quiver in Figure \ref{ref:II.4.1-61} are of the form
\begin{equation}
\label{ref:II.4.2-59} e_{v,u}f_{u,v}=f_{u+1,v}e_{v+1,u}
\end{equation}
Note that from these relations, or directly, it easily follows
\begin{equation}
\label{ref:II.4.3-60} \Hom(\tau^{-c}P_{w_m},\tau^{-c'}P_{w_{m'}})=
\begin{cases}
k&\text{if $c\le c'$ and $m\ge m'$}\\ 0&\text{otherwise}
\end{cases}
\end{equation}
\begin{figure}
\begin{center}
$
\xymatrix{ P_{w_0}\ar[r]^-{e_{00}} & \tau^{-1}P_{w_1}
\ar[r]^-{e_{01}} &\tau^{-2}P_{w_2}\ar@{.}[r] &
\tau^{-b+1}P_{w_{b-1}}\ar[r]^-{e_{0,b-1}} &
\tau^{-b}P_{w_b}\\ %end of first row
P_{w_1}\ar[r]^-{e_{10}} \ar[u]^-{f_{00}}\ar@{--}[ur]&
\tau^{-1}P_{w_2}\ar[u]^-{f_{10}} \ar[r]^-{e_{11}}\ar@{--}[ur]
&\tau^{-2}P_{w_3}\ar[u]^-{f_{20}}\ar@{.}[r] &
\tau^{-b+1}P_{w_b}\ar[u]^-{f_{b-1,0}}\ar@{--}[ur]\ar[r]^-{e_{1,b-1}}
&
\tau^{-b}P_{w_{b+1}}\ar[u]^-{f_{b0}}\\ %end of second row
P_{w_2}\ar[r]^-{e_{20}} \ar[u]^-{f_{01}}\ar@{--}[ur]&
\tau^{-1}P_{w_3}\ar[u]^-{f_{11}} \ar[r]^-{e_{21}}\ar@{--}[ur]
&\tau^{-2}P_{w_4}\ar[u]^-{f_{21}}\ar@{.}[r] &
\tau^{-b+1}P_{w_{b+1}}\ar[u]^-{f_{b-1,1}}\ar@{--}[ur]\ar[r]^-{e_{2,b-1}}
& \tau^{-b}P_{w_{b+2}}\ar[u]^-{f_{b1}}
\\ %end of third row
\ar@{.}[u]& \ar@{.}[u] &\ar@{.}[u]&\ar@{.}[u]&\ar@{.}[u] %end of fourth row
}
$
\end{center}
\caption{Part of the preprojective component of $\wrep(Q)$}
\label{ref:II.4.1-61}
\end{figure}
To show that $\tau^{-b}P_{w_b}$ is noetherian we have to show that there does not exist an
infinite sequence of non-zero maps (unique up to a scalar by \eqref{ref:II.4.3-60})
$\phi_i:\tau^{-c_i}P_{w_{n_i}} \r \tau^{-b} P_{w_b}$, such that the image of $\phi_{m+1}$ is
not contained in the image of $\bigoplus_{i\le m}\phi_i$.

Assume to the contrary that such a sequence $(\phi_i)_i$ does
indeed exist. Let $C\subset\{0,\ldots,b\}$ be the set of all $c$
such that $c=c_i$ for some $i$. For every $c\in C$ let $v_c$ be
the  $w_{n_i}$ closest to $w_0$ with the property that $c=c_i$ and let
$\theta_c$ be the corresponding $\phi_i$.  It is now clear from
the quiver given in Figure \ref{ref:II.4.1-61} as well as
\eqref{ref:II.4.3-60} that every $\phi_i$ factors through one of the
$\theta_c$, contradicting the choice of $(\phi_i)_i$. This
finishes the proof.
\end{proof}
Summarizing, we have the following main result of this chapter.
\begin{theorem}\label{ref:II.4.9-62}
The connected hereditary abelian noetherian $\Ext$-finite
categories $\Cscr$ with Serre functor and nonzero projective
objects are exactly the categories $\wrep(Q)$, where $Q$ is a
connected quiver which is a star.
\end{theorem}

\chapter[Sources of hereditary abelian categories]{Sources of hereditary abelian categories with no
projectives or injectives.}\label{ref:III-63} In this chapter we give various sources of
examples of  hereditary abelian categories. The main focus is on the $\Ext$-finite categories
which are noetherian and have a Serre functor, and we shall see in the next chapter that our
discussion includes all possible examples with no nonzero projectives. It is also interesting
to see how these examples fit into more general classes of hereditary abelian categories.

In Section \ref{ref:III.1-64} we investigate hereditary abelian categories whose objects have
finite length, and classify those having a Serre functor. \VDB
\begin{comment}
In Section \ref{ref:III.2-68}
we investigate when various quotient categories of module
categories, graded or not, give rise to hereditary abelian categories, and
when they satisfy the central properties of this paper. Of
particular interest is the graded two-dimensional case, where we
give an interpretation in terms of sheaves of hereditary orders.
\end{comment}

In Section \ref{ref:III.2-68} we describe the hereditary abelian categories that arise as
$\qgr(S)$ for $S$ a graded ring finite over its center (see \ref{ref:III.2-68} for definition
of $\qgr (S)$).

In Section \ref{ref:III.3-73} we give examples of noetherian hereditary abelian categories
obtained by tilting with respect to torsion pairs, performed inside the bounded derived
category, as discussed in \ref{ref:III.2-68}.

\section[Hereditary abelian categories with Serre functor]{Hereditary abelian categories with  Serre functor
 and all objects of finite length.} \label{ref:III.1-64} When $\Cscr$ is a hereditary abelian category, the
subcategory f.l.\ $\Cscr$ whose objects are those of finite length is again hereditary
abelian, by Proposition \ref{ref:2.2a}. In this section we classify the hereditary abelian
categories $\Cscr$ having a Serre functor, and where $\Cscr =f.l. \Cscr$.

Recall that $\tilde{A}_{n-1}$ denotes the graph which is a cycle with $n$ vertices. Then we
have the following.
\begin{theorem}
\label{ref:III.1.1-65} Let $\Cscr$ be a connected
  $\Ext$-finite noetherian hereditary abelian category in which every object
  has finite length and
    which has almost split
  sequences and no nonzero projectives or injectives. Then $\Cscr$
  is equivalent to the category of nilpotent finite dimensional
  representations of the quiver $\tilde{A}_{n-1}$ or of the quiver$A_{\infty}^{\infty}$, with all arrows
  oriented in the same direction. In the first case $n$
  is the number of  simple
  objects in $\Cscr$ and in the second case there is an infinite number of simple objects.

Conversely, the category of nilpotent finite dimensional representations of the quiver
$\tilde{A}_{n-1}$ or of the quiver $A_{\infty}^{\infty}$, with all arrows oriented in the same
direction, is a noetherian hereditary abelian $\Ext$-finite category with almost split
sequences, and hence also with a Serre functor.
\end{theorem}
\begin{proof}
By Theorem \ref{ref:I.3.2-20} the category $\Cscr$ has a Serre functor. As usual we denote the
Serre functor by $F$, and we let $\tilde{\tau}=F[-1]$. This now induces an autoequivalence
$\tau : \Cscr\r \Cscr$.

 Let
$S$ be an object in $\Cscr$\, and let $\Cscr'$ be the component of the AR-quiver of $\Cscr$\
containing $S$. Since $\tau \colon \Cscr \rightarrow \Cscr$ is an equivalence, it preserves
length. Hence we have an almost split sequence $0 \rightarrow \tau S \rightarrow E \rightarrow
S \rightarrow 0$, where $E$ is uniserial of length $2$. Consider the almost split sequence  $0
\rightarrow \tau E \rightarrow \tau S \bigoplus F \rightarrow E \rightarrow 0$. Since $E$ and
$\tau E$ are uniserial of length $2$ and $\tau S$ is simple, it is easy to see that $F$ is
uniserial of length $3$. Continuing the process, we see that $\Cscr'$ contains only uniserial
objects. It is a tube if it contains only a finite number of nonisomorphic simple objects, and
of the form $\ZZ A_\infty$ otherwise.

Assume there is some indecomposable object $X$ in $\Cscr$\ which
is not in $\Cscr'$. If $\Hom(X,\Cscr')\neq 0$, then $\Hom(X,S)
\neq 0$ for some simple object $S$ in $\Cscr'$. Using the
properties of almost split sequences and that all objects in
$\Cscr'$ are uniserial, we can lift the map $f \colon X
\rightarrow S$ to get an epimorphism from $X$ to an object of
arbitrary length in $\Cscr'$. Hence $\Hom(X,\Cscr')=0$, and
similarly $\Hom(\Cscr',X)=0$. It follows that all indecomposable
objects in $\Cscr$\ are in $\Cscr'$.

It is easy to see from the above that the component of $\Cscr$
uniquely determines $\Cscr$ and that $\Cscr$ is as in the
statement of the theorem.

The category of nilpotent finite dimensional representations of $\tilde{A}_{n-1}$ or of $A^\infty_\infty$ over $k$ is a
hereditary abelian $k$-category, and it is easy to see that it has almost split sequences, and
hence a Serre functor by Theorem  \ref{ref:I.3.2-20}.
\end{proof}
Note that in \cite{SOS} the connected finite quivers $Q$ where f.l.\ $(\Rep Q)$ has almost
split sequences are described to be the above $ \tilde{A}_n$  and the quivers $Q$ with no
oriented cycles.

There are various alternative but equivalent descriptions of the category $\Cscr$. For further
reference we now give a discussion of those.

 For $n>0$  we define $C$ as the ring of $n\times
n$-matrices of the form
\[
\begin{pmatrix}
R & Rx &\cdots& Rx\\ \vdots &\ddots &\ddots &\vdots\\ \vdots &&
\ddots & Rx\\ R&\cdots & \cdots & R
\end{pmatrix}
\]
with $R=k[[x]]$.  We put the standard $(x)$-adic topology on $C$.
Since $C$ is the completed path algebra of $\tilde{A}_{n-1}$, it
follows that $\Cscr$ is equivalent to the category of finite
dimensional representations over $C$.

Similarly for $A_{\infty}^{\infty}$  let $R=k[x]$ considered as graded ring with $\deg x=1$.
It is easily seen that $\Cscr$ is now equivalent to $\tors(R)$, which consists of the finite
dimensional $R$-modules.

 In order to have compatibility with the first case
let $C$ be the ring of lower triangular $\ZZ\times \ZZ$-matrices over $k$,
equipped with the product of the discrete topologies on $k$. Thus
$C$ is a pseudo-compact ring in the sense that its topology is
complete and separated and in addition is generated by ideals of
finite colength. For details on pseudo-compact rings see
\cite{Gabriel,VdB19,VdBVG}. Since $C$ is the completed path
algebra of ${A}_\infty^\infty$ it follows that $\Cscr$ is
equivalent to the category of pseudo-compact finite dimensional
representations over $C$. Note that one can prove (as in the proof
of \cite[Thm 1.1.3]{VdBVG}) that all finite dimensional
$C$-representations equipped with the discrete topology are
pseudocompact (this depends on $|k|=\infty$).

Let $\dis(C)$ be the category of finite dimensional $C$-modules (with the discrete topology).
In the sequel we will denote the equivalence $\Cscr\r \dis(C)$ by $(-)\hat{}$. To finish the
description of $\Cscr$ we describe the functor $\tau$.
\begin{proposition}
\label{ref:III.1.2-66}
 There exists an invertible bipseudocompact \cite{VdB19}  $C$-bimodule
  $\omega_C$ such that for $T\in \Cscr$ we have
\[
\widehat{\tau T}=\hat{T}\otimes_{C} \omega_C
\]
Furthermore non-canonically we have $\omega_C=\rad(C)$.
\end{proposition}
\begin{proof} Since Serre functors are unique we may assume
  $\Cscr=\dis(C)$. It is easy to see that $\omega_C$ exists and is given by
  the following formula.
\[
\omega_C=\projlim_{I\subset C} \tau (C/I)
\]
where $I$ runs through the ideals of $C$ with the property that
$C/I$ is finite dimensional. In fact such a formula would hold for
any autoequivalence of $\Cscr$.

To describe $\omega_C$ explicitly consider first the case $n<\infty$. In that case local
duality implies that $\omega_C\cong \Hom_R(C,\omega_R)$ where $\omega_R=\Omega^1_R=R\,dx\cong
R$. An explicit computation reveals that $\omega_C\cong \rad(C)$.

In the case of $A_{\infty}^{\infty}$ we could extend local duality theory to the
pseudo-compact ring $C$ to arrive at the same formula. However it is easier to use that
$\Cscr\cong \tors(R)$ with $R=k[x]$. In this case graded local duality implies that $\tau$
coincides with tensoring with $\omega_R=\Omega^1_R=R\,dx\cong R(-1)$. Translating this to $C$
we see that indeed $\tau$ coincides with tensoring with $\rad(C)$.
\end{proof}
We denote by $\PC(C)$ the category of pseudo-compact $C$-modules.
$\PC(C)$ is an abelian category with enough projectives and exact
inverse limits. Note the following fact which will be used in the
sequel.
\begin{lemma}
\label{ref:III.1.3-67} The  homological dimension of $\PC(C)$ is
one.
\end{lemma}
\begin{proof} In \cite[Cor 4.7,4.11]{VdB19} it has been shown that it is sufficient
  to prove that the projective dimension of each pseudocompact simple
  is equal to one. This is an easy direct verification.
\end{proof}

\section{Sheaves of hereditary orders and graded rings}
\label{ref:III.2-68} Let $X$ be a non-singular projective curve over $k$. Let $K$ be the
function field of $X$ and let $A=M_n(K)$ for some $n>0$. Let $\Oscr$ be a sheaf of hereditary
$\Oscr_X$-orders in $A$. Thus locally $\Oscr$ is a hereditary order over a Dedekind ring (in
the sense of \cite{reiner}). Let $\Cscr=\coh(\Oscr)$ be the category of coherent
$\Oscr$-modules. Then $\Cscr$ is obviously hereditary. Put
$\omega=\uHom_{\Oscr_X}(\Oscr,\omega_X)$. Then exactly as in the commutative case one proves
for $\Escr,\Fscr\in \Cscr$:
\[
\Hom_{\Oscr}(\Escr,\Fscr)\cong\Ext^1_\Oscr(\Fscr,\Escr\otimes_{\Oscr}\omega)^\ast
\]
Hence $\Cscr$ satisfies Serre duality.

{}{}From the structure of hereditary orders \cite{reiner} it follows that every point $x\in X$
corresponds to a unique $-\otimes\omega$ orbit of simples (those simples with support in $x$).
Hence in particular $\Cscr$ has an infinite number of $\tau$-orbits of simples.

At some point in the proof of Theorem \ref{ref:IV.5.2-152}  we need that if $S$ is a
commutative graded ring such that $\qgr(S)$ is connected hereditary abelian, then $\qgr(S)$ is
of the form $\coh(\Oscr)$ for $\Oscr$ as above. Here $\qgr(S)$ denotes the quotient category
$\gr(S)$/ finite length. We prove this in the larger generality that $S$ is finite over its
center. Once we have the classification in Theorem \ref{ref:IV.5.2-152} available we will be
able to give a proof under even more general conditions on $S$. See Corollary
\ref{ref:V.2.2-159}.

We start with the following preliminary result.
\begin{proposition}\label{ref:III.2.1-69} Let $S$ be a noetherian
  $\ZZ$-graded $k$-algebra finite dimensional over $k$ in every degree with left
  limited grading and which is finitely
  generated as a module over a central commutative ring $C$. Let
  $X=\Proj C$ (see \cite{H} for the definition of $\Proj$). Then
  $\qgr(S)$ is equivalent to the category of coherent modules over a
  sheaf $\Oscr$ of $\Oscr_X$-algebras which is coherent as
  $\Oscr_X$-module.
\end{proposition}
\begin{proof}
Without changing the category $\qgr(S)$ we may (and we will)
replace
 $S$ by $S=k+S_1+S_2+\cdots$.
 By the Artin--Tate lemma $C$
  is finitely generated. Let $x_1,\ldots,x_n$ be homogeneous
  generators for $C$,
  respectively of degree $a_1,\ldots,a_n>0$. Let $m_1,\ldots,m_t$ be
  homogeneous generators of $S$ as a module over $C$, and let  $e$ be the
  maximum of the degrees of the $m$'s. Let $p$ be the product of the
  $a_i$'s. We claim that there exists $b_0\in \NN$ such that for $b\ge
  b_0$ we have $S_p S_b=S_{p+b}$.

  We prove instead that for $a$ large we have $S_a=S_p S_{a-p}$. This
  is clearly equivalent. A
  general element of $S_a$ is a $k$-linear combination of elements of
  the form $s=x_1^{p_1}\cdots x_n^{p_n} m_j$. It is sufficient to have
  that for some $i$~: $p_i\ge a_1\cdots \hat{a}_i\cdots a_n$. Assume
  this is false. Then it follows that the degree of $s$ (which is $a$) is
  less than $np+e$. This proves what we want.

Define
\[
T=\begin{pmatrix} S& S(1) &\cdots &S(p-1)\\ S(-1)& S &\cdots
&S(p-2)\\ \vdots &\vdots &\ddots &\vdots\\ S(-p+1)& S(-p+2)
&\cdots& S
\end{pmatrix}
\]
  One  verifies directly
  that for $b\ge p+b_0-1$ one has $T_1 T_b=T_{b+1}=T_b T_1$. Since
  $S$ and $T$ are Morita equivalent we clearly have $\qgr(S)=\qgr(T)$.

  By definition $X$ is covered by the affine open sets $U_i=\Spec
  (C_{x_i})_0$. Let $\Tscr$ be the
   sheaf of graded rings on $X$ associated to $T$. Thus on $U_i$ the
  sections of $\Tscr$ are given by $T_{x_i}$.

  By localization theory $\qgr(T)$ is equivalent to the category of
  coherent graded $\Tscr$-modules. We claim that $\Tscr$ is strongly
  graded ($\Tscr_{a}\Tscr_b=\Tscr_{a+b}$) and that $\Oscr=\Tscr_0$ is
  coherent over $\Oscr_X$.  It is clearly sufficient to check this
  locally. That $T_{x_i}$ is strongly graded follows
  from Lemma \ref{ref:III.2.2-70} below.

  The sheaf $\Oscr$ is nothing but $\tilde{T}$ in the notation of
  \cite{H}. By the results in loc. cit. $\tilde{T}$ is coherent (it
  is easy to see this directly).

Since $\Tscr$ is strongly graded we have that the category of
coherent graded $\Tscr$-modules is equivalent to the category of
coherent $\Oscr$-modules. This finishes the proof.
\end{proof}
\begin{lemma}
\label{ref:III.2.2-70} Assume that $B$ is a $\ZZ$-graded ring such that
$B_1
  B_b=B_{b+1}=B_b B_1$ for large $b$, and assume furthermore that there exists
  a non-zero $r>0$ such that  $B_{-br}
  B_{br}=B_{0}$ for all $b$. Then $B$ is strongly graded.
\end{lemma}
\begin{proof} For $b$ large we have
\[
B_{1} B_{-1}\supset B_1  B_{br-1} B_{-br}=B_{br} B_{-br}=B_0
\]
A similar argument shows that $B_{-1}B_1=B_0$.
\end{proof}
We can now prove the desired result.
\begin{proposition}
\label{ref:III.2.3-71} Let $S$ be as in Proposition
\ref{ref:III.2.1-69}, and assume
  in addition that $\qgr(S)$ is hereditary. Then
  $\qgr(S)$ is a finite direct sum of hereditary abelian categories which are
  either of the form $\mod(\Lambda)$ for a finite dimensional algebra $\Lambda$,
  or else of the form $\coh(\Oscr)$ where $\Oscr$ is a sheaf of
  classical hereditary orders
  over a non-singular  irreducible projective curve $Y$.
\end{proposition}
\begin{proof} This follows easily from the structure theory of
  hereditary noetherian rings.

By Proposition \ref{ref:III.2.1-69} it follows that  $\qgr(S)$ is equivalent
  to the category of coherent $\Oscr$-modules  where $\Oscr$ is a
  coherent sheaf of $\Oscr_X$-algebras and $X=\Proj C$.

  Let $\Zscr$ be the center of $\Oscr$, and let $\pi:Y\r X$ be the
  $X$-scheme such that $\pi_\ast \Oscr_Y=\Zscr$. Now for the purposes of
  the proof we may replace $Y$ by its connected components. So we
  assume that $Y$ is connected.  Let $U\subset Y$ be  affine open. Then by
  Proposition \ref{ref:2.4a} and \cite[p. 431]{Gabriel} it follows that  $\Oscr(U)$ is
  hereditary. So according to \cite[p 151]{MR} we have that $\Oscr(U)$ will be a
  direct sum of a finite dimensional hereditary algebra $\Lambda$ and
  a hereditary $\Oscr(U)$-algebra $H$ which is a direct sum of
  infinite dimensional prime
  rings. We identify $\Lambda$ with a sheaf of finite support on $Y$.
  If $\Lambda\neq 0$ then $\Lambda$ corresponds to a central
  idempotent in $\Oscr(U)$ which is zero outside the support of
  $\Lambda$. Hence we can extend it to a central idempotent in
  $\Oscr$.  Since we had assumed that the center of $\Oscr$ is
  connected, this yields $\Oscr=\Lambda$.

  Hence we only have to consider the case where $\Oscr(U)$ is a direct sum
  of prime hereditary rings  for all affine $U\subset Y$. In our case
the $\Oscr(U)$ are then direct sums
of orders in
  central simple algebras.

  Now it follows from \cite{RS1} that
 the center of
  a hereditary order is a  Dedekind ring.
 In particular $Y$ is a non-singular curve. Since $Y$ is also
 connected,
  it follows that $Y$ is irreducible. Thus if $U\subset Y$ is
  affine,
  then $\Oscr_Y(U)$ is a domain, whence $\Oscr(U)$ is actually prime.
  So $\Oscr(U)$ is a classical hereditary order. This finishes the argument.
\end{proof}
One of the conditions of the previous proposition was that $\qgr(S)$ should be hereditary. The
following lemma which is almost a tautology (using Proposition \ref{ref:2.3a}), tells us when
this condition holds.
\begin{lemma}
\label{ref:III.2.4-72} With the above hypotheses on $S$ we have that $\qgr(S)$ is hereditary
if and only if $S_P$ is hereditary
 for every non-maximal graded prime ideal $P$ in the center of $S$.
\end{lemma} \VDB
\section[Infinite Dynkin and tame quivers]{Hereditary abelian categories
associated to infinite Dynkin and tame  quivers} \label{ref:III.3-73} In this  section we
construct  some particular hereditary abelian noetherian categories which are obtained from
categories of type $\rep(Q)$ via tilting with respect to some torsion pair.

We start by computing the AR-quivers for the finite dimensional representations of the quivers
$A^\infty_\infty$, $D_\infty$ and $A_\infty$ (with respect to  a convenient orientation). The
first two are needed for the construction of noetherian hereditary abelian categories with no
nonzero projectives  and the third will be useful for giving an example that in the
non-noetherian case a connected hereditary abelian $\Ext$-finite category $\Cscr$ with Serre
duality and nonzero projectives is not necessarily determined by its associated quiver.

\subsection{The case of $A_\infty^\infty$}
\label{ref:III.3.1-74} 
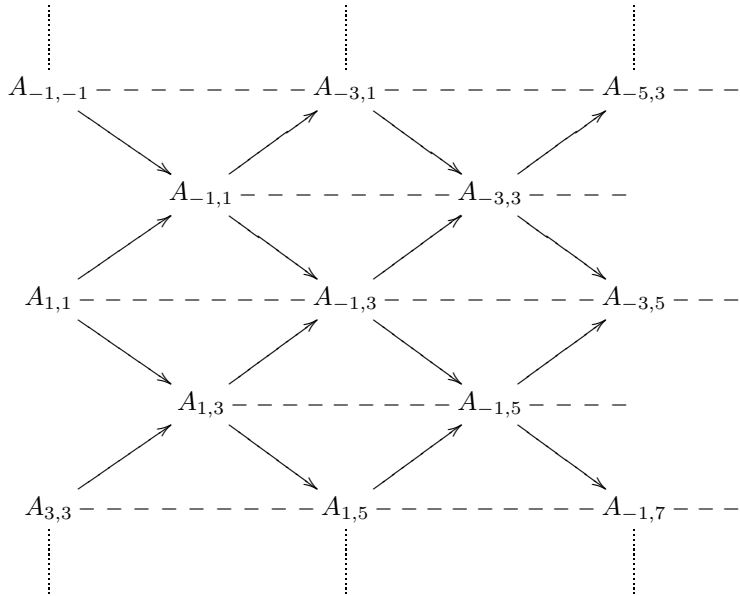
\begin{figure}
$$ \xymatrix{
\ar@{.}[d]&&\ar@{.}[d]&&\ar@{.}[d]&\\       %end first row
 A_{-1,-1}\ar[dr] \ar@{--}[rr] && A_{-3,1}\ar[dr]
\ar@{--}[rr] && A_{-5,3}\ar@{--}[r]&\\ %end second row
 & A_{-1,1} \ar[dr]\ar[ur]\ar@{--}[rr]
&& A_{-3,3}\ar[dr]\ar[ur]  \ar@{--}[r]&&\\ %end third row
 A_{1,1}\ar[ur]\ar[dr] \ar@{--}[rr] &&
A_{-1,3}\ar[ur]\ar[dr] \ar@{--}[rr]
&& A_{-3,5}\ar@{--}[r]&\\ %end fourth row
 & A_{1,3} \ar[dr]\ar[ur]\ar@{--}[rr]
&& A_{-1,5}\ar[dr]\ar[ur]  \ar@{--}[r]&&\\ %end fifth row
A_{3,3}\ar@{.}[d]\ar[ur] \ar@{--}[rr] && A_{1,5}\ar@{.}[d]\ar[ur]
\ar@{--}[rr]
&& A_{-1,7}\ar@{.}[d]\ar@{--}[r]&\\ %end sixth row
  &&
&& & %end seventh row
} $$ \caption{The preprojective component of $A^\infty_\infty$}
\label{ref:III.3.1-75}
\end{figure}

We consider the category of finite
dimensional representations of the quiver $$ \xymatrix{
\ar@{.}[dr]&& -2\ar[dl]\ar[dr] && 0\ar[dl]\ar[dr] && 2\ar[dl]\ar[dr] && 4\ar[dl]\ar@{.}[dr]&\\%
%end first row
&-3 && -1 && 1 && 3&& %end second row
} $$ of type $A^\infty_\infty$. Since each representation can be viewed as a representation of
some quiver of type $A_n$, we know the structure of the indecomposable representations. We
recall what they look like and at the same time we introduce a convenient notation:

$A_{n,m}$ ($n\le m$): $V_i=k$ for $n\le i\le m$ and $V_i=0$ otherwise. Then the AR-quiver
contains components given by Figures
\ref{ref:III.3.1-75},\ref{ref:III.3.2-77},\ref{ref:III.3.3-78},\ref{ref:III.3.4-79}
respectively.

\begin{figure}
$$ \xymatrix{
&\ar@{.}[d]&&\ar@{.}[d]&&\ar@{.}[d]\\       %end first row
 \ar@{--}[r] & A_{-6,2}\ar[dr] \ar@{--}[rr] && A_{-4,0}\ar[dr]
\ar@{--}[rr] && A_{-2,-2}\\ %end second row
& \ar@{--}[r]& A_{-4,2} \ar[dr]\ar[ur]\ar@{--}[rr]
&& A_{-2,0}\ar[dr]\ar[ur] \\ %end third row
 \ar@{--}[r] & A_{-4,4}\ar[ur]\ar[dr] \ar@{--}[rr] &&
A_{-2,2}\ar[ur]\ar[dr] \ar@{--}[rr]
&& A_{0,0}\\ %end fourth row
& \ar@{--}[r]& A_{-2,4} \ar[dr]\ar[ur]\ar@{--}[rr]
&& A_{0,2}\ar[dr]\ar[ur]  \\ %end fifth row
 \ar@{--}[r] & A_{2,6}\ar@{.}[d]\ar[ur] \ar@{--}[rr] &&
A_{0,4}\ar@{.}[d]\ar[ur] \ar@{--}[rr]
&& A_{2,2}\ar@{.}[d]\\ %end sixth row
 &  &&
&& %end seventh row
} $$ \caption{The preinjective component of $A^\infty_\infty$}
\label{ref:III.3.2-77}
\end{figure}
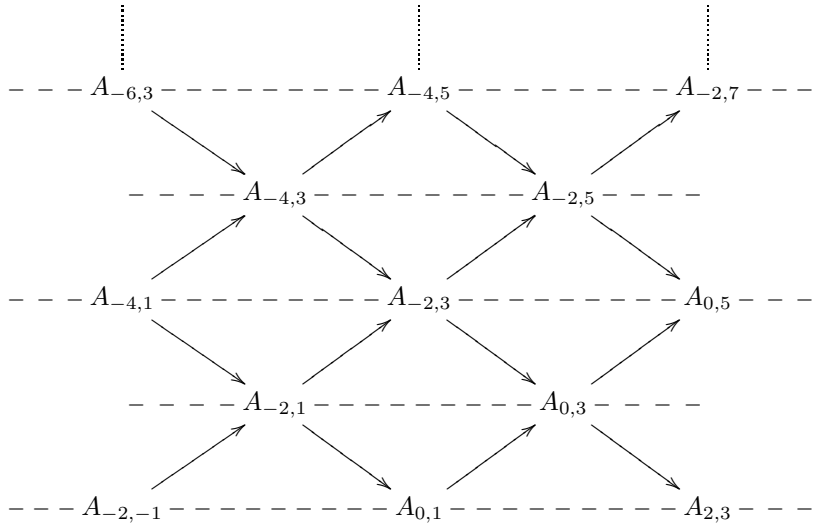
\begin{figure}
$$ \xymatrix{
&\ar@{.}[d]&&\ar@{.}[d]&&\ar@{.}[d]&\\       %end first row
 \ar@{--}[r] & \ar[dr]A_{-6,3} \ar@{--}[rr] && \ar[dr]
A_{-4,5}\ar@{--}[rr] && A_{-2,7}\ar@{--}[r]&\\ %end second row
& \ar@{--}[r]&  \ar[dr]\ar[ur]A_{-4,3}\ar@{--}[rr]
&& \ar[dr]\ar[ur]  A_{-2,5}\ar@{--}[r]&&\\ %end third row
 \ar@{--}[r] & \ar[ur]\ar[dr] A_{-4,1}\ar@{--}[rr] &&
\ar[ur]\ar[dr]A_{-2,3} \ar@{--}[rr]
&& A_{0,5}\ar@{--}[r]&\\ %end fourth row
& \ar@{--}[r]&  \ar[dr]\ar[ur]A_{-2,1}\ar@{--}[rr]
&& \ar[dr]\ar[ur] A_{0,3} \ar@{--}[r]&&\\ %end fifth row
 \ar@{--}[r] & A_{-2,-1}\ar[ur] \ar@{--}[rr] &&
\ar[ur]A_{0,1} \ar@{--}[rr]
&& A_{2,3}\ar@{--}[r]& %end sixth row
} $$ \caption{The first $\ZZ A_\infty$ component of
$A^\infty_\infty$} \label{ref:III.3.3-78}
\end{figure}
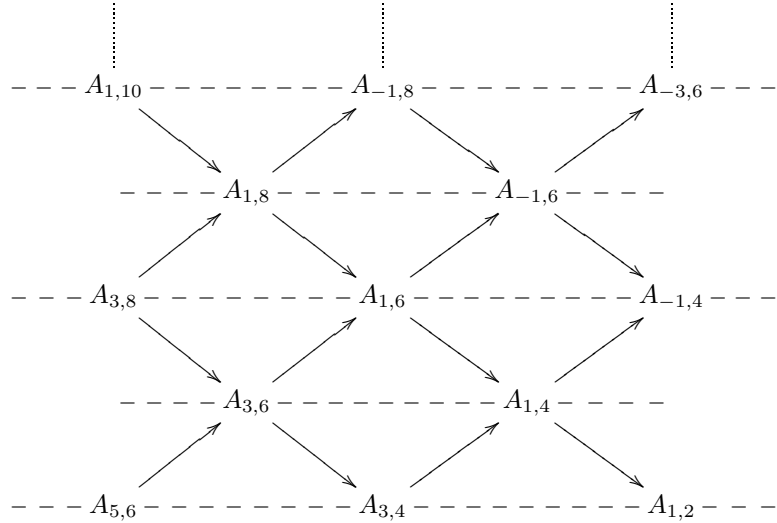
\begin{figure}
$$ \xymatrix{
&\ar@{.}[d]&&\ar@{.}[d]&&\ar@{.}[d]&\\       %end first row
 \ar@{--}[r] & \ar[dr]A_{1,10} \ar@{--}[rr] && \ar[dr]
A_{-1,8}\ar@{--}[rr] && A_{-3,6}\ar@{--}[r]&\\ %end second row
& \ar@{--}[r]&  \ar[dr]\ar[ur]A_{1,8}\ar@{--}[rr]
&& \ar[dr]\ar[ur]  A_{-1,6}\ar@{--}[r]&&\\ %end third row
 \ar@{--}[r] & \ar[ur]\ar[dr] A_{3,8}\ar@{--}[rr] &&
\ar[ur]\ar[dr]A_{1,6} \ar@{--}[rr]
&& A_{-1,4}\ar@{--}[r]&\\ %end fourth row
& \ar@{--}[r]&  \ar[dr]\ar[ur]A_{3,6}\ar@{--}[rr]
&& \ar[dr]\ar[ur] A_{1,4} \ar@{--}[r]&&\\ %end fifth row
 \ar@{--}[r] & A_{5,6}\ar[ur] \ar@{--}[rr] &&
\ar[ur]A_{3,4} \ar@{--}[rr]
&& A_{1,2}\ar@{--}[r]& %end sixth row
} $$ \caption{The second $\ZZ A_\infty$ component of
$A^\infty_\infty$} \label{ref:III.3.4-79}
\end{figure}
To see that the corresponding sequences are almost split we use
the known structure of almost split sequences for quivers of type
$A_n$.

We see that the preprojective component contains the $A_{n,m}$
with $m,n$ odd, the preinjective component contains those with
$m,n$ even, the $\ZZ A_\infty$ component in Figure \ref{ref:III.3.3-78}
contains the $A_{n,m}$ with $n$ even and $m$ odd and finally the
$\ZZ A_\infty$ component in Figure \ref{ref:III.3.4-79} contains the
$A_{n,m}$ with $n$ odd and
  $m$ even.

Hence all indecomposable representations occur in one of the four
above components, and so there are no other components.
\subsection{The case of $D_\infty$}
\label{ref:III.3.2-76} We consider the category of finite
dimensional representations for the quiver $$ \xymatrix{ &2\ar[dl]
\ar[d]\ar[dr] && 4\ar[dl]\ar[dr] && 6\ar[dl]\ar[dr] &&
8\ar[dl]\ar@{.}[dr] &\\ %end of first row
0 & 1 & 3 && 5 && 7 && } $$ of type $D_\infty$. Since each finite
dimensional representation can be viewed as a representation of
some quiver of type $D_n$, we know the structure of the
indecomposable representations. We recall what they look like, and
at the same time we introduce a notation convenient for describing
the distribution of modules in the AR-quiver (we always mean that
$V_j=0$ at the vertices which are not mentioned).
\begin{itemize}
\item[$A_{n,m}$:] $V_i=k$ for $2\le n\le i \le m$.
%(where $m\neq 1$ if
%  $n=0$).
\item[$A^{(1)}_{m}$:] $V_i=k$ for $i\le m$, $i\neq 1$, $V_1=0$.
\item[$A^{(0)}_m$:] $V_i=k$ for $i\le m$, $i\neq 0$, $V_0=0$.
\end{itemize}
%(Note that $A^{(0)}_m=A_{1,m}$ but it is convenient to use both
%notations.)
\begin{itemize}
\item[$B_{n,m}$:] $V_i=k^2$ for $2\le i\le n$, $V_i=k$ for
  $i=0,1$, and
  $n+1\le i\le m$ (for $1\le n\le  m$, $2\le m$).
\end{itemize}

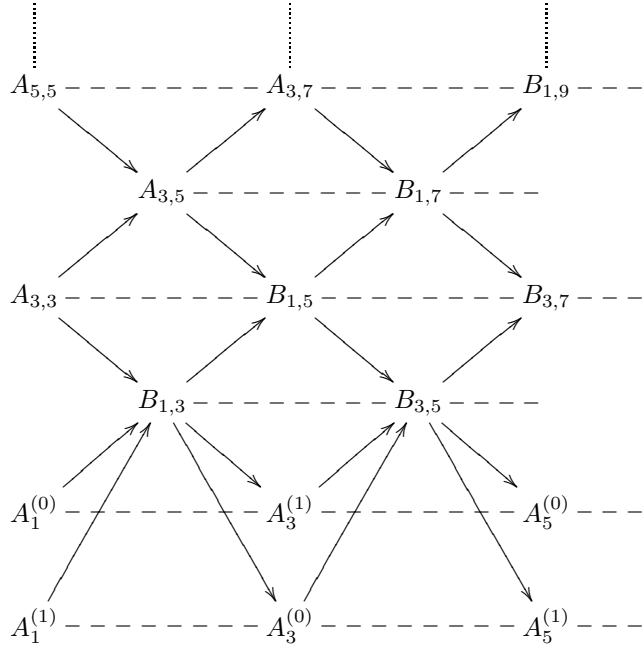
\begin{figure}
$$ \xymatrix{
\ar@{.}[d]&&\ar@{.}[d]&&\ar@{.}[d]&\\       %end first row
  A_{5,5}\ar[dr] \ar@{--}[rr] && A_{3,7}\ar[dr] \ar@{--}[rr]&
& B_{1,9}\ar@{--}[r]&\\ %end second row
& A_{3,5} \ar[dr]\ar[ur]\ar@{--}[rr]
&& B_{1,7}\ar[dr]\ar[ur]  \ar@{--}[r]&&\\ %end third row
 A_{3,3}\ar[ur]\ar[dr] \ar@{--}[rr] &&
B_{1,5}\ar[ur]\ar[dr] \ar@{--}[rr]
&&B_{3,7}\ar@{--}[r]&\\ %end fourth row
& B_{1,3} \ar[dr]\ar[ddr]\ar[ur]\ar@{--}[rr]
&& B_{3,5}\ar[dr]\ar[ddr]\ar[ur]  \ar@{--}[r]&&\\ %end fifth row
 A^{(0)}_1\ar[ur] \ar@{--}[rr] &&
A^{(1)}_3\ar[ur] \ar@{--}[rr]
&& A^{(0)}_5\ar@{--}[r]&\\ %end sixth row
 A_1^{(1)}\ar[uur] \ar@{--}[rr] &&
A_3^{(0)}\ar[uur] \ar@{--}[rr]
&& A^{(1)}_5\ar@{--}[r]&\\ %end seventh row
} $$ \caption{The preprojective component of $D_\infty$}
\label{ref:III.3.5-81}
\end{figure}

We now compute three components of the AR-quiver (see Figures
\ref{ref:III.3.5-81},\ref{ref:III.3.6-85},\ref{ref:III.3.7-86}). For those people that are
having trouble to spot the rule for Figure \ref{ref:III.3.7-86} it may help
to introduce the convention $B_{n,m}=B_{m,n}$.

It is not hard to see that the sequences corresponding to Figures
\ref{ref:III.3.5-81},\ref{ref:III.3.6-85},\ref{ref:III.3.7-86} are almost split by considering a
large enough subcategory of representations of a quiver of type
$D_n$.

\begin{figure}
$$ \xymatrix{
&\ar@{.}[d]&&\ar@{.}[d]&&\ar@{.}[d]\\       %end first row
 \ar@{--}[r] & A_{2,10}\ar[dr] \ar@{--}[rr] && A_{4,8}\ar[dr]
 \ar@{--}[rr] && A_{6,6}\\ %end second row
&\ar@{--}[r] & A_{2,8}\ar[ru]\ar[rd]\ar@{--}[rr]  && A_{4,6} \ar[ru]\ar[rd]%
\\ %end third row
 \ar@{--}[r] & B_{2,8}\ar[dr]\ar[ru] \ar@{--}[rr] && A_{2,6}\ar[dr]\ar[ru]
 \ar@{--}[rr] && A_{4,4}\\ %end fourth row
&\ar@{--}[r] & B_{2,6}\ar[ru]\ar[rd]\ar@{--}[rr]  && A_{2,4} \ar[ru]\ar[rd]%
\\ %end fifth row
\ar@{--}[r] & B_{4,6}\ar[dr]\ar[ddr]\ar[ru] \ar@{--}[rr] &&
B_{2,4}\ar[dr]\ar[ddr]\ar[ru]
 \ar@{--}[rr] && A_{2,2}\\ %end sixth row
&\ar@{--}[r] & A^{(1)}_4\ar[ru]\ar@{--}[rr]  && A_2^{(0)} \ar[ru]%
\\ %end seventh row
&\ar@{--}[r] & A^{(0)}_4\ar[ruu]\ar@{--}[rr]  && A_2^{(1)} \ar[ruu]%
 %end eight row
} $$ \caption{The preinjective component of $D_\infty$}
\label{ref:III.3.6-85}
\end{figure}
\begin{figure}
\tiny $$ \scriptscriptstyle \strut\hskip -1.5cm\xymatrix{
&\ar@{.}[d]&&\ar@{.}[d]&&\ar@{.}[d]&&\ar@{.}[d]&&\ar@{.}[d]&\\
%end first row
\ar@{--}[r]& B_{1,10}\ar[dr] \ar@{--}[rr] && \ar[dr]B_{3,8}
\ar@{--}[rr] && \ar[dr] B_{5,6}\ar@{--}[rr] &&
B_{4,7}\ar[dr]\ar@{--}[rr]&& B_{2,9}\ar@{--}[r]&\\
%end second row
& \ar@{--}[r] &B_{1,8} \ar[ur]\ar[dr] \ar@{--}[rr]&&
\ar[dr]\ar[ur]B_{3,6}\ar@{--}[rr] && \ar[dr]\ar[ur]
B_{4,5}\ar@{--}[rr]&& B_{2,7} \ar@{--}[r]\ar[ur]\ar[dr]&&\\
%end third row
\ar@{--}[r] & A_{3,8}
 \ar@{--}[rr] \ar[ur]\ar[dr] && \ar[ur]\ar[dr] B_{1,6}\ar@{--}[rr] &&
\ar[ur]\ar[dr]B_{3,4} \ar@{--}[rr]
&& B_{2,5}\ar@{--}[rr]\ar[ur]\ar[dr]&&A_{2,7} \ar@{--}[r]&\\ %end fourth row
&\ar@{--}[r]& A_{3,6} \ar@{--}[rr]\ar[dr]\ar[ur]&&
\ar[dr]\ar[ur]B_{1,4}\ar@{--}[rr] && \ar[dr]\ar[ur] B_{2,3}
\ar@{--}[rr]&&A_{2,5}\ar[ur]\ar[dr]\ar@{--}[r]&\\
%end fifth row
\ar@{--}[r]& \ar@{--}[rr] A_{5,6}\ar[ur]&& A_{3,4}\ar[ur]
\ar@{--}[rr] && \ar[ur]B_{1,2} \ar@{--}[rr]
&& A_{2,3}\ar@{--}[rr]\ar[ur]&& A_{4,5}\ar@{--}[r]& %end sixth row
} $$ \caption{The $\ZZ A_\infty$ component of $D_\infty$}
\label{ref:III.3.7-86}
\end{figure}
We see that the preprojective component contains the $A_{n,m}$, $A^{(0)}_m$ and $A^{(1)}_m$
with $m$ and $n$ odd or $0$ and the $B_{n,m}$ with $n,m$ odd. The preinjective component
contains the $A_{n,m}$, $A^{(0)}_m$, $A^{(1)}_m$ with $m$ even and $n$ even $\neq 0$ and the
$B_{n,m}$ with $n,m$ even. The third component contains the $A_{n,m}$ and $B_{n,m}$ with $n+m$
odd and the $A_{0,m}$ with $m$ even.

Hence all indecomposable representations lie in one of the above three components, and hence
there are no other components.
\subsection{The case of $A_\infty$}
\label{ref:III.3.3-80} We consider the category of finite
dimensional representations of the quiver $$ \xymatrix{
& 2\ar[dl]\ar[dr] && 4\ar[dl]\ar[dr] && 6\ar[dl]\ar[dr] && 8\ar[dl]\ar@{.}[dr]&\\%
%end first row
1 && 3 && 5 && 7&& %end second row
} $$ of type $A_\infty$. We use the same notation as in the case $A_\infty^\infty$ for the
indecomposable representations. Then the AR-quiver  contains components given by Figures
\ref{ref:III.3.8-87},\ref{ref:III.3.9-88} respectively.
\begin{figure}
$$ \xymatrix{
%\ar@{.}[d]&&\ar@{.}[d]&&\ar@{.}[d]&\\       %end first row
 A_{1,1}\ar[dr] \ar@{--}[rr] && A_{2,3}\ar[dr]
\ar@{--}[rr] && A_{4,5}\ar@{--}[r]&\\ %end first row
 & A_{1,3} \ar[dr]\ar[ur]\ar@{--}[rr]
&& A_{2,5}\ar[dr]\ar[ur]  \ar@{--}[r]&&\\ %end second row
 A_{3,3}\ar[ur]\ar[dr] \ar@{--}[rr] &&
A_{1,5}\ar[ur]\ar[dr] \ar@{--}[rr]
&& A_{2,7}\ar@{--}[r]&\\ %end third row
 & A_{3,5} \ar[dr]\ar[ur]\ar@{--}[rr]
&& A_{1,7}\ar[dr]\ar[ur]  \ar@{--}[r]&&\\ %end fourth row
A_{5,5}\ar@{.}[d]\ar[ur] \ar@{--}[rr] && A_{3,7}\ar@{.}[d]\ar[ur]
\ar@{--}[rr]
&& A_{1,9}\ar@{.}[d]\ar@{--}[r]&\\ %end fifth row
  &&
&& & %end seventh row
} $$ \caption{The preprojective component of $A_\infty$}
\label{ref:III.3.8-87}
\end{figure}
\begin{figure}
$$ \xymatrix{
%&\ar@{.}[d]&&\ar@{.}[d]&&\ar@{.}[d]\\       %end first row
% \ar@{--}[r] & A_{3,4}\ar[dr] \ar@{--}[rr] && A_{1,2}\ar[dr]
%\ar@{--}[rr] && A_{-2,-2}\\ %end first row
& \ar@{--}[r]& A_{3,4} \ar[dr]\ar@{--}[rr]
&& A_{1,2}\ar[dr] \\ %end first row
 \ar@{--}[r] & A_{3,6}\ar[ur]\ar[dr] \ar@{--}[rr] &&
A_{1,4}\ar[ur]\ar[dr] \ar@{--}[rr]
&& A_{2,2}\\ %end second row
& \ar@{--}[r]& A_{1,6} \ar[dr]\ar[ur]\ar@{--}[rr]
&& A_{2,4}\ar[dr]\ar[ur]  \\ %end third row
 \ar@{--}[r] & A_{1,8}\ar[ur]\ar[dr] \ar@{--}[rr] &&
A_{2,6}\ar[ur] \ar[dr]\ar@{--}[rr]
&& A_{4,4}\\ %end fourth row
& \ar@{--}[r]& A_{2,8}\ar@{.}[d]\ar[ur]\ar@{--}[rr]
&& A_{4,6}\ar@{.}[d]\ar[ur]&  \\ %end fifth row
 &  &&
&& %end six row
} $$ \caption{The preinjective component of $A_\infty$}
\label{ref:III.3.9-88}
\end{figure}
The $A_{n,m}$ with $m$ odd are in the preprojective component, and
those with $m$ even are in the preinjective component. Hence there
are no other components.

\subsection{The $\ZZ A_\infty^\infty$ category}
\label{ref:III.3.4-82} Let $\Cscr$ be $\rep(Q)$
 where $Q$ is the quiver $A_\infty^\infty$ with zigzag orientation (as in
\S\ref{ref:III.3.1-74}), and let $\Cscr_1$ be obtained from $\Cscr$ by tilting with respect to
the split torsion pair $(\Tscr,\Fscr)$, where  $\Tscr$ consists of the preinjectives (see
Figure \ref{ref:II.3.2-48}). Then $\Cscr_1$ is an $\Ext$-finite hereditary abelian category
with no nonzero projectives or injectives, derived equivalent to $\Cscr$. So $\Cscr_1$
satisfies Serre duality. The AR-quiver of $\Cscr_1$ consists of one $\ZZ A^\infty_\infty$
component $\Sigma$ and two $\ZZ A_\infty$ components.

 Below we show that $\Cscr_1$ is noetherian and we show in addition
 that $\Cscr_1$ has two $\tau$-orbits of simples. But first we give some preliminary results.

 In an $AR$-quiver the subquiver of the form
$$
\xymatrix{
    & B_1 \ar@{->}[rd]^{\alpha_1} \ar@{.}[dd]  & & & \\
    \tau C \ar[ur]^{\beta_1} \ar[dr]_{\beta_n} & & C & \text{given by an almost split sequence} &  \\
    & B_n \ar[ur]_{\alpha_n} & & &
    }
$$

$0\to\tau C\to B_1\oplus\cdots\oplus B_n\to C\to 0$,
with the $B_i$ indecomposable, (and which we here assume to be nonisomorphic) is called a {\it mesh}. The mesh category of a component
$\Cscr'$ is the path category of $\Cscr'$ modulo the relations
$\sum\limits^n_{i=1}\alpha_i\beta_i$ given by the meshes. The component $\Cscr '$ of an AR-quiver is said to be \emph{standard} if the full
 subcategory of indecomposable objects corresponding to the vertices of $\Cscr '$ is
 equivalent to the mesh category of $\Cscr'$.
\begin{lemmas}
\label{ref:III.3.4.1-83} Let  $\Cscr$ and $\Cscr_1$ be as above.
\begin{enumerate}
\item Let $X$ and $Y$ be indecomposable objects belonging to a component of the AR-quiver of
$\Cscr$ or $\Cscr_1$. If $f:X\r Y$ is a nonzero non-isomorphism, then $f$ is a composition of
irreducible maps. In particular there is a path between the corresponding vertices.
\item The components of the AR-quiver of $\Cscr$ and of $\Cscr_1$ are standard.
\end{enumerate}
\end{lemmas}
\begin{proof}
Part 1. for $\Cscr$ follows directly from the explicit description of the indecomposable
objects in $\Cscr$. The claim for $\Cscr_1$ can be reduced to the result for $\Cscr$, by
possibly applying high enough powers of $\tau^{-1}$.

To show that any component $\Cscr '$ of $\Cscr$ or $\Cscr_1$ is standard, one defines as usual
a functor from the path category of $\Cscr '$ to the subcategory of indecomposable objects of
$\Cscr$ given by $\Cscr '$, in such a way that there is induced a functor from the mesh
category of $\Cscr '$. Using 1., this will be an equivalence.
\end{proof}
The next result, which gives information on maps between components, is easily verified.
\begin{lemmas}
\label{lemIII.3.4.2} Let $\Cscr$ and $\Cscr_1$ be as above.
\begin{enumerate}
\item For $\Cscr$ there are no nonzero maps from the preinjective component to any other
component and to the preprojective component from any other component, and also no nonzero
maps between the two $\ZZ A_\infty$-components.
\item For $\Cscr_1$ there are no nonzero maps from the $\ZZ A_\infty$-components to the $\ZZ
A_\infty^\infty$-component, and no nonzero maps between the two $\ZZ A_\infty$-components.
\end{enumerate}
\end{lemmas}
We state explicitly the following special cases which will be used later.
\begin{lemmas}
\label{lemIII.3.4.3} Let $\Cscr_1$ be as above.
\begin{enumerate}
\item If $U$ and $W$ are indecomposable with $\Hom (W,U)\neq 0$ and $U$ is in the $\ZZ
A_\infty^\infty$-component of $\Cscr_1$, then there is a path from $W$ to $U$ in the
AR-quiver.
\item If $X$ is indecomposable in $\Cscr_1$, then $\End_{\Cscr_1} (X)\simeq k$.
\end{enumerate}
\end{lemmas}
\begin{propositions}
\label{ref:III.3.4.2-84}
 $\Cscr_1$ has exactly two ${{\tau}}$-orbits of
  simples.
\end{propositions}
\begin{proof}
We claim that the only simple objects in $\Cscr_1$ are the objects
of length two in the $\ZZ A_\infty$ components of $\Cscr$. This
clearly implies what we have to show.

First we note that the other indecomposables cannot be simple.
Indeed let $X$ be one of those other indecomposables. Assume first
that $X\in\Fscr$. Thus $X$ can be considered as an object of
$\Cscr$. Assume now in addition that $X$ is not a length one
projective. Then it is easy to see directly that we can make a
non-trivial short exact sequence
\begin{equation}
\label{ref:III.3.1-89} 0\r U\r X\r W\r 0
\end{equation}
with $U,W\in\Fscr$. Hence this sequence remains exact in
$\Cscr_1$, and thus $X$ is not simple.

If $X$ happens to lie in $\Tscr[-1]$, or is given by a length one projective in $\Cscr$, then
we first apply a power of ${{\tau}}^{-1}$, and then we apply the above reasoning.

Now we prove the converse. So let $X$ be an indecomposable object
of length  two in one of the $\ZZ A_\infty$ components of $\Cscr$.
We have to prove that $X$ becomes simple in $\Cscr_1$.

Assume now that $X$ is not simple in $\Cscr_1$. Thus there is an exact sequence in $\Cscr_1$
as in \eqref{ref:III.3.1-89}, but this time the indecomposable summands of $W$ must be in the
same component as $X$ (by lemma \ref{lemIII.3.4.2}). By applying a high power of
${{\tau}}^{-1}$ we may assume that the summands of $U$ are in $\Fscr$, and hence
\eqref{ref:III.3.1-89} represents an exact sequence in $\Cscr$. It is easy to see that this is
impossible.
\end{proof}
\begin{lemmas}
\label{ref:III.3.4.3-90} As in Figure \ref{ref:II.3.2-48}
let $\Sigma$ be the $\ZZ
  A_\infty^\infty$ component of
  $\Cscr_1$. Then the objects in $\Sigma$ generate $\Cscr_1$.
\end{lemmas}
\begin{proof} This follows directly from the fact that the projectives
  generate $\Cscr$.
\end{proof}

\begin{lemmas}
\label{ref:III.3.4.6-93} One has for $X,Y\in\Sigma$
\begin{equation}
\dim\Hom_{\Cscr_1}(X,Y)=
\begin{cases}
1&\text{if there is a path from $X$ to $Y$}\\ 0&\text{otherwise}
\end{cases}
\end{equation}
\end{lemmas}
\begin{proof}
  This can be shown in many ways. The most direct method is to compute
  in $\Cscr$, but it is most elegant to consider the
  function $\phi_X=\dim \Hom(X,-)$. This function is $1$ on $X$  and zero on the objects which are not on a path
  starting in $X$ (by lemma \ref{lemIII.3.4.3}). Furthermore $\phi_X$ is
  additive on almost split sequences not ending in $X$. This is enough
  to determine $\phi_X$ completely.
\end{proof}
We conclude:
\begin{lemmas}
\label{ref:III.3.4.7-94} If $X\in\Sigma$, then $X$ does not
contain a non-trivial direct sum.
\end{lemmas}
\begin{proof}
Assume $Y\oplus Z\subset X$ with $Y,Z$ indecomposable. Then we can
find an indecomposable $C$ which has a path to both $Y$ and $Z$.
Thus we obtain.
\[
1=\dim\Hom(C,X)\ge \dim\Hom(C,Y)+\dim \Hom(C,Z)=2
\]
which is a contradiction.
\end{proof}
\begin{propositions}
\label{ref:III.3.4.8-95} The category $\Cscr_1$ is
noetherian.
\end{propositions}
\begin{proof}
By lemma \ref{ref:III.3.4.3-90} it suffices to prove that
the objects in $\Sigma$ are noetherian. So let $C$ be an  object
in $\Sigma$, and assume that $C$ is not noetherian. Then there is
some infinite chain
\[
C_0\subsetneq C_1\subsetneq C_2\subsetneq \cdots \subset C
\]
By lemma \ref{ref:III.3.4.7-94} each of the $C_i$ is indecomposable, and by lemma
\ref{lemIII.3.4.2} each of the $C_i$ is in $\Sigma$. Thus by  lemma \ref{lemIII.3.4.3} each of
the $C_i$ is on a path from $C_0$ to $C$. Since there are only a finite number of objects on
such paths, this is a contradiction.
\end{proof}
\subsection{The $\ZZ D_\infty$ category}
\label{ref:III.3.5-96} This section closely parallels the previous one, but some of the
arguments are slightly more complicated. Let $\Cscr$ be $\rep(Q)$, where $Q$ is the quiver
$D_\infty$ with zigzag orientation (as in \ref{ref:III.3.2-76}), and let $\Cscr_1$ be obtained
from $\Cscr$ by tilting with respect to the split torsion pair $(\Tscr,\Fscr)$, where  $\Tscr$
consists of the preinjectives (see Figure \ref{ref:II.3.2-48}). Then $\Cscr_1$ is an
$\Ext$-finite hereditary abelian category with no nonzero projectives or injectives, derived
equivalent to $\Cscr$. So $\Cscr_1$ satisfies Serre duality. The AR-quiver of $\Cscr_1$
consists of one $\ZZ D_\infty$ component $\Sigma$ and one $\ZZ A_\infty$ component.

We want to show that $\Cscr_1$ is noetherian, and that $\Cscr_1$ has exactly one $\tau$-orbit
of simples. Some of the preliminary results we need for this are similar to the corresponding
results in \ref{ref:III.3.4-82}, and hence most of them are omitted.

We start with stating the following results, where the first one is of interest in itself and the
second one is different from the corresponding result for the $\ZZ A_\infty^\infty$ category.
\begin{lemmas}
\label{lemIII.3.5.1} The components of $\Cscr$ and $\Cscr_1$ are standard.
\end{lemmas}
\begin{propositions}
\label{lemIII.3.5.2} $\Cscr_1$ has exactly one $\tau$-orbit of simple objects.
\end{propositions}

\begin{figure}
\tiny $$ \scriptscriptstyle \strut\hskip -1.5cm\xymatrix{
&\ar@{.}[d]&&\ar@{.}[d]&&\ar@{.}[d]&&\ar@{.}[d]\\ %end first row
\ar@{--}[r]&1\ar[dr]\ar@{--}[rr]&& 1\ar@{--}[rr]\ar[dr]
&& 0\ar@{--}[rr]\ar[dr] && 0 \ar@{--}[r]&\\ %end second row
&\ar@{--}[r]&1\ar@{--}[rr]\ar[ur]\ar[dr]&&
1\ar@{--}[rr]\ar[ur]\ar[dr]&&0\ar@{--}[r]\ar[ur]\ar[dr]&&\\
%end third row
\ar@{--}[r]&1\ar[dr]\ar[ur]\ar@{--}[rr]&&
1\ar@{--}[rr]\ar[dr]\ar[ur]
&& 1\ar@{--}[rr]\ar[dr] \ar[ur]&& 0\ar@{--}[r]&\\ %end fourth row
&\ar@{--}[r]&1\ar@{--}[rr]\ar[ur]\ar[dr]\ar[ddr]&&
1\ar@{--}[rr]\ar[ur]\ar[dr]\ar[ddr]&&1\ar@{--}[r]\ar[ur]\ar[dr]\ar[ddr]&&\\
%end fifth row
\ar@{--}[r]&0\ar[ur]\ar@{--}[rr]&& 1\ar@{--}[rr]\ar[ur]
&& 0\ar@{--}[rr] \ar[ur]&& *+[o][F-]{1} \ar@{--}[r]&\\ %end sixth row
\ar@{--}[r]&1\ar[uur]\ar@{--}[rr]&& 0\ar@{--}[rr]\ar[uur]
&& 1\ar@{--}[rr] \ar[uur]&& 0 \ar@{--}[r]& %end seventh row
} $$ \caption{$\dim \Hom(X,Y)$ for $X,Y\in\Sigma$ with $X$ variable} \label{ref:III.3.11-102}
\end{figure}
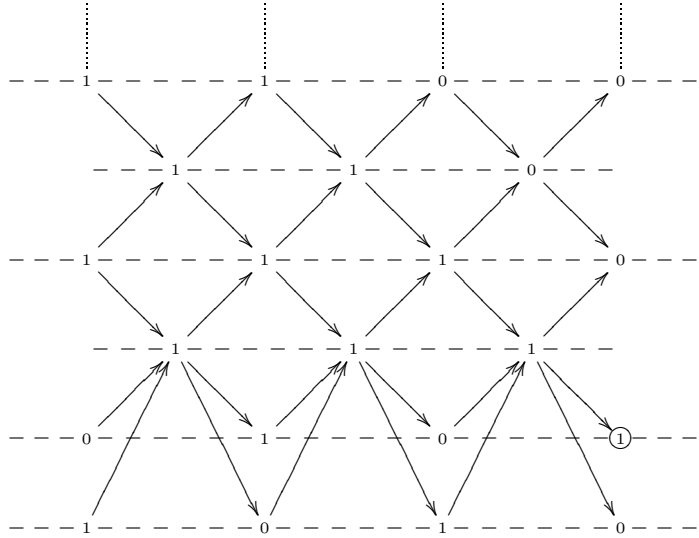

Now we need an analog of lemma \ref{ref:III.3.4.6-93}, but as is to be expected the result is
slightly more complicated in the present case. Therefore the result is most easily represented
in some diagrams. Figures \ref{ref:III.3.10-101},\ref{ref:III.3.11-102} give the values of
$\dim \Hom(X,Y)$ for $X,Y\in\Sigma$ where $X$ is varying and $Y$ is fixed. Figure
\ref{ref:III.3.11-102} corresponds to the case where $Y$ is on one of the ending vertices of
$D_\infty$, and Figure \ref{ref:III.3.10-101} corresponds to the general case. These diagrams
are obtained in the same way as lemma \ref{ref:III.3.4.6-93}. We have encircled the vertex
corresponding to $Y$ in each of the cases.
\begin{figure}
\tiny $$ \scriptscriptstyle \strut\hskip -1.5cm\xymatrix{
&\ar@{.}[d]&&\ar@{.}[d]&&\ar@{.}[d]&&\ar@{.}[d]\\ %end first row
\ar@{--}[r]&2\ar[dr]\ar@{--}[rr]&& 1\ar@{--}[rr]\ar[dr]
&& 1\ar@{--}[rr]\ar[dr] && 0 \ar@{--}[r]&\\ %end second row
&\ar@{--}[r]&2\ar@{--}[rr]\ar[ur]\ar[dr]&&
1\ar@{--}[rr]\ar[ur]\ar[dr]&&1\ar@{--}[r]\ar[ur]\ar[dr]&&\\
%end third row
\ar@{--}[r]&2\ar[dr]\ar[ur]\ar@{--}[rr]&&
2\ar@{--}[rr]\ar[dr]\ar[ur]
&& 1\ar@{--}[rr]\ar[dr] \ar[ur]&& *+[o][F-]{1} \ar@{--}[r]&\\ %end fourth row
&\ar@{--}[r]&2\ar@{--}[rr]\ar[ur]\ar[dr]\ar[ddr]&&
2\ar@{--}[rr]\ar[ur]\ar[dr]\ar[ddr]&&1\ar@{--}[r]\ar[ur]\ar[dr]\ar[ddr]&&\\
%end fifth row
\ar@{--}[r]&1\ar[ur]\ar@{--}[rr]&& 1\ar@{--}[rr]\ar[ur]
&& 1\ar@{--}[rr] \ar[ur]&& 0 \ar@{--}[r]&\\ %end sixth row
\ar@{--}[r]&1\ar[uur]\ar@{--}[rr]&& 1\ar@{--}[rr]\ar[uur]
&& 1\ar@{--}[rr] \ar[uur]&& 0 \ar@{--}[r]& %end seventh row
} $$ \caption{$\dim \Hom(X,Y)$ for $X,Y\in\Sigma$ with $X$ variable} \label{ref:III.3.10-101}
\end{figure}
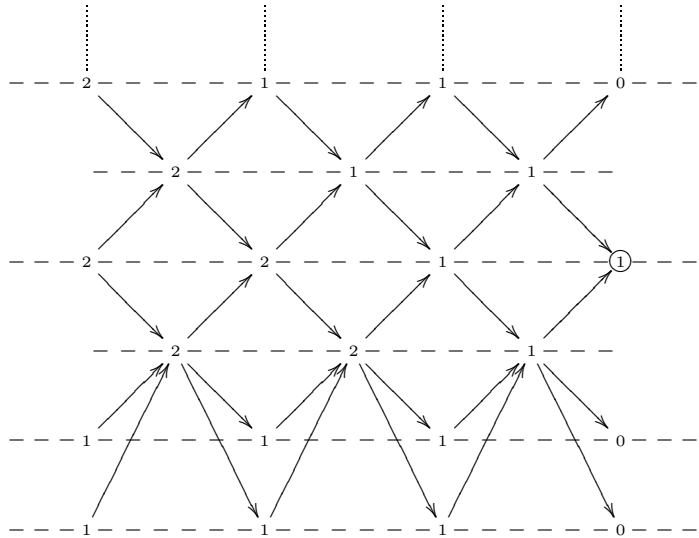

{}{}From these diagrams we compute.
\begin{lemmas}
\label{ref:III.3.5.5-103}
 If $X,Y\in\Sigma$ then $\dim \Hom(X,Y)\le 2$.
\end{lemmas}
This lemma allows us to prove the following:
\begin{lemmas}
\label{ref:III.3.5.6-104} If $X\in\Sigma$, then $X$ does not
contain a non-trivial direct sum of three objects.
\end{lemmas}
\begin{proof}
Assume $W\oplus Y\oplus Z\subset X$ with $W,Y,Z$ indecomposable.
Then according to the diagrams \ref{ref:III.3.10-101},\ref{ref:III.3.11-102} we
can find an indecomposable $C$ such that $\Hom(C,-)$ is non-zero
when evaluated on   $W,Y$ and $Z$. Thus we obtain
\[
2\ge\dim\Hom(C,X)\ge \dim\Hom(C,W)+\dim\Hom(C,Y)+\dim \Hom(C,Z)\ge
3
\]
which is a contradiction.
\end{proof}
\begin{propositions} The category $\Cscr_1$ is noetherian.
\end{propositions}
\begin{proof}
  As in \ref{ref:III.3.4-82} it suffices to prove that
  the objects in $\Sigma$ are noetherian. So let $C$ be an object in
  $\Sigma$, and assume that $C$ is not noetherian. Then there is some
  infinite chain
\[
C_0\subsetneq C_1\subsetneq C_2\subsetneq \cdots \subset C
\]
After possibly replacing $(C_i)_i$ by a subsequence we may by
lemma
 \ref{ref:III.3.5.6-104} assume that either each of the $(C_i)_i$ is
 indecomposable, or else that they all have exactly two indecomposable
 summands. The first case is dealt with exactly as in Proposition
 \ref{ref:III.3.4.8-95}, so we will concentrate on the second case.

Thus we now assume that $C_i=D_i\oplus E_i$ with $D_i$ and $E_i$ indecomposable. Using lemma
\ref{ref:III.3.5.6-104} it is easily seen that possibly after interchanging $D_{i+1}$ and
$E_{i+1}$ we may assume that there exist non-zero maps $D_i\r D_{i+1}$ and $E_i\r E_{i+1}$.
Thus as in the proof of Proposition \ref{ref:III.3.4.8-95} we must have $D_i\cong D_{i+1}$ and
$E_i\cong E_{i+1}$ for large~$i$, and hence also $C_i\cong C_{i+1}$. {}From the fact that
$\Cscr_1$ is $\Hom$-finite, it then easily follows that the inclusion $C_i\hookrightarrow
C_{i+1}$ must actually be an isomorphism. This finishes the proof.
\end{proof}
\subsection{An example}
\label{ref:III.3.6-105} We now give an example that shows that the important Theorem
\ref{ref:II.4.2-55} fails for non-noetherian hereditary abelian $\Ext$-finite categories with
Serre duality, even when they have noetherian injectives.

Let $\Cscr_1$ be the $\ZZ A_\infty^\infty$-category as defined in \S\ref{ref:III.3.4-82}. For
one of the $ZA_\infty$-components consider a section $Q_1\subset \ZZ A_\infty$ with zig-zag
orientation. Let $\Tscr_1$ consist of direct sums of indecomposable objects which are in the
chosen component to the right of $Q_1$, and let $\Fscr_1$ consist of direct sums of the other
indecomposables. Then it follows from  \ref{ref:III.3.4-82} together with Serre duality that
$(\Tscr_1,\Fscr_1)$ is a split torsion pair in $\Cscr_1$.

 Let
$\Dscr$ be the tilted category. Then $\Dscr$ is a hereditary abelian category whose
projectives are given by $Q_1$ and which is derived equivalent to $\Cscr_1$. It follows from
lemma \ref{ref:II.3.2-46} that $\Dscr$ has noetherian injectives. Hence if Theorem
\ref{ref:II.4.2-55} were to hold for $\Dscr$, then by Corollary \ref{ref:II.2.8-41} it would
follow that $\Dscr=\wrep(Q_1)=\rep(Q_1)$ (using the fact that $\Dscr$, being derived
equivalent to $\rep(Q)$, is indecomposable by construction). However this contradicts the
determination of the AR-quiver of $\rep Q_1$ which was made in \S\ref{ref:III.3.3-80}.

\subsection{Other interpretations of the $\ZZ A^\infty_\infty$
  and the $\ZZ D_\infty$-category}
\label{ref:III.3.7-106} In this section we define some $\Ext$-finite noetherian hereditary
abelian categories with almost split sequences and no nonzero projectives and injectives.
Ultimately when we have our general classification result (Theorem \ref{ref:IV.5.2-152}) it
will turn out that these are actually the same as the $\ZZ A_\infty^\infty$ and $\ZZ D_\infty$
category. The constructions in this section are very easy, but unfortunately for one of them
we have to assume that the characteristic is not two.

\noindent {\bf First construction.} Following \cite{SmithZhang} we consider $R=k[u,v]$ as
$\ZZ^2$-graded ring with $\deg u=(1,0)$ and $\deg v=(0,1)$ and let $\Cscr$ be the category of finitely
generated graded $R$-modules modulo the finite dimensional ones.  Then it is shown in
\cite{SmithZhang} that $\Cscr$ is a hereditary abelian category satisfying Serre duality (this
last fact follows easily from graded local duality) and it is also $\Ext$-finite. It is also
shown in \cite{SmithZhang} that there are precisely two $\tau$-orbits of simples. Hence
invoking Proposition \ref{ref:III.3.4.2-84} above and Theorem \ref{ref:IV.5.2-152} below it
follows that $\Cscr$ is the $\ZZ A_\infty^\infty$ category.

\noindent {\bf Second construction.} Now we assume that the characteristic of $k$ is not two.
Let $\Cscr$ be as in the previous construction. The category $\Cscr$ has a natural
automorphism $\sigma$ of order two which sends a graded module $(M_{ij})_{ij}$ to
$(M_{ji})_{ij}$ and which exchanges the $u$ and $v$ action. Let $\Dscr$ be the category of
$\ZZ/2\ZZ$ equivariant objects in $\Cscr$, i.e. pairs $(M,\phi)$ where $\phi$ is an
isomorphism $M\r \sigma(M)$ satisfying $\sigma(\phi)\phi= \Id_M$. It now follows easily from
the general results in \cite{ReitenRiedtmann} that $\Dscr$ is an $\Ext$-finite hereditary
abelian category with almost split sequences and no projectives or injectives, whence $\Cscr$
satisfies Serre duality by Corollary \ref{ref:I.3.2-20}. The two orbits of simples in $\Cscr$
collapse to a single orbit in $\Dscr$.  Hence invoking Proposition \ref{lemIII.3.5.2} above
and Theorem \ref{ref:IV.5.2-152} it follows $\Cscr$ is the $\ZZ D_\infty$ category.

\chapter[Hereditary noetherian abelian categories with no
projectives or \ldots]{Hereditary noetherian abelian categories with no
projectives or injectives.}\label{ref:IV-107}
%\subsection{Statement of the main result}
Let $\Cscr$ be a connected $k$-linear
  $\Ext$-finite noetherian hereditary abelian category
    having almost split
  sequences and no nonzero projectives or injectives, or equivalently, having
  a Serre functor and no nonzero projectives.
In this chapter we show that if  $\Cscr$ has some object which is
not of
  finite length, then
$\Cscr$ has one of the following
  forms.
\begin{enumerate}
\item If $\Cscr$ has an infinite number of $\tau$-orbits of simples
  then $\Cscr$ is the category of coherent sheaves over a hereditary
  order over a connected non-singular proper curve (as in
  \ref{ref:III.2-68}).
\item If $\Cscr$ has exactly two $\tau$-orbits of simples, then $\Cscr$
  is the $\ZZ A^\infty_\infty$-category (as in \ref{ref:III.3.4-82}).
\item If $\Cscr$ has exactly one $\tau$-orbit of simples then $\Cscr$
  is the $\ZZ D_\infty$-category (as in \ref{ref:III.3.5-96}).
\end{enumerate}

We first show that when $\Cscr$ is connected and  satisfies the above hypothesis, then either
each $\tau$-orbit of simple objects is finite or each such orbit is infinite. In the first
case we show that there is a unique simple object in the quotient category $\Cscr / \Tscr$,
where $\Tscr$ consists of the objects of finite length, coming from an object $F$ in $\Cscr$
with no summands in $\Tscr$. We are  able to construct an automorphism $t \colon \Cscr \r
\Cscr$ such that $(F,t)$ is ample, and then we use results of \cite{AZ} to show that $\Cscr$
is equivalent to some quotient category $\qgr S$. Then $\Cscr$ is equivalent to the category
of coherent sheaves over a hereditary order over a connected non-singular proper curve. In the
second case we show that $\Cscr$ is the $\mathbb{Z}A^\infty_\infty$-category if there is one simple
object in $\Cscr / \Tscr $, and the $\mathbb{Z}D_\infty$-category if there are two simple objects.

In this chapter $\Cscr$ will satisfy the above assumptions, except that we will  not a priori
assume that $\Cscr$ is connected. As usual we denote the Serre functor by $F$ and we write
$\tilde{\tau}=F[-1]$. The induced autoequivalence of $\Cscr$ is as usual denoted by $\tau$.

\section{Preliminaries}
\label{ref:IV.1-108} In this section we show that if $\Cscr$ is connected, then the quotient
category $\Cscr / \Tscr $ (where $\Tscr$ consists of the objects of finite length) contains at
most two simple objects, and if there are two, then they are permuted by $\tau$ modulo
$\Tscr$. We also show that if $A$ is not in $\Tscr$, there is some nonzero map from $A$ to a
simple object in every $\tau$-orbit.

 We call the objects in $\Tscr$ ``torsion''.
Analogously we call an object in $\Cscr$ torsion free if it contains no subobject in $\Tscr$.
Let $\Fscr$ be the full subcategory of $\Cscr$ consisting of torsion free objects. Clearly
$\Tscr$ and $\Fscr$ are stable under $\tau$ (as they are under any autoequivalence). In the
sequel we denote the quotient functor $\Cscr\r \Cscr/\Tscr$ by $\pi$.

Let $(S_i)_{i\in I}$ be the simple objects in $\Cscr$, and let $v$ be the permutation of $I$
defined by $\tau S_i\cong S_{v(i)}$. We say that $i\sim j$ if $i$ and $j$ are in the same
$v$-orbit. If $J\subset I$ then we denote by $\Tscr_J$ the full subcategory of $\Cscr$
consisting of objects whose Jordan-Holder quotients are in $(S_j)_{j\in J}$. We have
\begin{equation}
\label{ref:IV.1.1-109} \Tscr=\oplus_{o\in I/\!\sim} \Tscr_o
\end{equation}
where the $\Tscr_o$ are described by the results in
\ref{ref:III.1-64}.

We will now prove some elementary properties of  torsion free
objects. Note first that by Serre duality $\Ext^1(\Fscr,\Tscr)=0$.
In other words the torsion free objects behave as projectives with
respect to the torsion objects. {}From this we easily deduce:
\begin{lemma}
\label{ref:IV.1.1-110}
\begin{enumerate}
\item Every object in $\Cscr$ is the direct sum of a torsion object
  and a torsion free object.
\item The  categories $\Tscr_o$ for $o\in I/\sim$ are stable under
  essential extensions in $\Cscr$.
\end{enumerate}
\end{lemma}
If $J\subset I$ is $v$-stable and $A\in \Cscr$, then we define the $\Tscr_J$-topology on $A$
as the linear topology generated by the subobjects $A'\subset A$ which satisfy $A/A'\in
\Tscr_J$. An equivalent way of stating lemma \ref{ref:IV.1.1-110}.2 is the following.
\begin{lemma}
The $\Tscr_J$-topology satisfies the Artin-Rees condition.
\end{lemma}
Recall that the Artin-Rees condition means that if $A\subset B \subset D$ with $B/A\in
\Tscr_J$ then there exists $E\subset D$ with $D/E\in \Tscr_J$ and $E\cap B\subset A$.

 The following
technical result will be used repeatedly.
\begin{lemma}
\label{ref:IV.1.3-111} Let $A,D\in\Cscr$ and let $(D_i)_{i\in I}$
be an inverse system of subobjects of $D$ which has the property
that $\cap_i D_i=0$. Then there exists $i\in I$ such that for all
$H\subset D_i$ we have $\Ext^1(H,A)=0$.
\end{lemma}
\begin{proof} By Serre duality we have $\Ext^1(D_i,A)\cong \Hom(E,D_i)^\ast$
  with $E=\tau^{-1}A$. We will show that $\Hom(E,D_i)=0$ for some $i\in
  I$. By the vanishing of $\Ext^2$ this automatically implies that
  $\Ext^1(H,A)=0$ for all $H\subset D_i$. Since $\Hom(E,D_{j})\subset
  \Hom(E,D_i)$ if $D_j\subset D_i$ and $\Hom(E,D_i)$ is finite
  dimensional, it suffices to show that this inclusion is not equality
  for small $D_j$ if $\Hom(E,D_i)\neq 0$. Let $f:E\r D_i$ be a
  non-zero map. Since $\cap_i D_i=0$, there must exist some $l$ such that
  $\im f \not\subset D_{l}$ and hence some $D_j\subset D_i\cap D_l$ such
  that $\im f \not\subset D_{j}$. Thus $f\not\in \Hom(E,D_{j})$ and we
  are done.
\end{proof}
We have the following consequence.
\begin{corollary}
\label{ref:IV.1.4-112}
 The category $\Cscr/\Tscr$ is semisimple.
\end{corollary}
\begin{proof}
  It is sufficient to show that for any $A,D\in \Cscr$ there is a
  subobject $D'\subset D$ with $D/D'\in\Tscr$ such that
  $\Ext^1(D',A)=0$.  To prove this we verify the hypotheses for lemma
  \ref{ref:IV.1.3-111} for the inverse system of cotorsion subobjects of
  $D$, that is, subobjects $D'$ of $D$  with $D/D'\in \Tscr$.
  So let $E$ be a maximal subobject contained in all $D'\subset
  D$ with $D/D'\in\Tscr$.  If $E\neq 0$ then let $E'$ be a maximal
  subobject of $E$. Since $E/E'$ is simple we have $E/E'\in\Tscr$. By
  Artin-Rees we obtain the existence of $D''\subset D$ such that
  $D/D''\in\Tscr$ and $D''\cap E\subset E'$. This contradicts the
  definition of $E$. Hence $E=0$ and we are done.
\end{proof}
The following easy lemma is used several times.
\begin{lemma}
\label{ref:IV.1.5-113}
\begin{enumerate}
\item
Assume that $F\in\Fscr$ is simple modulo $\Tscr$. Then
$\Hom(F,F)\cong k$.
\item Assume that $F,G\in \Fscr$ are simple modulo $\Tscr$. If
  $\Hom(F,G)$ and $\Hom(G,F)$ are both non-zero, then any non-trivial map $\phi:F\r G$
  is an isomorphism.
\end{enumerate}
\end{lemma}
\begin{proof}
\begin{enumerate}
\item
{}{}From the fact that $F$ is torsion free we obtain
 $\Hom(F,F)\hookrightarrow \Hom_{\Cscr/\Tscr}(\pi F,\pi F)$, and since
 $F$ is simple modulo $\Tscr$, the latter is a division algebra. So
 $\Hom(F,F)$ has no zero divisors. Since by our general hypotheses $\Hom(F,F)$ is finite
 dimensional, we obtain that $\Hom(F,F)\cong k$.
\item Pick  non-trivial maps $\theta:G\r F$, $\phi:F\r G$. Then $\phi$ and $\theta$
  represent non-zero maps in $\Hom_{\Cscr/\Tscr}(\pi F,\pi G)$ and
  $\Hom_{\Cscr/\Tscr}(\pi G,\pi F)$. Since $\pi F$ and $\pi G$ are
  simple, these maps must be isomorphisms. In particular the
  compositions $\phi\theta$ and $\theta\phi$ are non-zero in
  $\Hom_{\Cscr/\Tscr}(\pi F,\pi F)$ and $\Hom_{\Cscr/\Tscr}(\pi G,\pi
  G)$, and hence they are also non-zero in $\Hom(F,F)$ and $\Hom(G,G)$.
  It follows from 1. that $\phi\theta$ and $\theta\phi$ are scalar
  multiples of the identity map. Thus $\phi$ is an isomorphism.\qed
\end{enumerate}
\def\qed{}
\end{proof}

Next we  examine more closely the relation between the category $\Tscr$ and the semisimple
category $\Cscr/\Tscr$.

If $D\in\Cscr$ is torsion free and $o$ is a $v$-orbit in $I$, then
we say that $o$ is linked to $D$ if and only if there exists $i\in
o$ together with a surjective map $D\r S_i$.
\begin{lemma}
\label{ref:IV.1.6-114}
\begin{enumerate}
\item If $D$ is linked to $o$ and if $D'$ is isomorphic to $D$ modulo $\Tscr$   then $D'$ is also linked to $o$.
\item If $D$ is linked to $o$ and $i\in o$ then there exists a torsion
  free object $D'$ isomorphic to $D$ modulo $\Tscr_o$ such that there is a
  surjective map $D'\r S_i$.
\end{enumerate}
\end{lemma}
\begin{proof} Assume that there is a surjective map $D\r S_i$ where $i\in o$.
\begin{enumerate}
\item
We can reduce to two cases.
\begin{enumerate}
\item $D\subset D'$ and $D'/D=S_p$ for some simple object $S_p$. Applying $\Hom(-,S_j)$ to the short
  exact sequence
\[
0\r D\r D'\r S_p\r 0
\]
yields
\[
0\r \Hom(S_p,S_j)\r \Hom(D',S_j)\r \Hom(D,S_j)\r \Ext^1(S_p,S_j)\r
0
\]
since $\Ext^1(D',S_j)=0$ using Serre duality. If $v(p)\neq i$ then $\Ext^1(S_p,S_i)=0$ and
hence $\Hom (D',S_i)\neq 0$ since $\Hom (D,S_i)\neq 0$. If $v(p)=i$ then $p\in o$, and we have
$\Hom(D',S_p)\neq 0$.
\item $D'\subset D$ and $D/D'=S_p$. Applying $\Hom(-,S_j)$ to the short
  exact sequence
\[
0\r D'\r D\r S_p\r 0
\]
yields
\[
0\r \Hom(S_p,S_j)\r \Hom(D,S_j)\r \Hom(D',S_j)\r \Ext^1(S_p,S_j)\r 0.
\]
If $p\neq i$ then we take $j=i$ and we have $\Hom(S_p,S_j)=0$. If $p=i$ then we take $j=v(i)$
and we have $\Ext^1(S_p,S_j)\neq 0$. In both cases we find $\Hom(D',S_j)\neq 0$.
\end{enumerate}
\item Assume that there
  is a non-zero  map $\phi:D\r S_j$ (with $j\in o$). If $D'=\ker \phi$ then it follows as
  above that
  $D'$ maps non-trivially to $S_{v(j)}$. By Serre duality we have that
$\Ext^1(S_{v^{-1}(j)},D)\neq 0$. Let $D''$ be a non-trivial extension of $D$ and
$S_{v^{-1}(j)}$. Then $D''$ is torsion free and maps non-trivially to $S_{v^{-1}(j)}$.
\qed
\end{enumerate}
\def\qed{}
\end{proof}
>From the previous lemma we can now deduce the following perhaps
slightly surprising result.
\begin{lemma}
\label{ref:IV.1.7-115}
\begin{enumerate}
\item The operation $\tau^2$ is the identity on objects in $\Cscr/\Tscr$.
\item If $D,D'\in\Fscr$ and $D'=\tau D$ modulo
  $\Tscr$ then $D$ and $D'$ are linked to the same
  $o\in\Tscr/\sim$.
\item If $D,D'\in\Fscr$ are simple modulo $\Tscr$ then $D'=D$ or
  $D'=\tau D$ modulo $\Tscr$ if and only if $D$ and $D'$ are linked to a
  common $o\in\Tscr/\sim$.
\end{enumerate}
\end{lemma}
\begin{proof}
  Clearly 2. and one direction of 3. are trivial. We start by
  proving the non-trivial direction of 3..  So assume that  $D$ and $D'$ in $\Fscr$ are
  simple modulo $\Tscr$ and are linked to a common $o\in I/\sim$.
  By lemma \ref{ref:IV.1.6-114} we may assume that there exists $i\in o$
  together with surjective maps $\phi:D\r S_i$ and $\phi':D'\r S_i$. Let
  $E$ be the pullback of these maps and let $A$ be the kernel of
  $\phi'$. Thus we have a commutative diagram with exact rows.
\begin{equation}
\label{ref:IV.1.2-116}
\begin{CD}
0 @>>> A @>>> D' @>\phi' >> S_i @>>> 0\\ @. @| @AAA @A\phi AA @.\\
0 @>>> A @>>> E @>>> D@>>> 0
\end{CD}
\end{equation}
If the lower exact sequence splits, then $\phi$ factors through $\phi'$. In particular there
is a non-trivial map $\theta :D\r D'$. Since $D'$ and $D$ are simple modulo $\Tscr$, it
follows that $\theta$ must be an isomorphism modulo $\Tscr$ since it is clearly not zero.

Assume now that the lower exact sequence in \eqref{ref:IV.1.2-116} does not split. Then it
yields a non-trivial element of $\Ext^1(D,A)$. Hence by Serre duality $\Hom(\tau^{-1}A,D)\neq
0$. Since $A=D'$ modulo $\Tscr$, we obtain $\tau^{-1}D'=D$ modulo $\Tscr$. This proves what we
want.

Now 1. follows easily from 3.. For let $D\in \Fscr$ be simple modulo $\Tscr$. Then $\tau^2D$
is linked to the same $o\in I/\sim$ as $D$. Hence by 3. we have that $\tau^2D=D$ or
$\tau^2D=\tau D$ modulo $\Tscr$. Thus we find that in any case $\tau^2D=D$ modulo $\Tscr$. This
holds for all simples in $\Cscr/\Tscr$, and since this category is semisimple, it follows that
$\tau^2$ is the identity on objects.
\end{proof}

\begin{corollary}
\label{ref:IV.1.8-117} Assume that is $\Cscr$ connected.  Then
  $\Cscr/\Tscr$
  contains at most two simple objects, and if there are two they are permuted by $\tau$
  modulo $\Tscr$. Furthermore every $A\in \Fscr$ is linked to every
  $o\in I/\sim$.
\end{corollary}
\begin{proof} We first claim that every $o\in I/\sim$ is linked to at
  least one
  object in $\Fscr$. This follows from the fact that
  $\Hom(\Tscr_o,\Tscr_{o'})=\Hom(\Tscr_{o'},\Tscr_o)=\Hom(\Tscr_o,\Fscr)=0$
  for $o'\neq o$. So if $\Hom(\Fscr,\Tscr_o)=0$ then $\Tscr_o$ splits
  off as a factor from $\Cscr$, contradicting the connectedness of
  $\Cscr$.

  Let the simple objects in $\Cscr/\Tscr$ be represented by
  $F_1,\ldots, F_l\in\Fscr$. Assume that $F_1,\ldots,F_t$ constitute one
  $\tau$-orbit modulo $\Tscr$ (thus $t=1,2$). Let $J,K $ be respectively
  the union of the orbits which are linked to $F_1\oplus \cdots\oplus
  F_t$ and $F_{t+1}\oplus \cdots\oplus
  F_l$. By lemma  \ref{ref:IV.1.7-115}(3) we have $J\cap K=\emptyset$ and
  by the previous paragraph $J\cup K=I$.

Define $\Cscr_1$ as the full subcategory of $\Cscr$ consisting of
objects of the form $F\oplus T$ where $F$ is isomorphic to
$F_1^{a_1}\oplus\cdots\oplus F_t^{a_t}$ modulo $\Tscr_J$ and $T\in
\Tscr_J$. We define $\Cscr_2$ in a similar way but using
$F_{t+1},\ldots ,F_l$ and $K$. It is easy to see that
$\Cscr=\Cscr_1\oplus\Cscr_2$. This contradicts again the
hypothesis that $\Cscr$ is connected.

So we obtain that $l=t\le 2$. Since every $o\in I/\sim$ is linked to at least one $F_i$, and
since $\tau$ acts transitively, we obtain that every $o\in I/\sim$ is linked to every $F_i$.
\end{proof}

We will use the following.
\begin{lemma}
\label{ref:IV.1.9-118}
Assume that $\Cscr$ is connected. Let $F\in\Fscr$ and $o\in I/{\sim}$. Then the  $\Tscr_o$ topology on
$F$ is
separated.
\end{lemma}
\begin{proof} This follows from the Artin-Rees property together with
lemma
\ref{ref:IV.1.6-114}. See the proof of Corollary
\ref{ref:IV.1.4-112}.
\end{proof}
\section{Completion}
In this section we use the results from Section 1 to show that if $\Cscr$ is connected, then
the orbits of the simple objects are either all finite or they are all infinite, and in the
first case there is a unique simple object in $\Cscr / \Tscr $. To accomplish this, we use, as
otherwise in this paper, some arguments inspired by geometry.  For example we want to
understand the formal local structure of the objects in $\Cscr$ around a $\tau$-orbit of
simples. For this we use the notion of ``completion'' which was introduced in \cite{VdB19}.

Let $o\in I/{\sim}$. By Appendix A  $\Tscr_{o}$ is a connected
hereditary abelian
category with all objects of  finite length. Furthermore it is clear
that $\Hom$ and $\Ext^1$ in $\Cscr$ and $\Tscr_o$ coincide.
Thus $\Tscr_o$ has a Serre functor and  we can apply the results in \S\ref{ref:III.1-64} to it. We will denote by $C_o$ the
pseudo-compact ring which was denoted by $C$ in \S\ref{ref:III.1-64}.  We consider the
elements of $C_o$ as matrices whose entries are indexed by the elements of $o$. For ${{D}}\in
\Cscr$ we define
\[
\hat{{{D}}}_o=\projlim_{{{D}}/{{D}}'\in \Tscr_o}
({{D}}/{{D}}')\hat{}
\]
where $\hat{(-)}$ was defined in
\S\ref{ref:III.1-64}. If $o$ is clear from the context then
we write $\hat{(-)}$ for $\hat{(-)}_o$.

Clearly $\hat{(-)}$ is a functor from $\Cscr$ to $\PC(C_o)$.  As
in \cite{VdB19} one  proves the following.
\begin{lemma}
\label{ref:IV.2.1-119}
\begin{enumerate}
\item $\hat{(-)}$ is an exact functor. In particular we can extend
  $(\hat{-})$ to a functor $D^b(\Cscr)\r D^b(\PC(C_o))$. This allows us
  to define the action of $(\hat{-})$ on $\Ext$'s.
\item If we restrict $\hat{(-)}$ to $\Tscr$ then it is the
  composition of the projection $\Tscr\r \Tscr_{o}$ with the
  functor $\hat{(-)}$ from \S\ref{ref:III.1-64}.
\item If we restrict $\hat{(-)}$ to $\Fscr$ then its image consists
  of projective pseudo-compact $C_o$-modules containing every
  indecomposable projective pseudo-compact $C_o$-module at most a
  finite number of times as a factor.
\item Let ${{D}}\in \Cscr$ and $T\in \Tscr_{o}$. Then completion
  defines an isomorphism
\begin{equation}
\label{ref:IV.2.1-120} \Ext^i_\Cscr({{D}},T)\cong
\Ext^i_{\PC(C_o)}(\hat{{{D}}},\hat{T})
\end{equation}
\item For every $D\in \Cscr$ we have a natural isomorphism
\[
\widehat{\tau D}\cong \hat{D}_o\otimes_{C_o} \omega_{C_o}
\]
\end{enumerate}
\end{lemma}
\begin{proof}
  This theorem is proved in the same way as \cite[Thm 5.3.1. Prop.\
  5.3.4, Cor.\ 5.3.5]{VdB19}. The only thing that needs to be proved
  slightly differently is 3.

So let $D\in \Fscr$. To prove that $\hat{D}$ is projective we take a minimal projective
resolution
\[
0\r P\r Q\r \hat{D}\r 0.
\]
By the definition of a minimal resolution we have $P\subset
\rad{Q}$. This yields for $j\in o$:
$\Hom_{\PC(C_o)}(P,\hat{S}_j)\cong \Ext^1(\hat{D},\hat{S}_j) \cong
\Ext^1(D,S_j)=0$. Thus $P=0$.

Let $M$ be  indecomposable projective over $C_o$. Then $M/\rad(M)\cong \hat{S}_{j}$ for some
$j$ in $o$. The number of times that $M$ occurs in a direct sum decomposition of  $\hat{D}$ is
given by the dimension of $\Hom_{\PC(C_o)}(\hat{D},\hat{S}_{j})\cong \Hom_{\Cscr}(D,S_j)$,
which is finite by hypothesis.
\end{proof}
Our next aim is to sharpen  lemma \ref{ref:IV.2.1-119}
  by actually showing that $\hat{D}$ is a noetherian $C_o$-module.
We start by proving the following result.
\begin{lemma}
\label{ref:IV.2.2-121} Let ${{D}}\in \Cscr$ and $o\in I/\sim$.
Then there exist only a finite number of $j\in o$ such that there
is a non-zero map ${{D}}\r S_j$.
\end{lemma}
\begin{proof}
Assume that the lemma is false. So ${{D}}$ maps to an infinite
number of different $S_j$ with $j\in o$. In particular $o$ must be
infinite. Without loss of generality we may assume that ${{D}}$ is
torsion free and that $D$ has simple image in $\Cscr/\Tscr$.

Choose an element $i\in o$ and write  $P_l=e_{v^l(i)} C_o$ where for $j\in o$ the
corresponding diagonal idempotent in $C_o$ is denoted by $e_j$. Thus by the explicit structure
of $C_o$ given in \S\ref{ref:III.1-64} we find.
\[
\Hom(P_l,P_{l'})=
\begin{cases}
0 &\text{if $l<l'$}\\ k&\text{otherwise}
\end{cases}
\]
Since by lemma \ref{ref:IV.2.1-119} it follows that  $\hat{D}$ is projective with finite
multiplicities we have
\[
\hat{D}=\prod_{l\in \ZZ} P_{l} ^{a_l}
\]
for certain $l\in\ZZ$. By Proposition \ref{ref:III.1.2-66} we obtain
\begin{equation}
\label{ref:IV.2.2-122} \widehat{\tau^{-1}D}=\prod_{l\in \ZZ} P_{l-1}^{a_l}
\end{equation}
Define $\tau_{l,m}$ as  $P_l/\rad^m(P_l)$ and let $T_{l,m}$ be the
object in $\Tscr_{o}$ with the property that
$\hat{T}_{l,m}=\tau_{l,m}$. Thus $T_{l,m}$ is the unique object of
length $m$ with cosocle $S_{v^l(i)}$. We have
\begin{equation}
\label{ref:IV.2.3-123} \Hom(P_l,\tau_{{l'},m})=
\begin{cases}
k&\text{if $l'\le l\le l'+m-1$}\\ 0&\text{otherwise}
\end{cases}
\end{equation}
Put
\[
T_n=T_{n,1}^{a_n}\oplus T_{n-1,2}^{a_{n-1}}\oplus\cdots\oplus
T_{-n,2n+1}^{a_{- n}}
\]
and $\tau_n=\hat{T}_n$.  Then combining the canonical surjective
maps $P_l\r \tau_{l,m}$ yields a canonical surjective map
$\hat{D}\r \tau_n$ and hence by  lemma \ref{ref:IV.2.1-119} a
surjective map $D\r T_n$.

Let $A_n$ be the kernel of this map. Using the exactness of
completion we find
\begin{equation}
\label{ref:IV.2.4-124} \hat{A}_n=\prod_{l<-n} P_l^{a_l}\oplus
P_{n+1}^{a_{-n}+\cdots+a_n} \oplus \prod_{l>n} P_l^{a_l}
\end{equation}

Applying $\Hom(-,D)$ to the exact sequence
\begin{equation}
\label{ref:IV.2.5-125} 0\r A_n\r {{D}}\r T_n \r 0
\end{equation}
we obtain the exact sequence
\[
0\r \Hom(D,D)\r \Hom(A_n,D)\r \Ext^1(T_n,D)\r \Ext^1(D,D).
\]
Furthermore $\dim \Ext^1(T_n,D)=\dim \Hom(\tau^{-1}D, T_n)=\dim \Hom(\widehat{\tau^{-1}
D},\tau_n)$. This goes to infinity by \eqref{ref:IV.2.2-122}, \eqref{ref:IV.2.3-123} and the
definition of $\tau_n$.
 It follows that $\dim \Hom(A_n,D)$ also goes to
infinity. On the other hand we now show that in fact $\Hom(A_n,D)\cong k$. This is clearly a
contradiction.

We apply $\Hom(A_n,-)$ to the exact sequence
\eqref{ref:IV.2.5-125}. This yields the exact sequence
\[
0\r \Hom(A_n,A_n)\r \Hom(A_n,D)\r \Hom(A_n,T_n)
\]
{}{}From the fact that $A_n$ is simple modulo $\Tscr$ we obtain by lemma
 \ref{ref:IV.1.5-113}  that $\Hom(A_n,A_n)\cong k$.

 On the
other hand, by completing and using \eqref{ref:IV.2.4-124} and
\eqref{ref:IV.2.3-123} we find $\Hom(A_n,T_n)=0$. This proves what
we want.
\end{proof}
\begin{corollary} Denote by $\pc(C_o)$ the full subcategory of
  $\PC(C_o)$ consisting of noetherian objects. Then $(-)\hat{}_o$
  defines a functor $\Cscr\r \pc(C_o)$.
\end{corollary}
\begin{proof} It clearly suffices to show that if $F\in\Fscr$ then
  $\hat{F}$ is noetherian. Now we know that $\hat{F}$ is a projective
  object in $\PC(C_o)$
  mapping only to a finite number of different simples (counting
  multiplicities). Hence $\hat{F}$ is a direct sum of a finite number
  of indecomposable projectives. Since the indecomposable
  pseudocompact projectives
  over $C_o$ are noetherian, it follows that $\hat{F}$ is noetherian.
\end{proof}
Note that if $|o|=\infty$ then $\pc(C_o)$ is equivalent to $\gr(R)$ with $R=k[x]$.
\begin{corollary}
\label{ref:IV.2.4-126} Let ${{D}}\in \Cscr$ and $T\in
\Tscr_{o}$. Then completion defines an isomorphism
\begin{equation}
\label{ref:IV.2.6-127} \Ext^i_\Cscr({{T}},D)\r
\Ext^i_{\PC(C_o)}(\hat{{{T}}},\hat{D})
\end{equation}
\end{corollary}
\begin{proof} Let $\eta_T:\Ext^1(T,\tau T)\r k$ be the linear map
  corresponding to the identity map $T\r T$ under Serre duality. As in
  the proof of Theorem \ref{ref:I.2.3-15} it follows that $\eta_T$ can be
  chosen freely, subject to the condition that it must be
  non-vanishing on the almost split sequence ending in $T$.

Now define $\eta_{\hat{T}}:\Ext^1(\hat{T},\tau\hat{T})\r k$ as
$\hat{\eta}_T$.
 Then by
    functoriality we have a commutative diagram:
$$ \xymatrix{ \Ext^i(T,D)\ar[d]^{(\hat{-})}&\times &\Ext^{1-i}(D,\tau T)
\ar[d]^{(\hat{-})}\ar[r]&\Ext^1(T,\tau T)\ar[r]^-{\eta_T}\ar[d]^{(\hat{-})}&k\\
\Ext^i(\hat{T},\hat{D})&\times &\Ext^{1-i}(\hat{D},\tau \hat{T}) \ar[r]&\Ext^1(\hat{T},\tau
\hat{T})\ar[r]^-{\eta_{\hat{T}}}&k } $$ It follows that the pairings for $\Cscr$ and for
$\PC(C_o)$ are compatible. The upper pairing is non-degenerate by Serre duality and the lower
pairing is non-degenerate by local duality for $C_o$ (or graded local duality for $k[x]$ if
$|o|=\infty$). Using these dualities \eqref{ref:IV.2.6-127} follows from
\eqref{ref:IV.2.1-120}.
\end{proof}

Let $\dis(C_o)$ be the pseudo-compact finite length modules over
$C_o$ and  put $\qpc(C_o)=\pc(C_o)/\dis(C_o)$.  Using Morita
theory we obtain
\begin{equation}
\label{ref:IV.2.7-128} \qpc(C_o)\cong
\begin{cases}
 \mod((C_o)_x)\cong \mod{k((x))}&\text{if $|o|<\infty$}\\
\gr(k[x,x^{-1}])\cong \mod(k)&\text{if $|o|=\infty$}
\end{cases}
\end{equation}

Then completion defines an exact functor $\Cscr/\Tscr_o\r \qpc(C_o)$. So we obtain the
following commutative diagram
\begin{equation}
\label{ref:IV.2.8-129}
\begin{CD}
\Cscr/\Tscr_o @>\hat{(-)}>> \qpc(C_o)\\ @A\pi AA @A\pi AA\\ \Cscr
@>\hat{(-)}>> \pc(C_o)
\end{CD}
\end{equation}
where as usual $\pi$ denotes the quotient functor. In the next
section we will show that this is actually a pullback diagram.

\begin{lemma}
\label{ref:IV.2.5-130}
Let $D,E\in \Cscr$ and assume that $o\in I/\sim$ is such
  that the linear topology on $E$ defined by $\Tscr_o$ is
  separated. Then the canonical map
\[
\Hom_{\Cscr}(D,E)\r \Hom_{\PC(C_o)}(\hat{D},\hat{E})
\]
is injective.
\end{lemma}
\begin{proof}
This is trivial.
\end{proof}
\begin{corollary}
\label{ref:IV.2.6-131} Assume that $F\in\Fscr$ is simple modulo
$\Tscr$. Assume
  furthermore that $F$ is linked to an infinite orbit $o\in
  I/\sim$.
  Then $\Hom_{\Cscr/\Tscr}(\pi F,\pi F)\cong k$.
\end{corollary}
\begin{proof}  By hypotheses $\Hom_{\Cscr/\Tscr}(\pi F,\pi F)$ is a
  skew field.  Since $k$ is algebraically closed, it now suffices to show
  that $\Hom_{\Cscr/\Tscr}(\pi F,\pi F)$ is finite dimensional.

We have
\[
\Hom_{\Cscr/\Tscr}(\pi F,\pi
F)=\dirlim_{F/F_1\in\Tscr}\Hom_\Cscr(F_1, F)
\]
 By lemma
  \ref{ref:IV.2.5-130} the $\Tscr_o$ topology on $F_1$ and $F$ is
  separated. Hence
\[
\dirlim_{F/F_1\in\Tscr}\Hom_\Cscr(F_1, F)
\subset \dirlim_{F/F_1\in\Tscr}
  \Hom_{\pc(C_o)}(\hat{F}_1,\hat{F})
\]
If $F_1$ runs through the subobjects in $F$ with the property that
$F/F_1\in\Tscr$ then $\hat{F}_1$ runs through the subobjects $G$ of
$\hat{F}$
such that $\hat{F}/G\in \dis(C_o)$. Thus
\begin{align*}
\dirlim_{F/F_1\in\Tscr}
  \Hom_{\pc(C_o)}(\hat{F}_1,\hat{F})
  &= \dirlim_{\substack{G\in\pc(C_o)\\\hat{F}/G\in\dis(C_o)}}
  \Hom_{\pc(C_o)}(G,\hat{F})\\
  &=\Hom_{\qpc(C_o)}(\pi\hat{F},\pi\hat{F})
\end{align*}
By \eqref{ref:IV.2.7-128} this is finite dimensional.
\end{proof}
By contrast the situation in case of linkage to a finite orbit is
different.
\begin{lemma}
\label{ref:IV.2.7-132} Assume that $F,F'\in\Fscr$ are simple modulo $\Tscr$. Assume
furthermore that $F$ and $F'$ are linked to a finite orbit $o\in I/\sim$. Then
$\Hom_{\Cscr/\Tscr}(\pi F,\pi F')$ is infinite dimensional.
\end{lemma}
\begin{proof} Since $F$ and $F'$ are torsion free one has that
  $\Hom_{\Cscr/\Tscr}(\pi F,\pi F')$ is equal to the union of
  $\Hom(F_1,F')$ where $F/F_1\in \Tscr$. Hence it suffices to show
  that the dimension of $\Hom(F_1,F')$ goes to infinity if $F_1$ runs
  through the subobjects of $F$ satisfying $F/F_1\in\Tscr_o$. Put
  $T=F/F_1$. Applying $\Hom(-,F')$ to the exact sequence
\[
0\r F_1\r F\r T\r 0
\]
yields the exact sequence
\[
0\r \Hom(F,F')\r \Hom(F_1,F')\r \Ext^1(T,F')\r \Ext^1(F,F')
\]
Hence it suffices to show that the dimension of $\Ext^1(T,F')$ goes to infinity. By Serre
duality and completion we must show that the dimension of $\Hom(\hat{F'},\widehat{\tau T})$
goes to infinity. To this end it is sufficient to show that for indecomposable projectives $P$
and $ P'$ over $C_o$  the dimension of $\Hom_{C_o}(P', S\otimes\omega_{C_o})$ goes to infinity
where $S$ runs through the finite dimensional quotients of $P$. This is an easy direct
verification using the explicit structure of $C_o$ given
 in \S\ref{ref:III.1-64}.
\end{proof}
\begin{corollary}
\label{ref:IV.2.8-133} Assume that $\Cscr$ is connected. Then the
$v$-orbits in $I$ are either all finite, or all infinite.
Furthermore if they are finite then there exists only one simple
object in $\Cscr/\Tscr$.
\end{corollary}
\begin{proof} Assume that there are both finite and infinite
  orbits. Then by Corollary \ref{ref:IV.1.8-117} there exists
  $F\in \Fscr$, simple in $\Cscr/\Tscr$, such that $F$ is linked to
  both a finite and an infinite orbit. This is a contradiction
  according to Corollary \ref{ref:IV.2.6-131} and lemma \ref{ref:IV.2.7-132}.

  Assume now that the orbits are finite but that there exist objects
  $F,F'\in\Fscr$ which are
  distinct and simple modulo $\Tscr$.
  Clearly $\Hom_{\Cscr/\Tscr}(\pi F,\pi F')=0$. On the other
  hand by Corollary \ref{ref:IV.1.8-117} we know that $F$ and $F'$ are both
  linked to the same finite orbit, whence by lemma \ref{ref:IV.2.7-132}
  it follows that $\Hom_{\Cscr/\Tscr}(\pi F,\pi F')$ is infinite dimensional. This is
  clearly a contradiction.
\end{proof}

\section{Description in terms of a pullback diagram}
\label{ref:IV.3-134} In this section we give useful descriptions
of $\Cscr$ which allow us to construct autoequivalences of $\Cscr$
with desired properties. This will be used for the construction of
ample pair in the next section.

We start with the following description of $\Cscr$.

\begin{lemma}
\label{ref:IV.3.1-135} \eqref{ref:IV.2.8-129} is a pullback
diagram in the sense that
  $\Cscr$ is equivalent to the category of triples
$(F,Z,\phi)$ where $F\in \Cscr/\Tscr_o$ and  $Z\in \pc(C_o)$ and $\phi$ is an isomorphism
$\hat{F}\r \pi Z$.
\end{lemma}
\begin{proof} Let $\Cscr'$ be the category of triples defined in the
  statement of the lemma. Define
  $U:\Cscr\r\Cscr'$ by $UF=({{\pi}}F,\hat{F},I)$ where
  $I$ is the natural isomorphism ${{\pi}}\hat{F}\r ({{\pi}}F)\hat{}$. We have to
  show that $U$ is an equivalence. To do this we show that $U$ is
  respectively faithful, full and essentially surjective.

\noindent {\bf Faithfulness : } Assume that $F_1,F_2\in\Cscr$ and that $f:F_1\r F_2$ is a
homomorphism such that $Uf=0$. Thus ${{\pi}}f=0$ and $\hat{f}=0$. The fact that ${{\pi}}f=0$
means that $f$ factors as $F_1\xrightarrow{f_1} T\xrightarrow{f_2} F_2$ where $T\in\Tscr_o$,
the map $f_1$ is an epimorphism and $f_2$ is a monomorphism. Then $0=\hat{f}$ factors
likewise. This then implies that $\hat{T}=0$ and hence $T=0$. We conclude $f=0$.

\noindent {\bf Fullness : } Let $F_1$ and $F_2$ be as in the previous paragraph. Now assume
that we are given maps $f':{{\pi}}F_1\r {{\pi}}F_2$ and $f'':\hat{F}_1\r \hat{F}_2$ such that
$\hat{f'}=\pi f''$. We have to produce a map $f:F_1\r F_2$ such that $f'={{\pi}}f$ and
$f''=\hat{f}$. Let $g:G_1\r F_2$ be a representative of $f'$ in $\Cscr$ where $G_1$ is a
subobject of $F_1$ such that $F_1/G_1\in\Tscr_o$. Below we denote the inclusion map $G_1\r
F_1$ by $\alpha$.

We now have ${{\pi}}\hat{g}={{\pi}}f''$. This means that there
exist $G'_1\subset \hat{G}_1$ with $\hat{G}_1/G'_1\in \dis(C_o)$
such that the compositions $G'_1\r \hat{G}_1\xrightarrow{\hat{g}}
\hat{F}_2$ and $G'_1\r \hat{G}_1\r \hat{F}_1\xrightarrow{f''}
\hat{F}_2$ are equal.

{}{}From \eqref{ref:IV.2.1-120} it follows easily that there exists $K_1$ which is a subobject
of $G_1$ with cokernel in $\Tscr_o$ such that $\hat{K}_1=G'_1$. We now replace $G_1$ by $K_1$,
and $g$ by the composition $K_1\r G_1\xrightarrow{g} F_2$. Hence we obtain the following
commutative diagram
\[
\begin{CD}
\hat{G}_1 @>\hat{g} >> \hat{F}_2\\ @VVV @|\\ \hat{F}_1 @>f'' >>
\hat{F}_2
\end{CD}
\]
So our problem is now to lift $g$ to a homomorphism $f:F_1\r F_2$ such that $\hat{f}=f''$. Let
$H$ be the pushout of $G_1\xrightarrow{g} F_2$ and $G_1\xrightarrow{\alpha} F_1$. Then we have
a commutative diagram with exact rows
\[
\begin{CD}
0 @>>> G_1 @>\alpha>> F_1 @>>> T_1 @>>> 0\\ @. @Vg VV @VVV @VVV
@.\\ 0 @>>> F_2 @>>> H @>\gamma>> T_1 @>>> 0
\end{CD}
\]
with $T_1\in\Tscr_o$. It is easy to see that factorizations of $g$
through $\alpha$ are in one-one correspondence with splittings of
$\gamma$. Similarly factorizations of $\hat{g}$ through
$\hat{\alpha}$ correspond to splittings of $\hat{\gamma}$.

Now we know that $\hat{g}$ factors through $\hat{\alpha}$ and
hence this yields a splitting of $\hat{\gamma}$. By Corollary
\ref{ref:IV.2.4-126} and \eqref{ref:IV.2.1-120} this splitting
corresponds to a splitting of $\gamma$ and thus to a factorization
$f$ of $g$ through $\alpha$. It is now easy to see that this is
the $f$ we are looking for.

\noindent {\bf Essential surjectivity : } Let $(F,Z,\phi)\in \Cscr'$. We have to show that
this is in the essential image of $U$. Choose $G\in \Cscr$ representing $F$. By definition
$\phi$ is now an isomorphism ${{\pi}}\hat{G} \r {{\pi}}Z$.

As usual $\phi$ is represented by a injective map  $\theta':K'_1\r
Z$ with cokernel in $\dis(C_o)$ where $K'_1$ is a subobject of
$\hat{G}$ also with cokernel in $\dis(C_o)$.

  Thus we have the following maps
\begin{equation}
\label{ref:IV.3.1-136} \hat{G}\l K'_1 \xrightarrow{\theta'} Z
\end{equation}
By working from left to right and employing
\eqref{ref:IV.2.1-120}\eqref{ref:IV.2.6-127}  these arrows and
objects can be ``uncompleted''. That is, there are arrows and
objects in $\Cscr$
\begin{equation}
\label{ref:IV.3.2-137} G\l K_1 \xrightarrow{\theta''} L
\end{equation}
such that completion gives \eqref{ref:IV.3.1-136}.

Now $\hat{L}=Z$ and $\pi L\cong \pi G\cong F$. Checking the appropriate
commutative diagrams we find that indeed $UL\cong (F,Z,\phi)$.
\end{proof}
We may use the previous lemma to define for any orbit $o$ in $I$ a
canonical associated autoequivalence. Consider the following
commutative diagram of functors.
\[
\begin{CD}
\Cscr/\Tscr_o @>\hat{(-)}>> \qpc(C_o) @<\pi << \pc(C_o)\\ @V \Id VV
@V \Id VV @V \rad(C_o)\otimes_{C_o} - VV\\ \Cscr/\Tscr_o
@>\hat{(-)}>> \qpc(C_o) @<\pi << \pc(C_o)
\end{CD}
\]
The vertical arrows are clearly autoequivalences, and hence, using
lemma \ref{ref:IV.3.1-135}, they define an autoequivalence on
$\Cscr$ which we will denote by $s_o$.

The functor $s_o$ has a more direct description on torsion free
objects.
\begin{lemma} Assume that $F\in\Fscr$. Then $s_o F$ is the kernel of
  the universal map
\begin{equation}
\label{ref:IV.3.3-138} F\r \oplus_{i\in o} S_i\otimes
\Hom(F,S_i)^\ast
\end{equation}
which in addition is surjective.
\end{lemma}
\begin{proof}
Surjectivity of \eqref{ref:IV.3.3-138} follows from the fact that
there are only a finite number of $i\in o$ such that
$\Hom(F,S_i)\neq 0$ (see  lemma
\ref{ref:IV.2.2-121}).

To prove that \eqref{ref:IV.3.3-138} gives the correct result for
$s_o$ on $\Cscr$ one has to check that it gives the correct result
on $\Cscr/\Tscr_o$, $\pc(C_o)$ and $\qpc(C_o)$. In each of the
cases this is clear.
\end{proof}
{}{}From Proposition
\ref{ref:III.1.2-66} we also obtain the
following obvious fact.
\begin{lemma} The functor $s_o$ is (non-canonically) isomorphic to $\tau$
  when restricted
  to $\Tscr_o$. It is the identity on $\Tscr_{o'}$ for $o'\neq o$.  In
  particular
\[
s_o(S_i)\cong
\begin{cases}
S_{v(i)}&\text{if $i\in o$}\\ S_i&\text{if $i\not\in o$}
\end{cases}
\]
\end{lemma}
The foregoing results have an obvious extension in the case of a
finite number of $v$ orbits
 $o_1,\ldots,o_n$. To state this let
$J$ be the union of these orbits. Then there is a pullback diagram
\begin{equation}
\label{ref:IV.3.4-139}
\begin{CD}
\Cscr/\Tscr_J @>\oplus_i \hat{(-)}_{o_i}>> \oplus_i
\qpc(C_{o_i})\\ @A\pi AA @A\pi AA\\ \Cscr @>\hat{(-)}_{o_i}>>
\oplus_i\pc(C_o)
\end{CD}
\end{equation}
In the same way as for a single orbit we can define an associated
autoequivalence $s_J$. It is easy to see that $s_J=s_{o_1}\cdots
s_{o_n}$.  Furthermore if $F\in\Fscr$ then $s_J$ is given by the
kernel of the surjective universal map
\begin{equation}
\label{ref:IV.3.5-140} F\r \oplus_{i\in J} S_i\otimes
\Hom(F,S_i)^\ast
\end{equation}
Finally on simples we have:
\begin{equation}
\label{ref:IV.3.6-141} s_J(S_i)\cong
\begin{cases}
S_{v(i)}&\text{if $i\in J$}\\ S_i&\text{if $i\not\in J$}
\end{cases}
\end{equation}

\section{The finite orbit case}
\label{ref:IV.4-142} In this section we complete our
classification in the finite orbit case, by constructing an ample
pair $(F,t)$ and applying results from \cite{AZ}.

If $\Dscr$ is a noetherian $\Ext$-finite abelian category and $(F,t)$ is a pair
consisting of an object  $F\in\Dscr$ and an autoequivalence $t$ of
$\Dscr$, then following \cite{AZ} we will say that this pair is
ample if the following conditions hold.
\begin{itemize}
\item[(a)] For every $D\in\Dscr$ there exists a
  surjective map $\oplus_{i=1}^m t^{-n_i} F\r D$ with $n_i\ge 0$.
\item[(b)] For every surjective map $E\r D$ in $\Dscr$   the
  induced map $\Hom(t^{-n}F, E)\r \Hom(t^{-n} F,D)$ is surjective for
  large $n$.
\end{itemize}
Associated to a pair $(F,t)$ one can canonically associate a
graded ring $\Gamma^\ast(F)$ defined by
\[
\Gamma^\ast(F)=\bigoplus_{n\in\ZZ} \Hom(t^{-n}F, F)
\]
It is shown in \cite{AZ} that if $(F,t)$ is ample then
$\Gamma^\ast(F)_{\ge 0}$ is a noetherian ring and $\Dscr$ is
equivalent to $\qgr(\Gamma^\ast(F)_{\ge 0})$.

In this section we impose the usual hypotheses on $\Cscr$, but we
assume in addition that $\Cscr$ is connected and that all
$v$-orbits in $I$ are finite.

By Corollary \ref{ref:IV.2.8-133} there is exactly one simple object modulo $\Tscr$. Assume it
is represented by $F\in\Fscr$. Since $\tau F=F$ modulo $\Tscr$, there must exist a finite
number of orbits $o_1,\ldots,o_n$ in $I$ such that $\tau F=F$ modulo $\Tscr_J$ with
$J=o_1\cup\cdots \cup o_n$.  Define $s_J$ as in Section \ref{ref:IV.3-134}. Writing $s=s_J$
and  defining  $t=s^{-1}$, we have the following.
\begin{lemma}
\label{ref:IV.4.1-143} $(F,t)$ is an ample pair in $\Cscr$.
\end{lemma}
\begin{proof} The $\Tscr_J$-topology on $F$ is separated by
lemma \ref{ref:IV.1.9-118}. So by
  lemma \ref{ref:IV.1.3-111} there exists a subobject $F_1\subset F$
  with $F/F_1\in\Tscr_J$ and with the property $\Ext^1(F',F)=0$ for
  every subobject $F'$ of $F_1$. By \eqref{ref:IV.3.5-140} there
  exists an $n_0$ such that $s^nF \subset F_1$ for $n\ge n_0$. Hence
  we have
\[
\Ext^1(t^{-n}F,F)=\Ext^1(s^n F,F)=0
\]
for $n\gg 0$. {}From this we easily obtain that if we can show (a) then for
all $X\in \Dscr$ we have $\Ext^1(t^{-n} F,X)=0$ for $n\gg 0$. Thus
(b) follows.

Now let $D\in \Cscr$ and let $E$ be the largest subobject of $D$
which can be written as a quotient of a direct sum $\oplus_{i=1}^t
s^{n_i} F$ with $n_i \geq 0$. Assume $D\neq E$ and write $T=D/E$.
As above there exists $n_0$ such that $\Ext^1(s^nF,E)=0$ for $n\ge
n_0$. On the other hand we claim that there exists a non-trivial
map $s^m F\r T$ for some $m\ge n_0$.
 By the vanishing of
$\Ext^1(s^mF,E)$ this map then lifts to a non-trivial map $s^mF\r
D$, contradicting the choice of $E$.

It remains to show that there is a non-trivial map $ s^m F\r T$ for
some $m\geq n_0$. Assume first that $T$ contains a torsion free
subobject $T'$. We may assume that $T'$ is simple modulo $\Tscr$.  Then the
proof of Lemma \ref{ref:IV.2.7-132} together with Corollary
\ref{ref:IV.1.8-117} shows that there is a non-trivial map $F'\r T'$
with $F/F'\in\Tscr_o$. As above we have that there is some $s^n
F\subset F'$. Since $F/s^n F$ is torsion the composition $s^n F\r F'\r
T'$ cannot be zero. This proves what we want.

Assume now that $T$ is torsion and let $S_i$ for $i\in I$ be a simple subobject of $T$ and let $o$ be the $v$-orbit of $i$. Then
according to lemma \ref{ref:IV.1.6-114} there exists $j\in o$ together with a non-trivial map
$F\r S_j$. We now consider two possibilities:
\begin{enumerate}
\item $o\not\subset J$. Since $F= \tau^nF$ modulo  $\Tscr_J$ for all $n$, it
  follows that there will also be a non-trivial map $\tau^n F\r S_j$ for
  all $n$, whence a non-trivial map $F\r S_{v^{-n}(j)}$. Choose $n$ in
  such a way that $v^{-n}j=i$ and apply $s^m$. This yields
  a non-trivial map $s^mF\r S_i$ for all $m$.
\item $o\subset J$. Now we  have a surjective map $s^nF\r
  S_{v^n(j)}$. Since the $v$-orbit of $j$ is finite, it follows that we
  can always find $m\ge n_0$ such that $s^mF$ maps surjectively to
  $S_i$.
\end{enumerate}
Hence in both cases we obtain a non-trivial map $s^mF\r S_i$ with
$m\ge n_0$. By the above discussion we are done.
\end{proof}
What remains to be done in order to  determine $\Cscr$ explicitly
is to compute $R=\Gamma^\ast(F)_{\ge 0}$. We do this in a number of steps.
\begin{step} $R$ is positively graded and $R_0=k$. Indeed from the
  fact that $F$ is simple modulo $\Tscr$ it follows that
  $\Hom(F,F)\cong k$ (see lemma
\ref{ref:IV.1.5-113}).
\end{step}
\begin{step} $R$ is a domain. This follows from the fact that $R_n$
  embeds in $\Hom_{\Cscr/\Tscr}(\pi F,\pi F)$, and since $F$ is
  simple modulo $\Tscr$, the latter is a division algebra $D$. The
  multiplication $R_n\otimes R_m\r R_{m+n}$ corresponds to
  multiplication in $D$. Hence we obtain an inclusion
  $R\hookrightarrow D[c,c^{-1}]$ where $c$ represents the inclusion
  $sF\r F$. Since $c$ also represents the identity element in $
  \Hom_{\Cscr/\Tscr}(\pi F,\pi F)$,   we find in addition that $c$ is a
  central non-zero divisor in $R$.
\end{step}
\begin{step} $R$ has Gelfand-Kirillov dimension two. To see this apply
  $\Hom(-,F)$ to the exact sequence
\[
0\r s^{n+1} F\r s^n F \r s^n(F/sF)\r 0
\]
Then we obtain an inclusion $R(-1)\xrightarrow{c} R$ with cokernel contained in $\oplus_n
\Ext^1(s^n(F/sF),F)$. Now by Serre duality we have $\dim \Ext^1(s^n(F/sF),F)=\dim
\Hom(\tau^{-1}F,s^n(F/sF))$. Since by \eqref{ref:IV.3.5-140}, $F/sF=\oplus_{i\in J} S_i\otimes
\Hom(F,S_i)^\ast$ we find  that $\Hom(\tau^{-1}F,s^n(F/sF))\cong \oplus_{i\in
  J}\Hom(F,S_{v^{n+1}(i)})\otimes  \Hom(F,S_i)^\ast$. Thus the
dimension of $\Hom(\tau^{-1}F,s^n(F/sF))$ is bounded independently of $n$ (using the fact that
$s$ has finite order).
\end{step}
\begin{step} $R$ is commutative.
To see this we compute $D=\Hom_{\Cscr /\Tscr}(\pi F,\pi F)$ using the equivalence
$\qgr(R)=\Cscr$. We denote the quotient map $\gr(R)\r \qgr(R)$ also by $\pi$. We find
\begin{align*}
\Hom_{\Cscr/\Tscr}(\pi F,\pi F)&=\dirlim_{0\neq G\subset
  F}\Hom_\Cscr(G,F)\\
&= \dirlim_{0\neq I\subset
  R}\Hom_{\qgr(R)}(\pi I,\pi R)\\
&=\dirlim_{0\neq J\subset
  R}\Hom_{\gr(R)}(J,R)
\end{align*}
Hence it follows that $D$ is equal to the degree zero part of the
graded quotient field of $R$. Since $R$
  has GK-dimension 2, it follows from \cite{Staf5} that $D$ is
  commutative. Since $R$ is included in $D[c,c^{-1}]$, it follows that
  $R$ is also commutative.
\end{step}
Using Proposition \ref{ref:III.2.3-71} together with the fact that $\Cscr$ has no nonzero
projective objects, we have now proved the following.
\begin{theorem} \label{ref:IV.4.2-144}
  Assume that $\Cscr$ is a connected $\Ext$-finite noetherian
  hereditary abelian category with Serre functor and no nonzero projective
  objects, and assume that the $\tau$-orbits of the simple objects are
  finite. Then $\Cscr$ is equivalent to the category $\coh \Oscr$ of
  coherent sheaves over a sheaf $\Oscr$ of hereditary orders over a
  non-singular connected projective curve.
\end{theorem}
This theorem  is sufficient for our classification result. However with only a little more
work one can prove the following refinement.  We leave the proof to the reader.
\begin{theorem}\label{ref:IV.4.3-145}
The following are equivalent:
\begin{itemize}
\item[(a)]
 $\Cscr$ is a connected
$\Ext$-finite noetherian hereditary abelian category with Serre functor and no nonzero
projective objects, and  the orbits of the simple objects are finite.
\item[(b)] $\Cscr$ is equivalent to $\qgr(R)$ where
  $R=k+R_1+R_2+\cdots$ is a finitely generated commutative domain of
  Krull dimension two,
  $R_1\neq 0$ and $R$ is an isolated singularity.
\item[(c)]
$\Cscr$ is equivalent to the category $\coh \Oscr$
 of coherent sheaves over a sheaf $\Oscr$ of hereditary orders over a
 non-singular connected projective curve.\qed
\end{itemize}
\def\qed{}\end{theorem}

\section{The infinite orbit case}
\label{ref:IV.5-146} Let $\Cscr$ satisfy the usual hypotheses and assume in addition that
$\Cscr$ is connected and that all $v$-orbits in $I$ are infinite. We show that there are
\emph{at most two} hereditary abelian categories with these properties and that these
categories (if they exist) are distinguished by the property of having one or two
$\tau$-orbits of simples. Since the $\ZZ A_\infty^\infty$ and $\ZZ D_\infty$-categories are
examples of hereditary abelian categories with respectively two and one orbits of simples, the
proof of the main result in this chapter (Theorem \ref{ref:IV.5.2-152}) is then complete.

We refer to Remark \ref{ref:IV.5.1-151} below for other approaches
towards finishing the proof of Theorem \ref{ref:IV.5.2-152}.

Following \cite{SmithZhang} it will be convenient to use injective
resolutions of objects in $\Cscr$. Therefore following
\cite{Gabriel} we let $\tilde{\Cscr}$ be the closure of $\Cscr$
under direct limits. Then $\tilde{\Cscr}$ is also
hereditary according to Appendix A.

As before we will separate our analysis in two cases.

\noindent \textbf{Case I: One simple object modulo $\Tscr$} \\ Let $F\in\Fscr$ be simple
modulo $\Tscr$. Since $\Cscr/\Tscr$ is semisimple and $\tilde{\Cscr}$ is hereditary, the minimal
injective resolution of $F$ has the form
\begin{equation}
\label{ref:IV.5.1-147} 0\r F\r E_1\r E_2\r 0
\end{equation}
where $E_1$ is the  injective hull of $F$ and
$E_2\in\tilde{\Tscr}$.

If $E(S_i)$ denotes the injective hull of $S_i$, then
$E_2=\oplus_{i\in
  I} E(S_i)^{a_i}$, where
$a_i=\dim \Ext^1(S_i,F)$.

Now choose $0\neq T\hookrightarrow E_2$ in such way that $T\in\Tscr$. Let $G$ be the inverse
image of $T$ in $E_1$. Clearly $G\in\Fscr$. Applying $\Hom(G,-)$ to \eqref{ref:IV.5.1-147} we
find the exact sequence
\[
0\r \Hom(G,F)\r \Hom(G,E_1)\r \Hom(G,E_2)\r \Ext^1(G,F)\r 0.
\]
Since $F$ and $G$ are simple modulo $\Tscr$ and  there exists a non-isomorphism $F\r G$ (the
inclusion), it follows from lemma \ref{ref:IV.1.5-113} that $\Hom(G,F)=0$.
 Furthermore by Corollary
\ref{ref:IV.2.6-131} we also have $\Hom(G,E_1)=k$. Finally by Serre duality we have $\dim
\Ext^1(G,F)=\dim \Hom(\tau^{-1}F,G)$.

Now since $\tau^{-1}F= F$ modulo $\Tscr$, it follows that $E_1$ is also the injective hull of
$\tau^{-1}F$. Thus there is an inclusion $\tau^{-1}F\r E_1$. By choosing $T$, and hence $G$,
large enough we may assume $\tau^{-1}F\subset G$, and hence by Corollary
\ref{ref:IV.2.6-131} we have $\Hom(\tau^{-1}F,G)\cong k$.

 Thus for $G$ large we must have
\begin{equation}
\label{ref:IV.5.2-148} \dim \Hom(T,E_2)\le \dim \Hom(G,E_2)=\dim
\Hom(G,E_1)+\dim \Ext^1(G,F)=2
\end{equation}
If $o$ is a $v$-orbit in $I$ then let $b_o=\sum_{i\in o}a_i$. Note
 that by lemma \ref{ref:IV.2.2-121} and Serre duality it follows that
 $b_o$ is finite.

{}{}From the description of the
category $\Tscr_o$ as $\mathrm{f.l. gr}(k[x])$ we obtain that we
can make $\dim \Hom(T,E_2)$ equal to $\sum_{o\in
  S} b_o^2$ for any finite subset $S$ of orbits in $I$. In particular
 there are at most two orbits for which $b_o\neq 0$ and for these orbits
 $b_o=1$.  Said differently, $F$ has at most two simple quotients
 (using Serre duality), and
 if there are two then they must lie in distinct $\tau$-orbits.

Assume now that $F$ has exactly one simple quotient $S_i$.
Applying $\Hom(-,E_2)$ to the exact sequence
\[
0\r F\r G\r T\r 0
\]
we find the exact sequence
\[
0\r \Hom(T,E_2)\r \Hom(G,E_2)\r \Hom(F,E_2)\r 0.
\]
 By Serre
duality $F$ has a unique non-trivial extension by the simple $S_{v^{-1}i}$. By the structure
of $\Tscr$ it follows that $E_2$ is uniserial with subquotients $S_{v^{-n}i}$ for $n>0$,
whence $\Hom(F,E_2)=0$. But again using the structure of $\Tscr$ we also have
$\Hom(T,E_2)\cong k$. Thus $\Hom(G,E_2)\cong k$, contradicting \eqref{ref:IV.5.2-148}.

So $F$ has two simple quotients which lie in different $\tau$-orbits. Translating
\eqref{ref:IV.3.4-139} to our situation using the description of $\pc(C_o)=\gr(k[x])$ and
$\qpc(C_o)=\mod(k)$ we find that $\Cscr$ is the pullback of
\begin{equation}
\label{ref:IV.5.3-149}
\begin{CD}
\mod(k) @>>> \mod(k)\oplus \mod(k)\\ @. @AAA\\ @.\gr(k[u])\oplus
\gr(k[v])
\end{CD}
\end{equation}
where the horizontal map is the diagonal map and the vertical map is obtained from localizing
at $u$ and $v$ and restricting to degree zero. Hence it follows that $\Cscr$ is indeed
determined up to equivalence.

\textbf{Case II : Two simple objects modulo $\Tscr$}\\ Let
$F_1,F_2\in\Fscr$ be the representatives of the simple objects in
$\Cscr/\Tscr$. As in Case I we start with the minimal injective
resolution of $F_1$:
\[
0\r F_1 \r E_1\r E_2\r 0
\]
With a similar reasoning as in Case I, but now using that by Serre duality $\Ext^1(G,F_1)=0$
for $G$ torsion free and isomorphic to $F_1$ modulo $\Tscr$, we find that $F_1$ has a unique
simple quotient. The same holds for $F_2$, and by Corollary \ref{ref:IV.1.8-117} these simple
quotients must be in the same $\tau$-orbit.

Now translating \eqref{ref:IV.3.4-139} to our situation we find
that $\Cscr$ is the pullback of
\begin{equation}
\label{ref:IV.5.4-150}
\begin{CD}
\mod(k)\oplus \mod(k) @>>> \mod(k)\\ @. @AAA\\ @. \gr(k[x])
\end{CD}
\end{equation}
where the horizontal map sends $(V_1,V_2)$ to $V_1\oplus V_2$ and
the vertical map is localizing at $x$ and restricting to degree
zero. Thus $\Cscr$ is again uniquely determined.

By the argument given in the beginning of this section the proof
of the main result is now complete.
\begin{remark}
\label{ref:IV.5.1-151} In this remark we indicate some alternative
arguments
  that could have been used in this section.
\begin{enumerate}
\item Consider the case where there  is one simple object modulo
  $\Tscr$.
In that case we could have observed (once we have proved that
  there are exactly two $\tau$-orbits of simples) that the functor $s$
  defined in Section \ref{ref:IV.3-134} defines a ``hypersurface'' in
  $\Cscr$ in the sense of \cite{SmithZhang}. Then it follows
  automatically from the results in \cite{SmithZhang} that $\Cscr$ is determined
  up to equivalence and is given by the category which was denoted by
  $\Cscr$ in \S\ref{ref:III.3.7-106}.  Unfortunately there is
  no corresponding approach for the case where there are two simple
  objects modulo $\Tscr$.
\item Another approach is that instead of constructing the $\ZZ
  A_\infty^\infty$ and $\ZZ D_\infty$ category beforehand we may
  actually show that the categories defined by the pullback diagrams
  (\ref{ref:IV.5.3-149}) and (\ref{ref:IV.5.4-150}) have the expected
  properties (hereditary abelian, existence of almost split sequences,
  etc\dots). This is fairly easy for the diagram
  (\ref{ref:IV.5.3-149})
  since one easily shows that in that case the resulting category is
  equivalent to the one which was denoted by $\Cscr$ in
  \S\ref{ref:III.3.7-106} (invoking the results in
  \cite{SmithZhang}, or directly). If the characteristic is different
  from two then diagram (\ref{ref:IV.5.4-150}) is obtained by a
  skew group construction from diagram (\ref{ref:IV.5.3-149}), so this
  case is easy too by the results in \cite{ReitenRiedtmann}.

  Unfortunately if one wants to include the characteristic 2  case
  there is no alternative but to give a direct proof for diagram
\eqref{ref:IV.5.4-150}. This can indeed be done, but requires some
  work.
\end{enumerate}
\end{remark} We now summarize the main result in this section.
\begin{theorem}\label{ref:IV.5.2-152}
The connected hereditary abelian noetherian $\Ext$-finite
categories with Serre functor and some object of infinite length
and no nonzero projective object are exactly those of one of the
following forms.
\begin{enumerate}
\item The category of coherent sheaves over a hereditary order
over a connected non-singular projective curve.
\item The $\ZZ A_\infty$-category.
\item The $\ZZ D_\infty$-category.
\end{enumerate}
\end{theorem}

Combining Theorem \ref{ref:IV.5.2-152} with Theorem
\ref{ref:III.1.1-65} and Theorem \ref{ref:II.4.9-62}
gives the proof of Theorem \ref{theoremb}.

\chapter{Applications}\label{ref:V-153}
In this chapter we  give some application of our classification theorem
to saturated categories and to graded rings.

An $\Ext$-finite abelian category of finite homological dimension is said to be saturated if
every cohomological functor of finite type $D^b(\Cscr)\r \mod(k)$ is representable \cite{Bondal4}. We will
use our classification result to show that connected saturated noetherian hereditary abelian
categories are either equivalent to $\mod \Lambda$ for  some indecomposable finite dimensional
hereditary $k$-algebra $\Lambda$ or to the category $\coh \Oscr$ for a sheaf of hereditary
orders over a non-singular connected projective curve. This result will also apply to some
categories of the form $\qgr(R)$ since these are often saturated.

\section{Saturatedness}
\label{ref:V.1-154} In this section   we give the application to saturated hereditary abelian
categories.

Let $\Ascr$ be a $\Hom$-finite triangulated category of finite homological dimension.
Following Bondal and Kapranov \cite{Bondal4} we say that $\Ascr$ is \emph{saturated} if every
cohomological functor of finite type $H:\Ascr\r \mod(k)$ is
representable.  Finite type means that for every $A\in\Ascr$ at most a
finite number of $H(A[n])$ are non-zero.

If $\Cscr$ is an $\Ext$-finite
abelian category of finite homological dimension then we say that $\Cscr$ is saturated if the
same holds for $D^b(\Cscr)$.

In \cite{Bondal4} it is shown that the following two categories
are saturated.
\begin{enumerate}
\item $\mod(\Lambda)$ where $\Lambda$ is a finite dimensional algebra
  of finite global dimension.
\item $\coh(X)$ where $X$ is a non-singular projective variety over $k$.
\end{enumerate}
We start with the following easy lemma which is a kind of converse
to 1.
\begin{lemma}
\label{ref:V.1.1-155}
 Assume that $\Cscr$ is an $\Ext$-finite
  abelian category of finite homological dimension in which every
  object has finite length.  Then  $\Cscr$  is
  saturated if and only if $\Cscr\cong \mod(\Lambda)$ for $\Lambda$ a
  finite dimensional algebra of finite global dimension.
\end{lemma}
\begin{proof} That $\mod(\Lambda)$ has the desired properties has
  already been stated, so we prove the converse.

  We embed $\Cscr$ into its closure under direct limits
  $\tilde{\Cscr}$ (see \cite{Gabriel}). Let $E$ be the direct sum of the injective
  hulls of the simples in $\Cscr$. If $S\in\Cscr$ is simple then
  clearly $\Hom(S,E)\cong k$ and hence
\begin{equation}
\label{ref:V.1.1-156} \dim\Hom(M,E)=\length M
\end{equation}
for $M$ in $\Cscr$. For use below we note that it is easy to see
that \eqref{ref:V.1.1-156} holds even for $M\in\tilde{\Cscr}$.

It follows from \eqref{ref:V.1.1-156} that $F=\Hom(-,E)^\ast$ defines a cohomological functor
$D^b(\Cscr)\r \mod(k)$. Since $\Cscr$ is saturated this functor is representable by some
object $P\in D^b(\Cscr)$, and since the restriction of $F$ is zero on $\Cscr[n]$ for $n\neq 0$
and exact on $\Cscr$ it follows that $P$ is actually a projective object in $\Cscr$.

Hence we now have $\Hom(P,-)^\ast\cong\Hom(-,E)$ as functors on $\Cscr$. Applying $\projlim$
it follows that this also holds as functors on $\tilde{\Cscr}$. Thus we have
$\Hom(P,E)\cong\Hom(E,E)^\ast$. Using \eqref{ref:V.1.1-156} we find that $\length E=\length
P<\infty$. Thus $E$ is an injective cogenerator of $\Cscr$ and hence a projective generator of
$\Cscr^{\text{opp}}$. It follows that $\Cscr^{\text{opp}}$ is equivalent to $
\mod(\Lambda^{\text{opp}})$ where $\Lambda=\End(E)$. Since $\Cscr^{\text{opp}}$ has finite
homological dimension, so does $\Lambda$. Furthermore
$\Cscr=\mod(\Lambda^{\text{opp}})^{\text{opp}}=\mod(\Lambda)$ (using the functor $(-)^\ast$).
This finishes the proof.
\end{proof}
We now have the following result for hereditary abelian categories. Note that it is not hard
to see that if $\Cscr$ is saturated, then $\Cscr$ has Serre duality.
\begin{theorem} Assume that $\Cscr$ is a  connected
  noetherian $\Ext$-finite hereditary abelian category. Then
  $\Cscr$ is saturated if and only if it has one of the following forms:
\begin{enumerate}
\item $\mod(\Lambda)$ where $\Lambda$ is a connected finite dimensional
  hereditary algebra.
\item $\coh(\Oscr)$ where $\Oscr$ is a sheaf of hereditary $\Oscr_X$-orders
 (see \ref{ref:III.2-68})   over a non-singular connected
  projective curve $X$.
\end{enumerate}
\end{theorem}
\begin{proof}
We already know that the category in 1. has the required
properties. That this is so for the category in 2. follows in the
same
 way as for $\coh(X)$ where $X$ is a non-singular proper curve.

We now prove the converse. Assume that $\Cscr$ is not of the form
2. This means by Theorem \ref{theoremb} that it  has  one of the
following forms
\begin{itemize}
\item[(a)] $\wrep(Q)$ for $Q$ a star.
\item[(b)] The $\ZZ A_\infty^\infty$ or the $\ZZ D_\infty$ category.
\item[(c)] Finite dimensional nilpotent representations over $\tilde{A}_n$
or $A_\infty^{\infty}$, with all arrows oriented in the same direction.
\end{itemize}
In the categories in (c)  all objects are of finite length but they are
clearly not of the form $\mod(\Lambda)$ since they have no projectives. Hence by lemma \ref{ref:V.1.1-155} it follows that
$\Cscr$ is not of type (c).

By construction the categories in (b)  are derived equivalent to
the finite dimensional representations over an infinite quiver.
Hence again by lemma \ref{ref:V.1.1-155} the categories in (b)
are not saturated.

So we are left with (a). Since $Q$ is a star, there is clearly a section $Q'$ in $\mathbb{Z}Q$
such that all paths in $Q'$ are finite. Then the indecomposable projective and injective
representations of $Q'$ have finite length. By Theorem \ref{ref:II.3.6-52} we know that $\wrep(Q)$ is derived
equivalent to $\rep(Q')$, and all objects in
$\rep(Q')$ have finite length. Again by lemma \ref{ref:V.1.1-155} we know that $\rep(Q')$ is
saturated if and only if $Q'$ is finite. But then the same holds for $\wrep(Q)$. We conclude
that $\Cscr$ is of the form $\wrep(Q)=\rep(Q)$ with $Q$ finite.
\end{proof}
\begin{remark} Abelian categories satisfying the hypotheses of the
  previous theorem (except noetherianness) can, in Bondal's
  terminology \cite{Bondal3}, be viewed as ``non-commutative non-singular
  proper curves''. It appears that non-commutative curves are actually
  quite close to commutative ones. Other manifestations of this
  principle are  \cite{SSW,Staf5,Staf6}.
\end{remark}

\section{Graded rings}
\label{ref:V.2-157} In this section we give the application to the hereditary abelian
categories of the form $\qgr(R)$. We start with recalling the definition of the condition
$\chi$ from \cite{AZ}. Let $R=k+R_1+R_2+\cdots$ be an $\NN$-graded ring which is right
noetherian. Let $m=R_1+R_2+\cdots$ be the graded maximal ideal of $R$. For $M$ in $\Gr(R)$ we
denote by $\Gamma_m(M)$ the maximal torsion submodule of $M$, that is, the maximal submodule
of $M$ which is locally annihilated by $m$. Then $R$ is said to satisfy condition $\chi$ if
for all $M$ in $\gr(R)$ and all $i$ we have that the right derived functor $R^i \Gamma_m (M)$
has right bounded grading.

The following is a special case of a result proved in \cite{BondalVdb}.
\begin{theorem}\cite{BondalVdb}\label{ref:V.2.1-158}
  Let $R=k\oplus R_1\oplus R_2\oplus\cdots$ be a  right noetherian graded ring satisfying
the following hypotheses:
\begin{enumerate}
\item
  $R$ satisfies $\chi$;
\item
$\Gamma_m$  has finite cohomological dimension and the same condition
on the left;
\item
$\Ext^n_R(k,{}_RR)$ is finite dimensional for all $n$ (this is a
special case of the $\chi$ condition for left modules);
\item
There is an $m\in\NN$ such that $\Ext^n_{\qgr(R)}(A,B)=0$ for all
$A,B\in\qgr(R)$ and $n>m$.
\end{enumerate}
Then $\qgr(R)$ is
  saturated.
\end{theorem}
{}{}From this result we obtain the following consequence.
\begin{corollary}
\label{ref:V.2.2-159}
Assume that $R$ satisfies the hypotheses of Theorem
\ref{ref:V.2.1-158} and that $\qgr(R)$ is in addition hereditary.
Then $\qgr(R)$ is a direct sum of categories of the following form
\begin{enumerate}
\item $\mod(\Lambda)$ where $\Lambda$ is a connected finite dimensional
  hereditary algebra.
\item $\coh(\Oscr)$ where $\Oscr$ is a sheaf of hereditary $\Oscr_X$-orders
 (see \S\ref{ref:III.2-68})   over a non-singular connected
  projective curve $X$.
\end{enumerate}
\end{corollary}

\appendix
\chapter{Some results on abelian categories}
\setcounter{section}{1}
\renewcommand{\thesection}{\Alph{section}}
\setcounter{lemma}{0} \setcounter{equation}{0} Here we give a few
homological results on abelian categories that have been used
before. Since they have a somewhat different flavor from the other
results in this paper we have decided to collect them here in an
appendix.

The first result   (Proposition \ref{ref:2.2a} below) implies in particular that if we take
all finite length objects in a hereditary abelian category we get
again a hereditary abelian category.  The
present proof of Proposition \ref{ref:2.2a} is based on the following lemma, where the proof is
included for the convenience of the reader.
\begin{lemma} Assume that $\Hscr$ is an abelian category. Then $\Hscr$
  is hereditary abelian if and only if $\Ext^1_\Hscr(X,-)$ preserves epis for
  all $X\in \Hscr$.
\end{lemma}
\begin{proof} Let us prove the non-obvious implication. Thus assume that
  $\Ext^1_\Hscr(X,-)$ preserves epis. Let
\begin{equation}
\label{ref:1a} 0\r A\xrightarrow{u} B\xrightarrow{v}
C\xrightarrow{w} D\r 0
\end{equation}
be an exact sequence. We have to show that it is Yoneda equivalent
to a split one.

Let $K=\ker w$. Then we have a surjective map $B\r K$, as well as
an exact sequence $0\r K\r C\r D\r 0$, representing an element of
$\Ext^1_\Hscr(D,K)$. By the fact that $\Ext^1_\Hscr(D,-)$
preserves epis, this can be lifted to an element of
$\Ext^1_\Hscr(D,B)$. Thus we obtain a commutative diagram with
exact rows:
\[
\begin{CD}
  @. 0   @>>> B @>i>> E @>>> D @>>> 0\\
@.   @.       @|     @VVV     @|\\ 0 @>>> A @>u>> B @>v>> C @>w>>
D @>>> 0
\end{CD}
\]
which can be transformed into the following diagram.
\[
\begin{CD}
0  @>>> A   @>(1_A,0)>> A\oplus B @>(0,i)>> E @>>> D @>>> 0\\ @.
@|       @V(u,1_B)VV     @VVV     @|\\ 0 @>>> A @>u>> B @>v>> C
@>w>> D @>>> 0
\end{CD}
\]
from which it follows that \eqref{ref:1a} represents a trivial
class in $\Ext^2_\Hscr(D,A)$.
\end{proof}
\begin{proposition}
\label{ref:2.2a}
  Let \h\ be a hereditary abelian category. Then any Serre subcategory
  $\h'$ of $\h$ is hereditary.
\end{proposition}
\begin{proof}
This now follows from the previous lemma using the fact that for
$X,Y\in\Hscr'$ one has
$\Ext^1_{\Hscr'}(X,Y)\cong\Ext^1_{\Hscr}(X,Y)$ (since $\Hscr'$ is
closed under extensions).
\end{proof}

If one has a noetherian abelian category $\Hscr$ of finite
homological dimension then one sometimes needs to go to the
closure $\tilde{\Hscr}$ of $\Hscr$ under direct limits. It follows
from the following proposition that this operation preserves
homological dimension.
\begin{proposition}
\label{ref:2.3a}
   Let \X\ be a noetherian abelian category where the
   isomorphism classes of objects form a set, and assume
   $\X \subset \tilde{\X}$ where $\tilde{\X}$ is an abelian
   category where $\tilde{\X}$ is the completion of \X\ with
   respect to direct limits of objects in \X (see \cite{Gabriel}). Then $\gd \X=\gd\tilde{\X}$.
\end{proposition}
\begin{proof}
(i) We first show $\gd\X \leq \gd\tilde{\X}$. We have that \X\  is
a thick (Serre) abelian subcategory of $\tilde{\X}$. Let $f \colon
B \rightarrow C$ be an epimorphism in $\tilde{X}$, with $C$ in \X.
We have by assumption $B=\dlim_{i}B_i$, where each $B_i$ is in \X.
Since $C$ is noetherian, there is some $j$ such that the induced
map $B_j \rightarrow C$ is an epimorphism. By the dual of
\cite[1.7.11]{KSI}  we have an isomorphism of derived categories
$D^b(\X)\simeq D_{\X}^b(\tilde{\X})$, and in particular
$\Ext_\X^t(X,Y)\simeq \Ext_{\tilde{\X}}^t(X,Y)$ for $X$ and $Y$ in
\X.  It follows that $\gd \X \leq \gd \tilde{\X}$.

(ii) Assume now that $\gd \X=r$. We want to show that for $M$ and
$N$  in $\tilde{\X}$ we have $\Ext_{\tilde{\X}}^{r+1}(M,N)=0$. If
$M$ and $N$ are both in \X, this follows from the above, since
$\Ext_{\tilde{\X}}^{r+1}(M,N)\simeq \Ext_\X^{r+1}(M,N)$, which is
0 by assumption. If $M$ is in \X, and hence by assumption is
noetherian, $\Hom_{\tilde{\X}}(M,\;)$ commutes with $\dlim_{ i}$,
and hence the same is true for the derived functors \cite{Groth1}.
It then follows that $\Ext_\X^{r+1}(M,N)=0$ if $M$ is in \X.

Let now $N$ be in $\tilde{\X}$, and consider  the exact sequence
\[0 \ra N \ra E_0 \ra E_1 \ra \cdots \ra E_{r-1}
\stackrel{f_{r-1}}{\longrightarrow} E_r
\stackrel{f_r}{\longrightarrow} Z_{r+1} \ra 0,\] where the $E_i$
are injective. Then we have $\Ext_{\tilde{\X}}^{r+1}(M,N)\simeq
\cok(\Hom_{\tilde{\X}}(M,E_r) \ra \Hom_{\tilde{\X}}(M,Z_{r+1}))$
for all $M$ in $\tilde{\X}$. In particular, for any inclusion map
$j \colon M \ra Z_{r+1}$ with $M$ in \X\ there is a commutative
diagram
\[ \xymatrix{& & & M \ar[dl]_{g_r}\ar[d]^{j} & \\ 0\ar[r] & Z_r \ar[r]^{f_{r-1}} & E_r \ar[r]^{f_r} &
Z_{r+1} \ar[r] & 0,}\] where $Z_r=\im f_{r-1}$. If we have  a
commutative diagram
\[ \xymatrix{M' \ar[dr]_{u} \ar[r]^{j'} & Z_{r+1} \\
    & M \ar[u]_{j},}\]
with $M$ and $M'$ in $\tilde{\X}$ and $u$ a monomorphism, we then have $g_r' \colon M' \ra
E_r$ with $f_rg_r'=j'$ and $g_r \colon M \ra E_r$ with $f_rg_r=j$, but we do not necessarily
have $g_ru=g_r'$. To adjust the map $g_r$ observe that $f_r (g_ru-g_r')=0$. Hence there is a
map $t \colon M'\ra Z_r$ with $f_{r-1}t=g_ru-g_r'$. Since $\Ext^1_{\tilde{\X}}(M/M',Z_r)
\simeq \Ext^{r+1}_{\tilde{\X}}(M/M',N)=0$ because $M/M'$ is in \X, there is a map  $s \colon M
\ra Z_r$ such that $su=t$, and hence $f_{r-1}su=g_ru-g_r'$. Letting $\bar{g_r}=-g_r+f_{r-1}s$
we see that $\bar{g_r}u=g_r'$ as desired.

We now use Zorn's lemma to choose a maximal pair $(j \colon T
\hookrightarrow Z_{r+1}, \phi \colon T \rightarrow E_r \text{ with
} f_r\phi=j)$. Assume that $T \neq Z_{r+1}$. Since $Z_{r+1}$ is
the direct limit of objects in \X, there is some $M \subset
Z_{r+1}$ with $M$ in \X\ and $M \not\subset T$. Consider $M \cap T
\subset T \subset Z_{r+1}$ and the map $\phi|_{M \cap T} \colon M
\cap T \ra E_r$. Then there is by the first part of the proof some
map $\phi' \colon M \ra E_r$ such that $f_{r} \phi' = i$ with $i
\colon M \ra Z_{r+1}$, and $\phi'|_{M\cap T}=\phi|_{M\cap T}$.
Hence the map $\phi'-\phi:M\oplus T\r E_r$ factors through a map
$\tilde{\phi} \colon M+T \ra E_r$. This map has the property that
$f_r\tilde{\phi}=i'$ with $i' \colon T+M \ra Z_{r+1}$. This
contradicts the maximality of $(j \colon T \hookrightarrow
Z_{r+1},\phi \colon T \ra E_r)$, and hence $T=Z_{r+1}$. This means
that $f_r \colon E_r \ra Z_{r+1}$ is a split epimorphism, and
hence $\id_{\tilde{\X}} N \leq r$, so that
$\Ext^{r+1}_{\tilde{\X}}(M,N)=0$ for all $M$ in $\tilde{\X}$.
Hence we have $\gd \tilde{\X} \leq r = \gd \X$.
\end{proof}
Now we consider quotient categories. Recall that a
\emph{Grothendieck category} is an abelian category with a
generator and exact direct limits. Such categories automatically
have products and injective hulls \cite{stenstrom}. Prominent
examples of Grothendieck categories are locally noetherian
categories. These are the abelian categories which are the closure
of  noetherian categories under direct limits (that is: categories
of the form $\tilde{\mathcal{X}}$ as considered in the previous
proposition). A localizing subcategory $\Cscr$ of a Grothendieck
category $\Ascr$ is a Serre subcategory which is closed under
direct limits. In that case $\Ascr/\Cscr$ is also a Grothendieck
category. $\Cscr$ is called \emph{stable} if it is closed under
injective hulls in $\Ascr$, or equivalently, if it is closed under
essential extensions. The following result is proved in
\cite{stenstrom} (and is easy to see).
\begin{proposition}
\label{ref:2.4a} Let $\Ascr$ be a Grothendieck category and
$\Cscr$ a
  stable localizing subcategory. Then the quotient functor $T:\Ascr\r
  \Ascr/\Cscr$ preserves injective hulls and in particular $\gldim
  \Ascr/\Cscr\le \gldim \Ascr$.
\end{proposition}

\backmatter

%\bibliographystyle{amsabbrv}
%\bibliography{mybibs}

\begin{thebibliography}{10}

\bibitem{Staf5}
M.~Artin and J.~T. Stafford, {\em Noncommutative graded domains with quadratic
  growth}, Invent. Math. {\bf 122} (1995), 231--276.

\bibitem{Staf6}
M.~Artin and J.~T. Stafford, {\em Semiprime graded algebras of dimension two},
  J. Algebra {\bf 227} (2000), no.~1, 68--123.


\bibitem{AZ}
M.~Artin and J.~J. Zhang, {\em Noncommutative projective schemes}, Adv. in
  Math. {\bf 109} (1994), no.~2, 228--287.

\bibitem{MAIR}
M.~Auslander and I.~Reiten, {\em Almost split sequences in dimension two}, Adv.
  in Math. {\bf 66} (1987), no.~1, 88--118.

\bibitem{ARS}
M.~Auslander, I.~Reiten, and S.~Smalo, {\em Representation theory of {A}rtin
  algebras}, Cambridge Studies in Advanced Mathematics, vol.~36, Cambridge
  University Press, 1995.

\bibitem{Beilinson1}
A.~Beilinson, {\em On the derived category of perverse sheaves}, K-theory,
  Arithmetic and Geometry, Lecture Notes in Mathematics, vol. 1289, Springer
  Verlag, 1987, pp.~27--41.

\bibitem{BBD}
A.~Beilinson, J.~Bernstein, and P.~Deligne, {\em Faisceaux pervers},
  Ast{\'e}risque, vol. 100, Soc. Math. France, 1983.

\bibitem{Bondal4}
A.~I. Bondal and M.~M. Kapranov, {\em Representable functors, {S}erre functors,
  and reconstructions}, Izv. Akad. Nauk SSSR Ser. Mat. {\bf 53} (1989), no.~6,
  1183--1205, 1337.

\bibitem{BondalVdb}
A.~I. Bondal and M.~Van~den Bergh, Generators of triangulated
categories and representability of functors, in preparation.

\bibitem{Bondal3}
A.~I. Bondal, {\em Non-commutative deformations and {P}oisson brackets on
  projective spaces}, MPI-preprint, 1993.

\bibitem{Gabriel}
P.~Gabriel, {\em Des cat\'egories ab\'eliennes}, Bull. Soc. Math. France {\bf
  90} (1962), 323--448.


\bibitem{GeigLen}
W.~Geigle and H. Lenzing, {\em A class of weighted projective curves arising in representation
theory of finite dimensional algebras, in: Singularities, representations of algebras and
vector bundles}, Springer Lecture Notes {\bf 1274} (1987), 265--297.


\bibitem{Groth1}
A.~Grothendieck, {\em Sur quelques points d'alg\`ebre homologiques}, T{\^o}hoku
  Math. J. (2) {\bf 9} (1957), 119--221.

\bibitem{Happel}
D.~Happel, {\em Triangulated categories in the representation theory of finite
  dimensional algebras}, London Mathematical Society lecture note series, vol.
  119, Cambridge University Press, 1988.

\bibitem{HappelReiten}
D.~Happel, {\em A characterization of hereditary categories with tilting object}, Preprint
2000.



\bibitem{HRS}
D.~Happel, I.~Reiten, and S.~Smalo, {\em Tilting in abelian categories and
  quasitilted algebra}, Memoirs of the AMS, vol. 575, Amer. Math. Soc., 1996.

\bibitem{H}
R.~Hartshorne, {\em Algebraic geometry}, Springer-Verlag, 1977.

\bibitem{KSI}
M.~Kashiwara and P.~Schapiro, {\em Sheaves on manifolds}, Die Grundlehren der
  mathematischen Wissenschaften, vol. 292, Springer Verlag, 1994.

\bibitem{Lenzing}
H.~Lenzing, {\em Hereditary noetherian categories with a tilting object}, Proc.
  Amer. Math. Soc. {\bf 125} (1997), 1893--1901.

\bibitem{MR}
J.~C. McConnell and J.~C. Robson, {\em Noncommutative {N}oetherian rings}, John
  Wiley {\&} Sons, New York, 1987.

\bibitem{reiner}
I.~Reiner, {\em Maximal orders}, Academic Press, New York, 1975.

\bibitem{ReitenRiedtmann}
I.~Reiten and C.~Riedtmann, {\em Skew group algebras in the representation
  theory of artin algbras}, J. Algebra {\bf 92} (1985), 224--282.

\bibitem{IRMB}
I.~Reiten and M.~Van~den Bergh, {\em Grothendieck groups and tilting objects},
  Algebras and Representation Theory, to appear.

\bibitem{Ri1}
C.~M. Ringel, {\em A ray quiver construction of hereditary abelian categories with Serre
duality}, Proc. ICRA IX, Beijing.

\bibitem{Ri2}
C.~M. Ringel, {\em The diamond category of a locally discrete ordered set}, Proc. ICRA IX,
Beijing.

\bibitem{RS1}
C.~Robson and L.~Small, {\em Hereditary prime {P.I.} rings are classical
  hereditary orders.}, J. London Math. Soc. (2) {\bf 8} (1974), 499--503.

\bibitem{ASch}
A.~Schofield, unpublished.

\bibitem{SSW}
L.~W. Small, J.~T. Stafford, and R.~B. Warfield, {\em Affine algebras of
  {G}elfand-{K}irillov dimension one are {PI}}, Math. Proc. Camb. Phil. Soc.
  {\bf 97} (1985), 407--414.

\bibitem{SOS}
S.~O. Smaloe, {\em Almost split sequences in categories of representations of
  quivers}, Preprint~3, Trondheim, 1999.

\bibitem{SmithZhang}
S.~P. Smith and J.~J. Zhang, {\em Curves on quasi-schemes}, Algebr. Represent.
  Theory {\bf 1} (1998), no.~4, 311--351.

\bibitem{stenstrom}
B.~Stenstr{\"o}m, {\em Rings of quotients}, Die {G}rundlehren der
  mathematischen {W}issenschaften in {E}inzeldarstellungen, vol. 217, Springer
  Verlag, Berlin, 1975.

\bibitem{VdB19}
M.~Van~den Bergh, {\em Blowing up of non-commutative smooth surfaces}, to
  appear in Memoirs of the AMS; math.QA/980911, 1998.

\bibitem{VdBVG}
M.~Van~den Bergh and M.~Van~Gastel, {\em Graded modules of {G}elfand-{K}irillov
  dimension one over three-dimensional {A}rtin-{S}chelter regular algebras}, J.
  Algebra {\bf 196} (1997), 251--282.

\bibitem{Verdier}
J.-L. Verdier, {\em Des cat\'egories d\'eriv\'ees des cat\'egories
  ab\'eliennes}, Ast\'erisque (1996), no.~239, xii+253 pp. (1997), With a
  preface by Luc Illusie, Edited and with a note by Georges Maltsiniotis.

\end{thebibliography}
\ifx\undefined\bysame
\newcommand{\bysame}{\leavevmode\hbox to3em{\hrulefill}\,}
\fi

\end{document}